\documentclass[12pt]{amsart}
\pdfoutput=1

\usepackage[utf8]{inputenc}
\usepackage[T1]{fontenc}
\usepackage{latexsym}
\usepackage{amsxtra}
\usepackage{amsmath}
\usepackage{amsfonts}
\usepackage{amssymb}
\usepackage{blkarray}
\usepackage{multirow}
\usepackage{pgf,tikz}
\usepackage{mathrsfs}
\usepackage{mathdots}
\usepackage{pdfpages}
\usepackage[colorlinks,linkcolor=blue,anchorcolor=black,citecolor=blue,filecolor=blue,menucolor=blue,runcolor=blue,urlcolor=blue]{hyperref}
\usetikzlibrary{arrows}

\tikzstyle{vertex}=[circle, draw, inner sep=0pt, minimum size=4pt]
\newcommand{\vertex}{\node[vertex]}

\newtheorem{theorem}{Theorem}[section]
\newtheorem{lemma}[theorem]{Lemma}
\newtheorem{corollary}[theorem]{Corollary}
\newtheorem{proposition}[theorem]{Proposition}

\newtheorem{hypothesis}[theorem]{Hypothesis}
\newtheorem{claim}[theorem]{Claim}
\newenvironment{subproof}[1][\proofname]{%
  \begin{proof}[#1]%
}{%
  \end{proof}%
}
\numberwithin{equation}{section}

\theoremstyle{definition}
\newtheorem{definition}[theorem]{Definition}

\theoremstyle{remark}

\begin{document}
\title{The highly connected even-cycle and even-cut matroids}

\author{Kevin Grace}
\address{Mathematics Department\\
Louisiana State University\\
Baton Rouge, Louisiana}
\email{kgrace3@lsu.edu}

\author{Stefan H. M. van Zwam}
\address{Mathematics Department\\
Louisiana State University\\
Baton Rouge, Louisiana}
\email{svanzwam@math.lsu.edu}
\thanks{The research for this paper was supported by National Science Foundation grant 1500343. The first author was also supported by a Huel D. Perkins Fellowship from the Louisiana State University Graduate School.}

\date{\today}

\begin{abstract}
The classes of even-cycle matroids, even-cycle matroids with a blocking pair, and even-cut matroids each have hundreds of excluded minors. We show that the number of excluded minors for these classes can be drastically reduced if we consider in each class only the highly connected matroids of sufficient size.
\end{abstract}

\maketitle
\section{Introduction}
\label{introduction}
The complete lists of excluded minors for the classes of even-cycle matroids and even-cut matroids are currently unknown. Irene Pivotto and Gordon Royle \cite{pr16} have found nearly 400 different excluded minors for the class of even-cycle matroids. We will show that, subject to a certain hypothesis described below, a highly connected binary matroid $M$ of sufficient size is an even-cycle matroid if and only if it contains no minor isomorphic to one of three matroids. Similarly, subject to that same hypothesis, a highly connected binary matroid $M$ of sufficient size is an even-cut matroid if and only if it contains no minor isomorphic to one of two matroids.

Unexplained notation and terminology in this paper will generally follow Oxley \cite{o11}. One exception is that we denote the vector matroid of a matrix $A$ by $M(A)$, rather than $M[A]$.

An even-cycle matroid is a binary matroid of the form $M=M\binom{w}{D}$, where $D\in\mathrm{GF}(2)^{V\times E}$ is the vertex-edge incidence matrix of a graph $G=(V,E)$ and $w\in \mathrm{GF}(2)^E$ is the characteristic vector of a set $W\subseteq E$. The pair $(G,W)$ is an \textit{even-cycle representation} of $M$. The edges in $W$ are called \textit{odd} edges, and the other edges are \textit{even} edges. \textit{Resigning} at a vertex $u$ of $G$ occurs when all the edges incident with $u$ are changed from even to odd and vice-versa. It is easy to see that this corresponds to adding the row of the matrix corresponding to $u$ to the characteristic vector of $W$. Therefore, resigning at a vertex does not change an even-cycle matroid. A pair of vertices $u,v$ of $G$ is a \textit{blocking pair} of $(G,W)$ if $(G,W)$ can be resigned so that every odd edge is incident with $u$ or $v$. We will say that an even-cycle matroid has a blocking pair if it has an even-cycle representation with a blocking pair.

In her PhD thesis \cite{p11}, Pivotto gives several descriptions of even-cut matroids, each of which can serve as a definition. The most practical definition for our purposes follows. An \textit{even-cut matroid} is a matroid $M$ that can be represented by a binary matrix with a row whose removal results in a matrix representing a cographic matroid. One can also think of an even-cut matroid as arising from a \emph{graft}, which is a pair $(G,T)$, where $G$ is a graph and $T$ is a subset of $V(G)$ of even cardinality whose members are called \emph{terminals}. The collection of inclusion-wise minimal edge cuts $\delta(U)$, where $U\subseteq V(G)$ and $|U\cap T|$ is even is the collection of circuits for an even-cut matroid. The graft $(G,T)$ is an \emph{even-cut representation} of that matroid. 

This paper is a continuation of \cite{gvz17}, which is based on the work of Geelen, Gerards, and Whittle in \cite{ggw15}. The results announced in \cite{ggw15} rely on the Matroid Structure Theorem by these same authors \cite{ggw06}. Our results are based on hypotheses given by Geelen, Gerards, and Whittle in \cite{ggw15}, as modified in \cite{gvz-prob}. These hypotheses are believed to be true, but their proofs are still forthcoming in future papers by Geelen, Gerards, and Whittle. One of these hypotheses, when restricted to the binary case, is Hypothesis \ref{hyp:connected-template}. We delay the precise statement of Hypothesis \ref{hyp:connected-template} to Section \ref{Binary Frame Templates} due to its technical nature.

Before we can state our main results, we need some additional definitions. If $\mathcal{F}$ is a collection of matroids, let $\mathcal{EX}(\mathcal{F})$ denote the class of binary matroids with no minor contained in $\mathcal{F}$. We denote by $\mathrm{PG}(3,2)\backslash e$, or $\mathrm{PG}(3,2)_{-2}$, or $\mathrm{PG}(3,2)\backslash L$, respectively, the matroid obtained by deleting from $\mathrm{PG}(3,2)$ one element, or two elements, or the three points of a line. Note that $\mathrm{PG}(3,2)\backslash L$ is the vector matroid of the following matrix and, therefore,
is the even-cycle matroid represented by the graph in Figure ~\ref{fig:PG32L}, \begin{figure}
 \begin{center}

\[\begin{tikzpicture}[x=1.3cm, y=1cm,
    every edge/.style={
        draw,
        }
]
	\vertex[fill] (1) at (2,3) {};
	\vertex[fill] (2) at (1,1) {};
	\vertex[fill] (3) at (3,1) {};
	\vertex[fill] (4) at (2,1.75) {};
	\path
		(1) edge (2)
		(1) edge[bend right=10, line width=1.pt] (2)
		(2) edge (3)
		(2) edge[bend right=10, line width=1.pt] (3)
		(3) edge[bend left=10] (4)
		(3) edge[bend right=10, line width=1.pt] (4)
		(1) edge[bend left=10] (4)
		(1) edge[bend right=10, line width=1.pt] (4)
		(2) edge[bend left=10] (4)
		(2) edge[bend right=10, line width=1.pt] (4)
		(1) edge (3)
		(1) edge[bend left=10, line width=1.pt] (3)
	;
\end{tikzpicture}\]
\end{center}
\caption{Even-cycle representation of $\mathrm{PG}(3,2)\backslash L$}
  \label{fig:PG32L}
\end{figure} with odd edges printed in bold.
\[
\left[
\begin{array}{cccccccccccc}
0&0&0&0&0&0&1&1&1&1&1&1\\
1&0&0&1&1&0&1&0&0&1&1&0\\
0&1&0&1&0&1&0&1&0&1&0&1\\
0&0&1&0&1&1&0&0&1&0&1&1\\
\end{array}
\right].
\]
We define $L_{19}$ to be the dual of the cycle matroid of the graph obtained from $K_7$ by deleting two adjacent edges, and we define $L_{11}$ to be the vector matroid of the following matrix.
\[
\left[
\begin{array}{ccccccccccc}
1&0&0&0&0&0&1&0&1&0&1\\
0&1&0&0&0&0&1&0&0&1&1\\
0&0&1&0&0&0&0&1&1&1&0\\
0&0&0&1&0&0&0&1&1&0&1\\
0&0&0&0&1&0&0&1&0&1&1\\
0&0&0&0&0&1&0&0&1&1&1\\
\end{array}
\right].
\]
Finally, let $H_{12}$ be the matroid with the even-cycle representation given in Figure ~\ref{fig:H12}\begin{figure}
 \begin{center}
\[\begin{tikzpicture}[x=1.3cm, y=1cm,
    every edge/.style={
        draw,
        }
]
	\vertex[fill] (1) at (2,3) {};
	\vertex[fill] (2) at (1,2) {};
	\vertex[fill] (3) at (1,1) {};
	\vertex[fill] (4) at (3,1) {};
	\vertex[fill] (5) at (3,2) {};
	\path
		(1) edge[bend left=10] (2)
		(1) edge[line width=1.pt] (2)
		(2) edge[bend left=20] (3)
		(2) edge[line width=1.pt] (3)
		(3) edge[bend left=20] (4)
		(3) edge[line width=1.pt] (4)
		(4) edge[bend left=20] (5)
		(4) edge[line width=1.pt] (5)
		(1) edge[bend right=10] (5)
		(1) edge[line width=1.pt] (5)
		(2) edge[bend left=20] (5)
		(2) edge[line width=1.pt] (5)
	;
\end{tikzpicture}\]
\end{center}
\caption{Even-cycle representation of $H_{12}$}
  \label{fig:H12}
\end{figure}
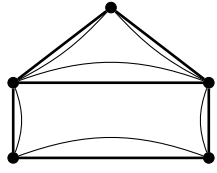. Again, odd edges are printed in bold.

We will prove the following theorems in Section ~\ref{Even-Cycle Matroids}.
\begin{theorem}
 \label{connectedevencycle}
Suppose Hypothesis \ref{hyp:connected-template} holds. Then there exists $k\in\mathbb{Z}_+$ such that a $k$-connected binary matroid with at least $2k$ elements is contained in $\mathcal{EX}(\mathrm{PG}(3,2)\backslash e, L_{19}, L_{11})$  if and only if it is an even-cycle matroid.
\end{theorem}
\begin{theorem}
 \label{blockingpair}
Suppose Hypothesis \ref{hyp:connected-template} holds. Then there exists $k\in\mathbb{Z}_+$ such that a $k$-connected binary matroid with at least $2k$ elements is contained in $\mathcal{EX}(\mathrm{PG}(3,2)\backslash L, M^*(K_6))$ if and only if it is an even-cycle matroid with a blocking pair.
\end{theorem}

We will prove the following theorem in Section \ref{Even-Cut Matroids}.
\begin{theorem}
 \label{connectedevencut}
Suppose Hypothesis \ref{hyp:connected-template} holds. Then there exists $k\in\mathbb{Z}_+$ such that  a $k$-connected binary matroid with at least $2k$ elements is contained in $\mathcal{EX}(M(K_6), H^*_{12})$ if and only if it is an even-cut matroid.
\end{theorem}

Pivotto \cite[Section 2.4.2]{p11} showed that the class of even-cycle matroids with a blocking pair consists of the duals of the members of the class of matroids with an even-cut representation with at most four terminals. Moreover, it is well-known that, for every positive integer $k$, a matroid $M$ is $k$-connected if and only if $M^*$ is $k$-connected. Therefore, Theorem \ref{blockingpair} immediately implies the following result.
\begin{corollary}
 \label{cor:four-terminals}
Suppose Hypothesis \ref{hyp:connected-template} holds. Then there exists $k\in\mathbb{Z}_+$ such that a $k$-connected binary matroid with at least $2k$ elements is contained in $\mathcal{EX}((\mathrm{PG}(3,2)\backslash L)^*, M(K_6))$ if and only if it has an even-cut representation with at most four terminals.
\end{corollary}

A \emph{rank-$(\leq t)$ perturbation} of a binary matroid $M$ is the vector matroid of the matrix obtained by adding a binary matrix of rank at most $t$ to a binary matrix representing $M$. The work in this paper, as well as the work in \cite{ggw15}, \cite{gvz17}, and \cite{gvz-prob}, is based on the hypothesis that the highly connected members of a minor-closed class of binary matroids are rank-$(\leq t)$ perturbations of graphic or cographic matroids. Geelen, Gerards, and Whittle \cite{ggw15} introduced the notion of a \emph{template} to describe the perturbations in more detail. Hypothesis \ref{hyp:connected-template} states that the sufficiently highly connected members of a proper minor-closed class of binary matroids can be constructed using finitely many of these templates. Moreover, every matroid constructed using one of these templates must be contained in the minor-closed class.

In Sections \ref{Binary Frame Templates} and \ref{sec:template minors}, we will review the necessary definitions and results found in \cite{ggw15}, \cite{gvz17}, and \cite{gvz-prob}. In Section \ref{Binary Frame Templates}, we give the definition of a template and state Hypothesis \ref{hyp:connected-template}. In Section \ref{sec:template minors}, we define a preorder $\preceq$ on the set of templates and list several reduction operations that produce from a template $\Phi$ another template $\Phi'$ such that $\Phi'\preceq\Phi$. If $\Phi'\preceq\Phi$, then every matroid constructed using $\Phi'$ is a minor of some matroid constructed using $\Phi$.

In Section \ref{sec:Excluded Minors}, we prove that $\mathrm{PG}(3,2)\backslash e$, $L_{19}$, $L_{11}$, $\mathrm{PG}(3,2)\backslash L$, $M(K_6)$, and $H^*_{12}$ are indeed excluded minors for the respective classes given in the theorems above.

The rest of the paper is devoted to the converse statements. Much of the work required to prove these statements involves analysis of specific templates, showing that either one of these excluded minors can be constructed using the template or that the template is highly structured -- to the point that only even-cycle matroids (or even-cycle matroids with a blocking pair, or even-cut matroids) can be constructed using the template. The finite case checks involved in this process are by and large carried out using the SageMath software system ~\cite{sage}. The code for the computations can be found in the Appendix, and the technical lemmas proved by the computations will be given in Sections \ref{sec:technical-cycle} and \ref{sec:technical-cut}. In Section \ref{Even-Cycle Matroids}, we prove Theorems \ref{connectedevencycle} and \ref{blockingpair}, and in Section \ref{Even-Cut Matroids}, we prove Theorem \ref{connectedevencut}. Finally, in Section \ref{sec:vertical}, we prove variations of Theorems \ref{connectedevencycle}-\ref{connectedevencut} with weaker notions of connectivity.

Our main results and the techniques used to prove them give no indication of how large the value for $k$ must be. The sets of matroids in our theorems are not unique, and their members do not necessarily need to be excluded minors for the classes we study. For example, $L_{19}$ and $M^*(K_6)$ can be replaced with $M^*(K_n)$ for $n>6$, and $M(K_6)$ can be replaced with $M(K_n)$ for $n>6$. We chose the small matroids that we did because they are actually excluded minors for the various classes. Presumably, this comes at the cost of a larger value for $k$.

The following notation will be used throughout this paper. We denote an empty matrix by $[\emptyset]$. We denote a group of one element by $\{0\}$ or $\{1\}$, if it is an additive or multiplicative group, respectively. If $S'$ is a subset of a set $S$ and $G$ is a subgroup of the additive group $\mathbb{F}^S$, we denote by $G|S'$ the projection of $G$ into $\mathbb{F}^{S'}$. Similarly, if $\bar{x}\in G$, we denote the projection of $\bar{x}$ into $\mathbb{F}^{S'}$ by $\bar{x}|S'$.

\section{Binary Frame Templates}
\label{Binary Frame Templates}
In this section, we will define the notion of a frame template. We will simplify the definition found in \cite{ggw15} by restricting ourselves to the binary case.

A \textit{binary frame matrix} is a binary matrix in which every column has at most two nonzero entries. The vector matroid of any such matrix is a graphic matroid.

The structure theorem given by Geelen, Gerards, and Whittle \cite{ggw15} is somewhat technical. They introduced the concept of a template in order to facilitate the description of the theorem. A \textit{binary frame template} is a tuple $\Phi=(C,X,Y_0,Y_1,A_1,\Delta,\Lambda)$ such that the following hold:
\begin{itemize}
 \item [(i)] $C$, $X$, $Y_0$, and $Y_1$ are disjoint finite sets.
 \item [(ii)] $A_1$ is binary matrix with rows indexed by $X$ and columns indexed by $Y_0\cup Y_1\cup C$.
 \item [(iii)] $\Lambda$ is a subgroup of the additive group of $\mathrm{GF}(2)^X$.
 \item [(iv)] $\Delta$ is a subgroup of the additive group of $\mathrm{GF}(2)^{Y_0 \cup Y_1\cup C}$.
\end{itemize}
A binary frame template is in \emph{standard form} if there are partitions $C=C_0\cup C_1$ and $X=X_0\cup X_1$ such that $A_1[X_0,C_0]$ is an identity matrix and such that $A_1[X_1,C]$ is a zero matrix. Note that $|C_0|=|X_0|$ and that some or all of $C_0,C_1,X_0,$ and $X_1$ may be empty. It was proved in \cite[Lemma 3.11]{gvz17} that every binary frame template can be put into standard form; therefore, for the rest of the paper, we will assume that all binary frame templates are in standard form. Furthermore, the term \emph{template} will always refer to a binary frame template in standard form, unless otherwise specified.

Let $\Phi=(C,X,Y_0,Y_1,A_1,\Delta,\Lambda)$ be a template. Let $B$ and $E$ be disjoint finite sets, and let $A'$ be a binary matrix with rows indexed by $B$ and columns indexed by $E$. We say that $A'$ \textit{respects} $\Phi$ if the following hold:
\begin{itemize}
 \item [(i)] $X\subseteq B$ and $C, Y_0, Y_1\subseteq E$.
 \item [(ii)] $A'[X, C\cup Y_0\cup Y_1]=A_1$.
 \item [(iii)] There exists a set $Z\subseteq E-(C\cup Y_0\cup Y_1)$ such that $A'[X,Z]=0$, each column of $A'[B-X,Z]$ is a unit vector, and $A'[B-X, E-(C\cup Y_0\cup Y_1\cup Z)]$ is a binary frame matrix.
 \item [(iv)] Each column of $A'[X,E-(C\cup Y_0\cup Y_1\cup Z)]$ is contained in $\Lambda$.
 \item [(v)] Each row of $A'[B-X, C\cup Y_0\cup Y_1]$ is contained in $\Delta$.
\end{itemize}

Figure \ref{fig:A'} \begin{figure}
\begin{center}
\begin{tabular}{ r|c|c|cccc| }
\multicolumn{2}{c}{}&\multicolumn{1}{c}{$Z$}&\multicolumn{1}{c}{$Y_0$}&\multicolumn{1}{c}{$Y_1$}&\multicolumn{1}{c}{$C_0$}&\multicolumn{1}{c}{$C_1$}\\
\cline{2-7}
$X_0$&columns from $\Lambda|X_0$&0&\multicolumn{2}{c|}{\multirow{2}{*}{$*$}}&\multicolumn{1}{c|}{$I$}&$*$\\
\cline{2-3} \cline{6-7}
$X_1$&columns from $\Lambda|X_1$&0&&&\multicolumn{2}{|c|}{0}\\
\cline{2-7}
&\multirow{5}{*}{frame matrix}&\multirow{5}{*}{unit or zero columns}&\multicolumn{4}{c|}{\multirow{5}{4em}{rows from  $\Delta$}}\\
&&&&&&\\
&&&&&&\\
&&&&&&\\
&&&&&&\\
\cline{2-7}
\end{tabular}
\end{center}
\caption{The matrix $A'$}
  \label{fig:A'}
\end{figure} shows the structure of $A'$.

Suppose that $A'$ respects $\Phi$ and that $Z$ satisfies (iii) above. Now suppose that $A\in \mathrm{GF}(2)^{B\times E}$ satisfies the following conditions:
\begin{itemize}
\item [(i)] $A[B,E-Z]=A'[B,E-Z]$
\item [(ii)] For each $i\in Z$ there exists $j\in Y_1$ such that the $i$-th column of $A$ is the sum of the $i$-th and the $j$-th columns of $A'$.
\end{itemize}
We say that $A$ \textit{conforms} to $\Phi$. We say that a matroid $M$ \textit{conforms} to $\Phi$ if there is a matrix $A$ that conforms to $\Phi$ such that $M$ is isomorphic to the vector matroid of $M(A)/C\backslash Y_1$. A matroid $M$ \textit{coconforms} to a template $\Phi$ if its dual $M^*$ conforms to $\Phi$.

Let $\mathcal{M}(\Phi)$ and $\mathcal{M}^*(\Phi)$ respectively denote the sets of matroids that conform and coconform to $\Phi$.  The next hypothesis is from \cite{gvz-prob} and is a modification of a hypothesis of Geelen, Gerards, and Whittle \cite{ggw15}. (See \cite{gvz-prob} for details about this modification, which was necessary to correct a problem for the non-binary case.) The hypothesis is believed to be true, but its proof is forthcoming in future papers by Geelen, Gerards, and Whittle.

\begin{hypothesis}[{\cite[Hypothesis 4.3, binary case]{gvz-prob}}]
\label{hyp:connected-template}
Let $\mathcal M$ be a proper minor-closed class of binary matroids. Then there exist $k\in\mathbb{Z}_+$ and templates $\Phi_1,\ldots,\Phi_s,\Psi_1,\ldots,\Psi_t$ such that
\begin{enumerate}
\item
$\mathcal{M}$ contains each of the classes $\mathcal{M}(\Phi_1),\ldots,\mathcal{M}(\Phi_s)$,
\item
$\mathcal{M}$ contains each of the classes $\mathcal{M}^*(\Psi_1),\ldots,\mathcal{M}^*(\Psi_t)$, and
\item
if $M$ is a $k$-connected member of $\mathcal M$ with at least $2k$ elements, then either $M$ is a member of at least one of the classes $\mathcal{M}(\Phi_1),$ $\ldots,$ $\mathcal{M}(\Phi_s)$, or $M^*$ is a member of at least one of the classes $\mathcal{M}(\Psi_1),\ldots,\mathcal{M}(\Psi_t)$.
\end{enumerate}
\end{hypothesis}

\section{Template Minors}
\label{sec:template minors}
The definitions and results in this section are from ~\cite{gvz17}. Proofs of the results can be found there as well. To simplify the proofs in ~\cite{gvz17}, it was helpful to expand the concept of conforming slightly.
\begin{definition}[{\cite[Definition 2.2]{gvz17}}]
 \label{virtual}
Let $A'$ be a matrix that respects $\Phi$, as defined in Section \ref{Binary Frame Templates}, except that we allow columns of $A'[B-X,Z]$ to be either unit columns or zero columns. Let $A$ be a matrix that is constructed from $A'$ as described in Section \ref{Binary Frame Templates}. Let $M$ be isomorphic to $M(A)/C\backslash Y_1$. We say that $A$ and $M$ \textit{virtually conform} to $\Phi$ and that $A'$ \textit{virtually respects} $\Phi$. If $M^*$ virtually conforms to $\Phi$, we say that $M$ \textit{virtually coconforms} to $\Phi$.
\end{definition}
We will denote the set of matroids that virtually conform to $\Phi$ by $\mathcal{M}_v(\Phi)$ and the set of matroids that virtually coconform to $\Phi$ by $\mathcal{M}^*_v(\Phi)$.

\begin{definition}[{\cite[Definition 3.1]{gvz17}}]
 A \textit{reduction} is an operation on a frame template $\Phi$ that produces a frame template $\Phi'$ such that $\mathcal{M}(\Phi')\subseteq \mathcal{M}(\Phi)$.
\end{definition}

\begin{proposition}[{\cite[Proposition 3.2]{gvz17}}] \label{reductions}
The following operations are reductions on a template $\Phi$:
\begin{itemize}
\item[(1)] Replace $\Lambda$ with a proper subgroup.
\item[(2)] Replace $\Delta$ with a proper subgroup.
\item[(3)] Remove an element $y$ from $Y_1$. (More precisely, replace $A_1$ with $A_1[X, Y_0\cup (Y_1-y)\cup C]$ and replace $\Delta$ with $\Delta|(Y_0\cup (Y_1-y)\cup C)$.)
\item[(4)] For all matrices $A'$ respecting $\Phi$, perform an elementary row operation on $A'[X, E]$. (Note that this alters the matrix $A_1$ and performs a change of basis on $\Lambda$.)
\item[(5)] If there is some element $x\in X$ such that, for every matrix $A'$ respecting $\Phi$, we have that $A'[\{x\},E]$ is a zero row vector, remove $x$ from $X$. (This simply has the effect of removing a zero row from every matrix conforming to $\Phi$.)
\item[(6)] Let $c\in C$ be such that $A_1[X,\{c\}]$ is a unit column whose nonzero entry is in the row indexed by $x\in X$, and let either $\lambda_x=0$ for each $\lambda\in\Lambda$ or $\delta_c=0$ for each $\delta\in\Delta$. We contract $c$ from every matroid conforming to $\Phi$ as follows. Let $\Delta'$ be the result of adding $-\delta_cA_1[\{x\},Y_0\cup Y_1\cup C]$ to each element $\delta\in\Delta$. Replace $\Delta$ with $\Delta'$, and then remove $c$ from $C$ and $x$ from $X$. (More precisely, replace $A_1$ with $A_1[X-x, Y_0\cup Y_1\cup (C-c)]$, replace $\Lambda$ with $\Lambda|(X-x)$, and replace $\Delta$ with $\Delta'|(Y_0\cup Y_1\cup (C-c))$.)
\item[(7)] Let $c\in C$ be such that $A_1[X,\{c\}]$ is a zero column and $\delta_c=0$ for all $\delta\in\Delta$. Then remove $c$ from $C$. (More precisely, replace $A_1$ with $A_1[X, Y_0\cup Y_1\cup (C-c)]$, and replace $\Delta$ with $\Delta|(Y_0\cup Y_1\cup (C-c))$.)
\end{itemize}
\end{proposition}

Since we always have $Y_0\subseteq E(M)$ for every matroid $M$ conforming to $\Phi$, operations (9)-(11) listed in the definition below are not reductions as defined above, but we continue the numbering from Proposition ~\ref{reductions} for ease of reference.

\begin{definition}[{\cite[Definition 3.3]{gvz17}}]
\label{weaklyconforming}
A template $\Phi'$ is a \textit{template minor} of $\Phi$ if $\Phi'$ is obtained from $\Phi$ by repeatedly performing the following operations:
\begin{itemize}
\item[(8)] Performing a reduction of type 1-7 on $\Phi$.
\item[(9)] Removing an element $y$ from $Y_0$, replacing $A_1$ with $A_1[X,(Y_0-y)\cup Y_1\cup C]$, and replacing $\Delta$ with $\Delta|((Y_0-y)\cup Y_1\cup C)$. (This has the effect of deleting $y$ from every matroid conforming to $\Phi$.)
\item[(10)] Let $x\in X$ with $\lambda_x=0$ for every $\lambda\in\Lambda$, and let $y\in Y_0$ be such that $(A_1)_{x,y}\neq0$. Then contract $y$ from every matroid conforming to $\Phi$. (More precisely, perform row operations on $A_1$ so that $A_1[X, \{y\}]$ is a unit column with $(A_1)_{x,y}=1$. Then replace every element $\delta\in\Delta$ with the row vector $-\delta_y A_1[\{x\}, Y_0\cup Y_1\cup C]+\delta$. This induces a group homomorphism $\Delta\rightarrow\Delta'$, where $\Delta'$ is also a subgroup of the additive group of $\mathrm{GF}(2)^{C\cup Y_0 \cup Y_1}$. Finally, replace $A_1$ with $A_1[X-x,(Y_0-y)\cup Y_1\cup C]$, project $\Lambda$ into $\mathrm{GF}(2)^{X-x}$, and project $\Delta'$ into $\mathrm{GF}(2)^{(Y_0-y)\cup Y_1\cup C}$. The resulting groups play the roles of $\Lambda$ and $\Delta$, respectively, in the new template.)
\item[(11)] Let $y\in Y_0$ with $\delta_y=0$ for every $\delta\in\Delta$. Then contract $y$ from every matroid conforming to $\Phi$. (More precisely, if $A_1[X, \{y\}]$ is a zero vector, this is the same as simply removing $y$ from $Y_0$. Otherwise, choose some $x\in X$ such that $(A_1)_{x,y}\neq0$. Then for every matrix $A'$ that respects $\Phi$, perform row operations so that $A_1[X,\{y\}]$ is a unit column with $(A_1)_{x,y}=1$. This induces a group isomorphism $\Lambda\rightarrow\Lambda'$ where $\Lambda'$ is also a subgroup of the additive group of $\mathrm{GF}(2)^X$. Finally, replace $A_1$ with $A_1[X-x,(Y_0-y)\cup Y_1\cup C]$, project $\Lambda'$ into $\mathrm{GF}(2)^{X-x}$, and project $\Delta$ into $\mathrm{GF}(2)^{(Y_0-y)\cup Y_1\cup C}$. The resulting groups play the roles of $\Lambda$ and $\Delta$, respectively, in the new template.)
\end{itemize}
\end{definition}

Throughout the rest of this paper, we will refer to the operations listed in Proposition ~\ref{reductions} and Definition ~\ref{weaklyconforming} as operations (1)-(11).

Let $\Phi'$ be a template minor of $\Phi$. Let $A$ be a matrix that virtually conforms to $\Phi'$, and let $M$ be a matroid that virtually conforms to $\Phi'$. We say that $A$ and $M$ \textit{weakly conform} to $\Phi$. Let $\mathcal{M}_w(\Phi)$ denote the set of matroids that weakly conform to $\Phi$, and let $\mathcal{M}^*_w(\Phi)$ denote the set of matroids whose duals weakly conform to $\Phi$. If $M\in\mathcal{M}^*_w(\Phi)$, we say that $M$ \textit{weakly coconforms} to $\Phi$.

\begin{lemma}[{\cite[Lemma 3.4]{gvz17}}] \label{minor}
 If a matroid $M$ weakly conforms to a template $\Phi$, then $M$ is a minor of a matroid that conforms to $\Phi$.
\end{lemma}

Since Hypothesis \ref{hyp:connected-template} deals with minor-closed classes, it can be generalized using the preceding lemma.

\begin{corollary}[{\cite[Modification of Corollary 3.5]{gvz17}}] \label{weakframe}
Suppose Hypothesis \ref{hyp:connected-template} holds. Let $\mathcal M$ be a proper minor-closed class of binary matroids. Then there exist $k\in\mathbb{Z}_+$ and templates $\Phi_1,\ldots,\Phi_s,\Psi_1,\ldots,\Psi_t$ such that
\begin{enumerate}
\item
$\mathcal{M}$ contains each of the classes $\mathcal{M}_w(\Phi_1),\dots,\mathcal{M}_w(\Phi_s)$,
\item
$\mathcal{M}$ contains each of the classes $\mathcal{M}^*_w(\Psi_1),\ldots,\mathcal{M}^*_w(\Psi_t)$, and
\item
if $M$ is a simple $k$-connected member of $\mathcal M$ with at least $2k$ elements, then either $M$ is a member of at least one of the classes $\mathcal{M}_v(\Phi_1),\ldots,\mathcal{M}_v(\Phi_s)$, or $M^*$ is a member of at least one of the classes $\mathcal{M}_v(\Psi_1),\ldots,\mathcal{M}_v(\Psi_t)$.
\end{enumerate}
\end{corollary}

If $\mathcal{M}_{w}(\Phi)=\mathcal{M}_{w}(\Phi')$, we say that $\Phi$ is \textit{equivalent} to $\Phi'$ and write $\Phi\sim\Phi'$. Note that $\sim$ is indeed an equivalence relation.

\begin{definition}[{\cite[Definition 3.6]{gvz17}}]
We define a preorder $\preceq$ on the set of templates as follows. We say $\Phi\preceq\Phi'$ if $\mathcal{M}_w(\Phi)\subseteq\mathcal{M}_w(\Phi')$. This is indeed a preorder since reflexivity and transitivity follow from the subset relation.
\end{definition}

Let $\Phi_0$ be the template with all groups trivial and all sets empty. We call this template the \textit{trivial template}. In general, we say that a template $\Phi$ is \textit{trivial} if $\Phi\preceq\Phi_0$. It is easy to see (by using operations (1), (2), (3), (6), and (9) as many times as needed) that for any template $\Phi$, we have $\Phi_0\preceq \Phi$. Therefore, if $\Phi\preceq\Phi_0$, then actually $\Phi\sim\Phi_0$.

We now describe a collection of minimal nontrivial templates.

\begin{definition}[{\cite[Definition 3.7]{gvz17}}]
\leavevmode
\begin{itemize}
\item Let $\Phi_C$  be the template with all groups trivial and all sets empty except that $|C|=1$ and $\Delta\cong\mathbb{Z}/2\mathbb{Z}$.
\item Let $\Phi_X$  be the template with all groups trivial and all sets empty except that $|X|=1$ and $\Lambda\cong\mathbb{Z}/2\mathbb{Z}$.
\item Let $\Phi_{Y_0}$  be the template with all groups trivial and all sets empty except that $|Y_0|=1$ and $\Delta\cong\mathbb{Z}/2\mathbb{Z}$. 
\item Let $\Phi_{CX}$ be the template with $Y_0=Y_1=\emptyset$, with $|C|=|X|=1$, with $\Delta\cong\Lambda\cong\mathbb{Z}/2\mathbb{Z}$, and with $A_1=[1]$.
\item Let $\Phi_{Y_1}$ be the template with all groups trivial, with $C=Y_0=\emptyset$, with $|Y_1|=3$ and $|X|=2$, and with $A_1=
\begin{bmatrix}
1& 0 &1\\
0& 1 & 1
\end{bmatrix}$.
\end{itemize}
\end{definition}

\begin{lemma}[{\cite[Lemma 3.8]{gvz17}}]
\label{YCD}
 The following relations hold:
\begin{itemize}
\item[(1)] $\Phi_{Y_1}\preceq\Phi_X$
\item[(2)] $\Phi_{Y_1}\preceq\Phi_C$
\item[(3)] $\Phi_{Y_0}\preceq\Phi_C$
\item[(4)] $\Phi_C\preceq\Phi_{CX}$
\item[(5)] $\Phi_X\preceq\Phi_{CX}$
\end{itemize}
\end{lemma}

\begin{lemma}[{\cite[Lemma 3.9]{gvz17}}] \label{yshift}
Let $\Phi$ be a template with $y\in Y_1$. Let $\Phi'$ be the template obtained from $\Phi$ by removing $y$ from $Y_1$ and placing it in $Y_0$. Then $\Phi'\preceq \Phi$.
\end{lemma}

We call the operation described in Lemma ~\ref{yshift} a \textit{$y$-shift}.

\begin{lemma}[{\cite[Lemma 3.12]{gvz17}}]
 \label{PhiD}
If $\Phi=(C,X,Y_0,Y_1,A_1,\Delta,\Lambda)$ is a binary frame template with $\Lambda|X_1$ nontrivial, then $\Phi_X\preceq\Phi$.
\end{lemma}

\begin{lemma}[{\cite[Lemma 3.17]{gvz17}}]
 \label{simpleY1}
If $\Phi$ is a template with $\Delta$ trivial, then $\Phi$ is equivalent to a template $\Phi'$ where $M(A_1[X,Y_1])$ is a simple matroid.
\end{lemma}

The next lemma can be extracted from the proof of \cite[Lemma 4.12]{gvz17}.

\begin{lemma}
\label{PG32Phi}
Let $\Phi$ be a template such that $\mathcal{M}_w(\Phi)\subseteq\mathcal{EX}(\mathrm{PG}(3,2))$. Then either $\Phi\preceq\Phi_X$ or $\Phi$ is equivalent to a template with $C=\emptyset$ and with $\Lambda$ and $\Delta$ trivial.
\end{lemma}

\begin{definition}[{\cite[Definition 4.4]{gvz17}}]
\label{Xr}
Let $X_r$ be the largest simple matroid of rank $r$ that virtually conforms to $\Phi_{Y_1}$.
\end{definition}

Note that $X_1=U_{1,1}$, and if $r\geq2$, then $X_r$ is the vector matroid of the following binary matrix, where we choose for the frame matrix the matrix representation of $M(K_{r-1})$, so that the identity matrices are $(r-2)\times(r-2)$ matrices. We will call this matrix $A_r$:

\begin{center}
\begin{tabular}{ |c|ccc|c|c|c|}
\hline
\multirow{2}{*}{0}&1&0&1&$1\cdots1$&$0\cdots0$&$1\cdots1$\\
&0&1&1&$0\cdots0$&$1\cdots1$&$1\cdots1$\\
\hline
frame matrix&\multicolumn{3}{c|}{0}&$I$&$I$&$I$\\
\hline
\end{tabular}
\end{center}

\begin{lemma}[{\cite[Lemma 4.5]{gvz17}}]
 \label{Y1minors}
The class $\mathcal{M}_v(\Phi_{Y_1})$ is the class of even-cycle matroids with a blocking pair. This class is minor-closed.
\end{lemma}

\section{Excluded Minors}
\label{sec:Excluded Minors}
In this Section, we will establish that the matroids from the introduction are indeed excluded minors for the various classes of matroids. Computations 1 and 2 in the Appendix give the SageMath code for functions that test whether a binary matroid is an excluded minor for the class of even-cycle and even-cut matroids, respectively. The code is based on the fact that an even-cycle matroid $M$ can be represented by a binary matrix with a row whose removal results in a matrix representing a graphic matroid. Thus, there is some binary extension $N$ of $M$ on ground set $E(M)\cup\{e\}$ such that $N/e$ is graphic. Therefore, to check if a binary matroid $M$ is even-cycle, it suffices to check if $N/e$ is graphic for some binary extension $N$ of $M$. If this is false for all such $N$, then $M$ is not even-cycle. The even-cut case is analogous.

\begin{theorem}
\label{proveexminors}
 Each of the matroids $\mathrm{PG}(3,2)\backslash e$, $L_{19}$, and $L_{11}$ is an excluded minor for the class of even-cycle matroids.
\end{theorem}

\begin{proof}
The largest even-cycle matroid of rank $r$ has a representation obtained by putting an odd edge in parallel with every even edge of $K_r$, and by adding an odd loop. Therefore, the size of the largest even-cycle matroid of rank $r$ is $2\binom{r}{2}+1=r^2-r+1$. Therefore, the matroid $\mathrm{PG}(3,2)\backslash e$, which has rank 4 and size 14 is too large to be even-cycle. Deletion of any element from $\mathrm{PG}(3,2)\backslash e$ results in the unique matroid (up to isomorphism) obtained from $\mathrm{PG}(3,2)$ by deleting two elements. This is exactly the largest simple even-cycle matroid of rank 4, as described above. Thus, deletion of any element from $\mathrm{PG}(3,2)\backslash e$ results in an even-cycle matroid. To see that contraction of any element from $\mathrm{PG}(3,2)\backslash e$ results in an even-cycle matroid, note that every binary matroid of rank 3 is even-cycle since removal of any row results in a matrix that obviously has at most two nonzero entries per column.

The fact that $L_{19}$ and $L_{11}$ are excluded minors for the class of even-cycle matroids was verified using SageMath in Computations 3 and 4 in the Appendix.
\end{proof}

\begin{theorem}
\label{evencutexcludedminors}
 The matroids $M(K_6)$ and $H^*_{12}$ are excluded minors for the class of even-cut matroids.
\end{theorem}

\begin{proof}
 This was verified using SageMath in Computations 5 and 6 in the Appendix.
\end{proof}

\begin{lemma}
\label{cosimplify}
A matroid is an even-cycle matroid with a blocking pair if and only if its cosimplification also is.
\end{lemma}

\begin{proof}
The class of even-cycle matroids with a blocking pair is minor-closed; therefore the cosimplification of an even-cycle matroid with a blocking pair will be such a matroid as well.

For the converse, let $M$ be even-cycle with a blocking pair, and consider an even-cycle representation of $M$ with a blocking pair. It suffices to consider coextensions $N$ of $M$, with $E(N)=E(M)\cup e$ and such that either $\{e,f\}$ is a series pair of $N$ or $e$ is a coloop of $N$. First, we consider the case where $\{e,f\}$ is a series pair. If $f$ is represented by an even edge in the even-cycle representation of $M$, then $e$ and $f$ in $N$ are represented by edges obtained by subdividing $f$ in $M$. This has no effect on the blocking pair. If $f$ is represented by an odd edge other than a loop, we resign at a vertex in the blocking pair that is incident with $f$. This maintains the blocking pair, but now $f$ is represented by an even edge as above. Now consider the case where $f$ is represented by an odd loop. Since $M$ is even-cycle with a blocking pair, recall from Lemma ~\ref{Y1minors} and Definition \ref{Xr} that $M$ is a restriction of a matroid represented by a matrix of the following form:
\begin{center}
\begin{tabular}{ |c|ccc|c|c|c|}
\multicolumn{1}{c}{}&$f$&\multicolumn{5}{c}{}\\
\hline
\multirow{2}{*}{0}&1&0&1&$1\cdots1$&$0\cdots0$&$1\cdots1$\\
&0&1&1&$0\cdots0$&$1\cdots1$&$1\cdots1$\\
\hline
frame matrix&\multicolumn{3}{c|}{0}&$I$&$I$&$I$\\
\hline
\end{tabular}
\end{center}
Since $N$ contains $\{e,f\}$ as a series pair, $N$ is a restriction of a matroid $N'$ represented by a matrix of the following form:
\begin{center}
\begin{tabular}{ |c|c|ccc|c|c|c|}
\multicolumn{1}{c}{$e$}&\multicolumn{1}{c}{}&$f$&\multicolumn{5}{c}{}\\
\hline
0&\multirow{2}{*}{0}&1&0&1&$1\cdots1$&$0\cdots0$&$1\cdots1$\\
0&&0&1&1&$0\cdots0$&$1\cdots1$&$1\cdots1$\\
\hline
0&frame matrix&\multicolumn{3}{c|}{0}&$I$&$I$&$I$\\
\hline
1&$0\cdots\cdots\cdots\cdots0$&1&0&0&$0\cdots0$&$0\cdots0$&$0\cdots0$\\
\hline
\end{tabular}
\end{center}
By Lemma ~\ref{Y1minors}, $N'$ is even-cycle with a blocking pair. Therefore, so is $N$.

Lastly, we consider the case where $e$ is a coloop of $N$. Then $N$ can be represented by a graph obtained from the graph representing $M$ by adding a new vertex and joining it to any other vertex with an even edge. The blocking pair is maintained.
\end{proof}

\begin{theorem}
\label{proveevencyclebp}
The matroids $\mathrm{PG}(3,2)\backslash L$ and $M^*(K_6)$ are excluded minors for the class of even-cycle matroids with a blocking pair.
\end{theorem}

\begin{proof}
 Recall from Definition ~\ref{Xr} and Lemma ~\ref{Y1minors} that $X_r$ is the largest simple matroid of rank $r$ that is even-cycle with a blocking pair. Note that $X_4$ is the matroid obtained from $\mathrm{PG}(3,2)$ by deleting an independent set of size 3. Therefore, $\mathrm{PG}(3,2)\backslash L$ is not a restriction of $X_4$. By Lemma ~\ref{Y1minors}, $\mathrm{PG}(3,2)\backslash L$ is not an even-cycle matroid with a blocking pair. However, since $X_3=\mathrm{PG}(2,2)=F_7$, all binary matroids of rank at most 3 are even-cycle matroids with blocking pairs. Therefore, $\mathrm{PG}(3,2)\backslash L/e$ is even-cycle with a blocking pair for each element $e$ of $\mathrm{PG}(3,2)\backslash L$. Moreover, by deleting any element from $\mathrm{PG}(3,2)\backslash L$, we obtain a restriction of $X_4$. Therefore, $\mathrm{PG}(3,2)\backslash L$ is an excluded minor for the class of even-cycle matroids with a blocking pair.

By Theorem \ref{evencutexcludedminors}, $M(K_6)$ is not an even-cut matroid. Recall from Section \ref{introduction} that the dual of an even-cycle matroid with a blocking pair is an even cut matroid. Therefore, $M^*(K_6)$ is not an even-cycle matroid with a blocking pair. It remains to show that $M^*(K_6\backslash e)$ and $M^*(K_6/e)$ are even-cycle with a blocking pair. By Lemma ~\ref{cosimplify}, $M^*(K_6/e)$ has an even-cycle representation with a blocking pair if and only if $M^*(K_5)$ does. Even-cycle representations of $M^*(K_5)$ and $M^*(K_6\backslash e)$, with odd edges printed in bold, are given in Figure \begin{figure}
\begin{center}
\begin{tikzpicture}[line cap=round,line join=round,>=triangle 45,x=1.0cm,y=1.0cm]
\clip(2.54,1.72) rectangle (11.42,5.48);
\draw (4.,5.)-- (3.,4.3);
\draw (3.,4.3)-- (3.5,3.2);
\draw (3.,4.3)-- (4.,4.);
\draw (4.,5.)-- (5.,4.3);
\draw (5.,4.3)-- (4.,4.);
\draw (5.,4.3)-- (4.5,3.2);
\draw [line width=2.pt] (4.,5.)-- (4.,4.);
\draw [line width=2.pt] (4.,4.)-- (3.5,3.2);
\draw [line width=2.pt] (4.,4.)-- (4.5,3.2);
\draw [line width=2.pt] (4.5,3.2)-- (3.5,3.2);
\draw (7.,5.)-- (9.,5.);
\draw (9.,5.)-- (11.,5.);
\draw (11.,5.)-- (11.,4.);
\draw (11.,4.)-- (11.,3.);
\draw (11.,3.)-- (9.,3.);
\draw (7.,3.)-- (7.,5.);
\draw [line width=2.pt] (9.,4.)-- (9.,3.);
\draw (9.,4.)-- (9.,5.);
\draw [line width=2.pt] (9.,4.)-- (10.,4.);
\draw [line width=2.pt] (10.,4.)-- (11.,5.);
\draw (10.,4.)-- (11.,4.);
\draw (3.48,2.84) node[anchor=north west] {$M^*(K_5)$};
\draw (8.12,2.82) node[anchor=north west] {$M^*(K_6\backslash e)$};
\begin{scriptsize}

\draw [fill=black] (4.,4.) circle (2.0pt);
\draw [fill=black] (4.,5.) circle (2.0pt);
\draw [fill=black] (3.5,3.2) circle (2.0pt);
\draw [fill=black] (4.5,3.2) circle (2.0pt);
\draw [fill=black] (3.,4.3) circle (2.0pt);
\draw [fill=black] (5.,4.3) circle (2.0pt);
\draw [fill=black] (7.,5.) circle (1.5pt);
\draw [fill=black] (9.,5.) circle (1.5pt);
\draw [fill=black] (11.,5.) circle (1.5pt);
\draw [fill=black] (11.,3.) circle (1.5pt);
\draw [fill=black] (11.,4.) circle (1.5pt);
\draw [fill=black] (10.,4.) circle (1.5pt);
\draw [fill=black] (9.,4.) circle (1.5pt);
\end{scriptsize}
\vertex [fill](a) at (9,3) {};
\vertex [fill](b) at (7,3) {};
\path
(a) edge[bend right=40]  (b)
(a) edge[line width=2.0pt]  (b);
\end{tikzpicture}
\end{center}
\caption{Even-cycle representations of $M^*(K_5)$ and $M^*(K_6\backslash e)$}
  \label{fig:M*K5K6e}
\end{figure} \ref{fig:M*K5K6e}. Each of these representations have blocking pairs.
\end{proof}

\section{Some Technical Lemmas Proved with SageMath: Even-Cycle Matroids}
\label{sec:technical-cycle}
In this section, we list several technical lemmas that we need to prove Theorems \ref{connectedevencycle} and \ref{blockingpair}. Many of the proofs will merely refer the reader to a computation in the Appendix. The computations use the SageMath software system \cite{sage}. In Lemmas \ref{A.3}-\ref{A.22}, $\Phi$ is a template with $C=\emptyset$ and all groups trivial. Moreover, there are binary matrices $P_0$ and $P_1$ and a partition $Y_1=Y'_1\cup Y''_1$ such that $A_1[X,Y'_1]$ is an identity matrix, $A_1[X,Y''_1]=P_1$, and $A_1[X,Y_0]=P_0$. The reader may prefer to move on to Section \ref{Even-Cycle Matroids}, referring to Section \ref{sec:technical-cycle} as necessary.

Computation 7 gives a SageMath function which builds the largest possible simple matroid $M$ of rank $r(M(A_1))+n-1$ that virtually conforms to $\Phi$. It does this by choosing for the frame matrix the matrix representation of $M(K_n)$ and by including all possible elements of $Z$ that can be constructed from elements of $Y_1$. For \texttt{AY0}, we input the matrix $A_1[X,Y_0]$. For \texttt{AY1}, we input the matrix $A_1[X,Y_1]$. In most of Computations 8-41, we then test if $M$ contains $\mathrm{PG}(3,2)\backslash e$ as a minor by looking for a subset $S$ of the ground set of $M$ such that $r(M/S)=4$ and $|$si$(M/S)|\geq14$. If this subset exists, then $M$ must contain $\mathrm{PG}(3,2)\backslash e$ as a minor. In the Python programming language, on which SageMath is based, a set of size $n$ has elements labeled $0,1,\dots,n-1$. Thus, for example, if $S=\{17,12,7\}$, then the 18th, 13th, and 8th columns are to be contracted. We can similarly test for a $\mathrm{PG}(3,2)_{-2}$ minor.

\begin{lemma}
\label{A.3}If $P_1$ contains a column with four or more nonzero entries, then $\mathcal{M}_w(\Phi)\nsubseteq\mathcal{EX}(\mathrm{PG}(3,2)\backslash e)$. More generally, if $M(A_1[X,Y_1])$ contains a circuit of size at least 5, then $\mathcal{M}_w(\Phi)\nsubseteq\mathcal{EX}(\mathrm{PG}(3,2)\backslash e)$.
\end{lemma}
\begin{proof}
See Computation 8 in the Appendix.
\end{proof}

\begin{lemma}\label{A.4}
If $P_1$ contains the submatrix $\left[
\begin{matrix}
1 & 0 \\
1 & 0 \\
0 & 1 \\
0 & 1 \\
\end{matrix}
\right]$, then $\mathcal{M}_w(\Phi)\nsubseteq\mathcal{EX}(\mathrm{PG}(3,2)\backslash e)$.
\end{lemma}
\begin{proof}
See Computation 9 in the Appendix.
\end{proof}

\begin{lemma}\label{A.5}
If $P_1$ contains the submatrix
$\left[
\begin{matrix}
1 & 0 &1 \\
1 & 1 &0 \\
0 & 1 &1 \\
\end{matrix}
\right]$, then $\mathcal{M}_w(\Phi)\nsubseteq\mathcal{EX}(\mathrm{PG}(3,2)\backslash e)$.
\end{lemma}
\begin{proof}
See Computation 10 in the Appendix.
\end{proof}

\begin{lemma}\label{A.6}
If $P_1$ contains the submatrix
$\left[
\begin{array}{ccc}
1 & 1 & 1 \\
1 & 0 & 1 \\
0 & 1 & 1 \\
\end{array}
\right]$, then $\mathcal{M}_w(\Phi)\nsubseteq\mathcal{EX}(\mathrm{PG}(3,2)\backslash e)$.
\end{lemma}
\begin{proof}
See Computation 11 in the Appendix.
\end{proof}

\begin{lemma}\label{A.7} If $P_1$ contains the submatrix $P_1'=
\left[
\begin{matrix}
1 & 0\\
1 & 0\\
1 & 1\\
0 & 1\\
\end{matrix}
\right]$, then $\mathcal{M}_w(\Phi)\nsubseteq\mathcal{EX}(\mathrm{PG}(3,2)\backslash e)$.
\end{lemma}
\begin{proof}In the matrix $[I_4|P_1']$, columns 1, 2, 4, 5, and 6 form a circuit of size 5, which is forbidden by Lemma \ref{A.3}.
\end{proof}

\begin{lemma}\label{A.8}
If $P_0$ contains a column with five nonzero entries or either of the following matrices below as a submatrix, then $\mathcal{M}_w(\Phi)\nsubseteq\mathcal{EX}(\mathrm{PG}(3,2)\backslash e)$:
\[
\left[
\begin{array}{cc}
1&0\\
1&0\\
1&0\\
0&1\\
0&1\\
0&1\\
\end{array}
\right],
\left[
\begin{array}{cc}
1&0\\
1&0\\
1&0\\
1&1\\
0&1\\
0&1\\
\end{array}
\right].
\]
\end{lemma}
\begin{proof}
If we contract a column from $Y_0$ (in other words, we perform operation (10) on $\Phi$), then a column from the identity matrix $A_1[X,Y_1']$ becomes a column in the resulting $A_1[X,Y_1'']$. We see then that if $P_0$ contains one of the submatrices listed in this result, then contraction of the columns of that submatrix results in the submatrix forbidden in Lemma \ref{A.3}, \ref{A.4}, or \ref{A.7}.
\end{proof}

The following conditions will be used as conclusions for Lemmas \ref{A.9}-\ref{A.14}.
\begin{enumerate}
 \item $\mathcal{M}_w(\Phi)\nsubseteq \mathcal{EX}(\mathrm{PG}(3,2)\backslash e, L_{11})$ and $\mathcal{M}_w(\Phi)\nsubseteq \mathcal{EX}(\mathrm{PG}(3,2)_{-2})$,
\item $P_0$ is of the following form, where each column of $[H_1|H_0]$ has at most two nonzero entries:
\[
\left[
\begin{array}{c|c}
1\cdots1&0\cdots0\\
\hline
H_1&H_0\\
\end{array}
\right], \textnormal{or}\]
\item $P_0$ is of the form 
\[
\left[
\begin{array}{c|c|c|c}
1\cdots1&1\cdots1&0\cdots0&0\cdots0\\
1\cdots1&0\cdots0&1\cdots1&0\cdots0\\
\hline
H_{1,1}&H_{1,0}&H_{0,1}&H_{0,0}\\
\end{array}
\right],
\]
where each column of $H_{1,1}$ and $H_{0,0}$ has at most two nonzero entries, and where each column of $H_{0,1}$ and $H_{1,0}$ is a unit column or a zero column.
\end{enumerate}

\begin{lemma}\label{A.9}
If $P_0$ contains the columns $\left[
\begin{matrix}
1 & 1\\
1 & 0\\
1 & 0\\
0 & 1\\
0 & 1\\
0 & 0\\
\vdots & \vdots\\
0 & 0\\
\end{matrix}
\right]$, then (1) or (2) holds.
\end{lemma}

\begin{proof}
Suppose neither (1) nor (2) holds. Since (1) does not hold, Lemma \ref{A.8} implies that $P_0$ contains none of the following submatrices, whether $x=0$ or $x=1$, and also contains no column with five nonzero entries:

\[
\left[
\begin{array}{ccc}
1 & 1 &x\\
1 & 0 &0\\
1 & 0 &0\\
0 & 1 &1\\
0 & 1 &1\\
0 & 0 &1\\
\end{array}
\right],
\left[
\begin{array}{ccc}
1 & 1 &x\\
1 & 0 &0\\
1 & 0 &0\\
0 & 1 &0\\
0 & 1 &1\\
0 & 0 &1\\
0 & 0 &1\\
\end{array}
\right],
\left[
\begin{array}{ccc}
1 & 1 &x\\
1 & 0 &0\\
1 & 0 &0\\
0 & 1 &0\\
0 & 1 &0\\
0 & 0 &1\\
0 & 0 &1\\
0 & 0 &1
\end{array}
\right].
\]

Therefore, since (2) does not hold, $P_0$ contains one of the following submatrices:
\[
\left[
\begin{array}{ccc}
1&1&1\\
1&0&1\\
1&0&0\\
0&1&1\\
0&1&0\\
0&0&1\\
\end{array}
\right],
\left[
\begin{array}{ccc}
1&1&0\\
1&0&1\\
1&0&0\\
0&1&1\\
0&1&0\\
0&0&1\\
\end{array}
\right],
\left[
\begin{array}{ccc}
1&1&0\\
1&0&1\\
1&0&1\\
0&1&1\\
0&1&0\\
\end{array}
\right],
\left[
\begin{array}{ccc}
1&1&1\\
1&0&1\\
1&0&1\\
0&1&1\\
0&1&0\\
\end{array}
\right],
\left[
\begin{array}{ccc}
1&1&0\\
1&0&1\\
1&0&1\\
0&1&1\\
0&1&1\\
\end{array}
\right].
\]
Computations 12-16 in the Appendix show that each of these submatrices results in a template to which virtually conforms a matroid which has $\mathrm{PG}(3,2)\backslash e$ as a minor or which has both $\mathrm{PG}(3,2)_{-2}$ and $L_{11}$ as minors, contradicting the assumption that (1) does not hold.
\end{proof}

\begin{lemma}\label{A.10} If $P_0$ contains the columns $\left[
\begin{array}{cc}
1 & 1\\
1 & 1\\
1 & 0\\
1 & 0\\
0 & 1\\
0 & 0\\
\vdots&\vdots\\
0 & 0\\
\end{array}
\right],$ then (1) or (3) holds.
\end{lemma}
\begin{proof} Suppose neither (1) nor (3) holds. Since (1) does not hold, Lemma \ref{A.8} implies that $P_0$ contains none of the following submatrices and also contains no column with five nonzero entries:
\[
\left[
\begin{array}{ccc}
1&1&0\\
1&1&0\\
1&0&1\\
1&0&1\\
0&1&0\\
0&0&1\\
\end{array}
\right],
\left[
\begin{array}{ccc}
1&1&0\\
1&1&0\\
1&0&1\\
1&0&0\\
0&1&1\\
0&0&1\\
\end{array}
\right],
\left[
\begin{array}{ccc}
1&1&0\\
1&1&0\\
1&0&1\\
1&0&0\\
0&1&0\\
0&0&1\\
0&0&1\\
\end{array}
\right],
\left[
\begin{array}{ccc}
1&1&0\\
1&1&0\\
1&0&0\\
1&0&0\\
0&1&1\\
0&0&1\\
0&0&1\\
\end{array}
\right],\]
\[
\left[
\begin{array}{ccc}
1&1&0\\
1&1&0\\
1&0&0\\
1&0&0\\
0&1&0\\
0&0&1\\
0&0&1\\
0&0&1\\
\end{array}
\right],
\left[
\begin{array}{ccc}
1&1&1\\
1&1&0\\
1&0&0\\
1&0&0\\
0&1&1\\
0&0&1\\
\end{array}
\right],
\left[
\begin{array}{ccc}
1&1&1\\
1&1&0\\
1&0&0\\
1&0&0\\
0&1&0\\
0&0&1\\
0&0&1\\
\end{array}
\right].
\]

Therefore, since (3) does not hold, $P_0$ contains one of the following submatrices:
\[
\left[
\begin{array}{ccc}
1&1&0\\
1&1&0\\
1&0&1\\
1&0&1\\
0&1&1\\
\end{array}
\right],
\left[
\begin{array}{ccc}
1&1&1\\
1&1&0\\
1&0&1\\
1&0&1\\
0&1&1\\
\end{array}
\right],
\left[
\begin{array}{ccc}
1&1&1\\
1&1&0\\
1&0&1\\
1&0&1\\
0&1&0\\
\end{array}
\right],
\left[
\begin{array}{ccc}
1&1&1\\
1&1&0\\
1&0&1\\
1&0&0\\
0&1&1\\
\end{array}
\right],
\left[
\begin{array}{ccc}
1&1&1\\
1&1&0\\
1&0&1\\
1&0&0\\
0&1&0\\
0&0&1\\
\end{array}
\right].
\]
Computations 17-21 in the Appendix show that each of these submatrices results in a template to which virtually conforms a matroid which has $\mathrm{PG}(3,2)\backslash e$ as a minor or which has both $\mathrm{PG}(3,2)_{-2}$ and $L_{11}$ as minors, contradicting the assumption that (1) does not hold.
\end{proof}

\begin{lemma}\label{A.11}
 If $P_0$ contains the columns $\left[
\begin{array}{cc}
1 & 1\\
1 & 1\\
1 & 0\\
0 & 1\\
0 & 0\\
\vdots & \vdots\\
0 & 0\\
\end{array}
\right],$ then (1), (2), or (3) holds.
\end{lemma}
\begin{proof}
Suppose neither (1), (2), nor (3) holds. First, suppose that $P_0$ contains one of the following submatrices:
\[
\left[
\begin{array}{ccc}
1&1&0\\
1&1&0\\
1&0&1\\
0&1&1\\
0&0&1\\
\end{array}
\right],
\left[
\begin{array}{ccc}
1&1&1\\
1&1&0\\
1&0&1\\
0&1&0\\
0&0&1\\
\end{array}
\right],
\left[
\begin{array}{ccc}
1&1&1\\
1&1&0\\
1&0&0\\
0&1&1\\
0&0&1\\
\end{array}
\right],
\left[
\begin{array}{ccc}
1&1&1\\
1&1&0\\
1&0&0\\
0&1&0\\
0&0&1\\
0&0&1\\
\end{array}
\right].
\]
For each of these four submatrices, since (1) does not hold, Lemma \ref{A.8} implies that the column of $P_0$ containing the third column of the submatrix contains no nonzero entries except those in the submatrix. But then these submatrices are forbidden by Lemma \ref{A.9} (with the third column playing the role of one of the two columns in Lemma \ref{A.9}) and the assumption that (2) does not hold. Moreover, the following submatrices are forbidden by Lemma \ref{A.8}:
\[
\left[
\begin{array}{ccc}
1&1&0\\
1&1&0\\
1&0&1\\
0&1&0\\
0&0&1\\
0&0&1\\
\end{array}
\right],
\left[
\begin{array}{ccc}
1&1&0\\
1&1&0\\
1&0&0\\
0&1&0\\
0&0&1\\
0&0&1\\
0&0&1\\
\end{array}
\right].
\]

Therefore, since (3) does not hold, $P_0$ must contain the following submatrix:
\[
\left[
\begin{array}{ccc}
1&1&1\\
1&1&0\\
1&0&1\\
0&1&1\\
\end{array}
\right].
\]
Since no column of $P_0$ can contain five nonzero entries, if this submatrix is contained in the following submatrix:
\[
\left[
\begin{array}{ccc}
1&1&1\\
1&1&0\\
1&0&1\\
0&1&1\\
0&0&1\\
\end{array}
\right],
\]
then the column of $P_0$ containing the third column of this matrix contains no other nonzero entries. Therefore, this submatrix is forbidden by Lemma \ref{A.10} and the assumption that (3) does not hold. Thus, $P_0$ contains the following three columns:
\[
\left[
\begin{array}{ccc}
1&1&1\\
1&1&0\\
1&0&1\\
0&1&1\\
0&0&0\\
\vdots&\vdots&\vdots\\
0&0&0\\
\end{array}
\right].
\]

Now, since (2) does not hold, $P_0$ contains one of the following submatrices, where $x$ is either $0$ or $1$:
\[
\left[
\begin{matrix}
1&1&1&x\\
1&1&0&1\\
1&0&1&1\\
0&1&1&1\\
\end{matrix}
\right],
\left[
\begin{matrix}
1&1&1&x\\
1&1&0&0\\
1&0&1&1\\
0&1&1&1\\
0&0&0&1\\
\end{matrix}
\right],
\left[
\begin{matrix}
1&1&1&x\\
1&1&0&0\\
1&0&1&0\\
0&1&1&1\\
0&0&0&1\\
0&0&0&1\\
\end{matrix}
\right],
\left[
\begin{matrix}
1&1&1&x\\
1&1&0&0\\
1&0&1&0\\
0&1&1&0\\
0&0&0&1\\
0&0&0&1\\
0&0&0&1
\end{matrix}
\right].
\]
The first of these four submatrices is forbidden by Computations 22 and 23 in the Appendix. If $x=0$, the third and fourth of these submatrices are forbidden by Lemma \ref{A.8}. The second submatrix is forbidden by Lemma \ref{A.8} if the column of $P_0$ containing the fourth column contains an additional nonzero entry, but it is forbidden by Lemma \ref{A.9} (and the assumption that (2) does not hold) otherwise. If $x=1$, the second matrix is forbidden by Lemma \ref{A.10} (along with the assumption that (3) does not hold and the fact that no column  of $P_0$ can contain five nonzero entries). The third and fourth submatrices are forbidden by Lemma \ref{A.8} (and the assumption that (1) does not hold). This completes the proof by contradiction.
\end{proof}

\begin{lemma}\label{A.12}
If $P_0$ contains the columns $\left[
\begin{array}{cc}
1 & 1\\
1 & 1\\
1 & 1\\
1 & 0\\
0 & 0\\
\vdots&\vdots\\
0 & 0\\
\end{array}
\right],$ then (1), (2), or (3) holds.
\end{lemma}
\begin{proof}
Suppose neither (1), (2), nor (3) holds. Consider the following matrices:
\[
\left[
\begin{matrix}
1&1&1\\
1&1&0\\
1&1&1\\
1&0&1\\
0&0&1\\
\end{matrix}
\right],
\left[
\begin{matrix}
1&1&1\\
1&1&0\\
1&1&1\\
1&0&0\\
0&0&1\\
0&0&1\\
\end{matrix}
\right],
\left[
\begin{matrix}
1&1&1\\
1&1&0\\
1&1&0\\
1&0&1\\
0&0&1\\
0&0&1\\
\end{matrix}
\right],
\left[
\begin{matrix}
1&1&1\\
1&1&0\\
1&1&0\\
1&0&0\\
0&0&1\\
0&0&1\\
0&0&1\\
\end{matrix}
\right],\]
\[\left[
\begin{matrix}
1&1&0\\
1&1&0\\
1&1&1\\
1&0&1\\
0&0&1\\
0&0&1\\
\end{matrix}
\right],
\left[
\begin{matrix}
1&1&0\\
1&1&0\\
1&1&1\\
1&0&0\\
0&0&1\\
0&0&1\\
0&0&1\\
\end{matrix}
\right],
\left[
\begin{matrix}
1&1&0\\
1&1&0\\
1&1&0\\
1&0&1\\
0&0&1\\
0&0&1\\
0&0&1\\
\end{matrix}
\right],
\left[
\begin{matrix}
1&1&0\\
1&1&0\\
1&1&0\\
1&0&0\\
0&0&1\\
0&0&1\\
0&0&1\\
0&0&1\\
\end{matrix}
\right].
\]

Since no column of $P_0$ can have five nonzero entries, the first two of these matrices are forbidden by Lemma \ref{A.10}. The other six are forbidden by Lemma \ref{A.8}. Now, for all but the first of these matrices, suppose it is not a submatrix of $P_0$, but the matrix formed by removing the last row is a submatrix of $P_0$. For the second matrix, this is impossible by Lemma \ref{A.11}. For the third, fourth, fifth, and sixth matrices, this is impossible by Lemma \ref{A.9}. For the seventh and eighth matrices, this is impossible by Lemma \ref{A.8}. Thus, since (3) does not hold, $P_0$ must contain $\left[
\begin{matrix}
1&1&1\\
1&1&0\\
1&1&1\\
1&0&1\\
\end{matrix}
\right]$ as a submatrix. Therefore, after swapping the second and third rows, we may assume without loss of generality that $P_0$ contains the following three columns:
\[
\left[
\begin{matrix}
1 & 1 & 1\\
1 & 1 & 1\\
1 & 1 & 0\\
1 & 0 & 1\\
0 & 0 & 0\\
\vdots&\vdots&\vdots\\
0 & 0 & 0\\
\end{matrix}
\right].
\]
However, since (3) does not hold, $P_0$ must contain a fourth column resulting in $P_0$ containing as a submatrix either one of the eight forbidden matrices above or $\left[
\begin{matrix}
1 & 1 & 1&0\\
1 & 1 & 1&1\\
1 & 1 & 0&1\\
1 & 0 & 1&1\\
\end{matrix}
\right]$, which is forbidden by Computation 22 in the Appendix. This completes the proof by contradiction.
\end{proof}

\begin{lemma}\label{A.13}
If $P_0$ contains the submatrix $\left[
\begin{array}{cc}
1 & 1\\
1 & 1\\
1 & 0\\
1 & 0\\
0 & 1\\
0 & 1\\
\end{array}
\right]$, then (1) or (3) holds.
\end{lemma}
\begin{proof} Suppose neither (1) nor (3) holds. The following submatrices are forbidden from $P_0$ by Lemma \ref{A.8}:
\[
\left[
\begin{array}{ccc}
1&1&0\\
1&1&0\\
1&0&1\\
1&0&0\\
0&1&1\\
0&1&1\\
\end{array}
\right],
\left[
\begin{array}{ccc}
1&1&0\\
1&1&0\\
1&0&1\\
1&0&0\\
0&1&1\\
0&1&0\\
0&0&1\\
\end{array}
\right],
\left[
\begin{array}{ccc}
1&1&0\\
1&1&0\\
1&0&0\\
1&0&0\\
0&1&1\\
0&1&1\\
0&0&1\\
\end{array}
\right],
\left[
\begin{array}{ccc}
1&1&0\\
1&1&0\\
1&0&0\\
1&0&0\\
0&1&0\\
0&1&1\\
0&0&1\\
0&0&1\\
\end{array}
\right],\]
\[\left[
\begin{array}{ccc}
1&1&0\\
1&1&0\\
1&0&0\\
1&0&0\\
0&1&0\\
0&1&0\\
0&0&1\\
0&0&1\\
0&0&1\\
\end{array}
\right],
\left[
\begin{array}{ccc}
1&1&1\\
1&1&0\\
1&0&1\\
1&0&1\\
0&1&0\\
0&1&0\\
\end{array}
\right],
\left[
\begin{array}{ccc}
1&1&1\\
1&1&0\\
1&0&1\\
1&0&0\\
0&1&0\\
0&1&0\\
0&0&1\\
\end{array}
\right],
\left[
\begin{array}{ccc}
1&1&1\\
1&1&0\\
1&0&0\\
1&0&0\\
0&1&0\\
0&1&0\\
0&0&1\\
0&0&1\\
\end{array}
\right].
\]

Therefore, since (3) does not hold and since no column of $P_0$ contains five nonzero entries, $P_0$ contains one of the following submatrices:
\[
\left[
\begin{array}{ccc}
1&1&0\\
1&1&0\\
1&0&1\\
1&0&1\\
0&1&1\\
0&1&1\\
\end{array}
\right],
\left[
\begin{array}{ccc}
1&1&1\\
1&1&0\\
1&0&1\\
1&0&1\\
0&1&1\\
0&1&0\\
\end{array}
\right],
\left[
\begin{array}{ccc}
1&1&1\\
1&1&0\\
1&0&1\\
1&0&0\\
0&1&1\\
0&1&0\\
\end{array}
\right].
\]

Computations 24-26 in the Appendix show that each of these submatrices results in a template to which virtually conforms a matroid which has $\mathrm{PG}(3,2)\backslash e$ as a minor or which has both $\mathrm{PG}(3,2)_{-2}$ and $L_{11}$ as minors, contradicting the assumption that (1) does not hold.
\end{proof}

\begin{lemma}\label{A.14}
If $P_0$ contains the submatrix $\left[
\begin{array}{cc}
1 & 1\\
1 & 1\\
1 & 1\\
1 & 0\\
0 & 1\\
\end{array}
\right]$, then (1) or (3) holds.
\end{lemma}
\begin{proof}
Suppose neither (1) nor (3) holds. The following submatrices are forbidden from $P_0$ by Lemma \ref{A.8}:
\[
\left[
\begin{matrix}
1&1&1\\
1&1&0\\
1&1&0\\
1&0&0\\
0&1&1\\
0&0&1\\
\end{matrix}
\right],
\left[
\begin{matrix}
1&1&1\\
1&1&0\\
1&1&0\\
1&0&0\\
0&1&0\\
0&0&1\\
0&0&1\\
\end{matrix}
\right],
\left[
\begin{matrix}
1&1&0\\
1&1&0\\
1&1&1\\
1&0&0\\
0&1&1\\
0&0&1\\
\end{matrix}
\right],
\left[
\begin{matrix}
1&1&0\\
1&1&0\\
1&1&1\\
1&0&0\\
0&1&0\\
0&0&1\\
0&0&1\\
\end{matrix}
\right],\]
\[\left[
\begin{matrix}
1&1&0\\
1&1&0\\
1&1&0\\
1&0&1\\
0&1&1\\
0&0&1\\
\end{matrix}
\right],
\left[
\begin{matrix}
1&1&0\\
1&1&0\\
1&1&0\\
1&0&0\\
0&1&1\\
0&0&1\\
0&0&1\\
\end{matrix}
\right],
\left[
\begin{matrix}
1&1&0\\
1&1&0\\
1&1&0\\
1&0&0\\
0&1&0\\
0&0&1\\
0&0&1\\
0&0&1\\
\end{matrix}
\right].
\]

The following submatrices are forbidden by Lemma \ref{A.13}. Moreover, for each of these matrices, if the given matrix is not a submatrix of $P_0$, then the submatrix obtained from this matrix by deleting the last row is forbidden from $P_0$ by Lemma \ref{A.10}:
\[
\left[
\begin{matrix}
1&1&1\\
1&1&0\\
1&1&1\\
1&0&0\\
0&1&1\\
0&0&1\\
\end{matrix}
\right],
\left[
\begin{matrix}
1&1&1\\
1&1&0\\
1&1&1\\
1&0&0\\
0&1&0\\
0&0&1\\
0&0&1\\
\end{matrix}
\right],
\left[
\begin{matrix}
1&1&1\\
1&1&0\\
1&1&0\\
1&0&1\\
0&1&1\\
0&0&1\\
\end{matrix}
\right],
\left[
\begin{matrix}
1&1&0\\
1&1&0\\
1&1&1\\
1&0&1\\
0&1&1\\
0&0&1\\
\end{matrix}
\right].
\]

Therefore, since (3) does not hold, $P_0$ contains a submatrix obtained from the following matrix by deleting either the third or fourth column:
\[\left[
\begin{matrix}
1&1&1&0\\
1&1&0&1\\
1&1&1&1\\
1&0&1&1\\
0&1&1&1\\
\end{matrix}
\right].
\]
This matrix itself is not a submatrix of $P_0$ because it contains a submatrix forbidden by Computation 23 in the Appendix. Thus, after swapping the second and third rows, we may assume without loss of generality that $P_0$ contains the following submatrix:
\[\left[
\begin{matrix}
1&1&1\\
1&1&1\\
1&1&0\\
1&0&1\\
0&1&1\\
\end{matrix}
\right].
\]

Since (3) does not hold and since we have already forbidden all other possibilities, $P_0$ must contain the following submatrix:
\[\left[
\begin{matrix}
1&1&1&1\\
1&1&1&0\\
1&1&0&1\\
1&0&1&1\\
0&1&1&1\\
\end{matrix}
\right],
\]
but this contains a submatrix forbidden by Computation 23 in the Appendix. This completes the proof by contradiction.
\end{proof}

\begin{lemma}\label{A.15} If $[P_1|P_0]$ contains  any of the following submatrices, with the submatrix to the left of the vertical line contained in $P_1$, with the column to the right of the vertical line contained in $P_0$, and with either $x=0$ or $x=1$, then $\mathcal{M}_w(\Phi)\nsubseteq\mathcal{EX}(\mathrm{PG}(3,2)\backslash e)$:
\[
\left[
\begin{array}{cc|c}
1&1&x\\
1&0&1\\
0&1&1\\
0&0&1\\
\end{array}
\right],
\left[
\begin{array}{cc|c}
1&1&x\\
1&0&0\\
0&1&1\\
0&0&1\\
0&0&1\\
\end{array}
\right],
\left[
\begin{array}{cc|c}
1&1&x\\
1&0&0\\
0&1&0\\
0&0&1\\
0&0&1\\
0&0&1\\
\end{array}
\right]
\]
\end{lemma}
\begin{proof}See Computations 27-32 in the Appendix. \end{proof}

\begin{lemma}\label{A.16}
If $[P_1|P_0]$ contains  any of the following submatrices, with the submatrix to the left of the vertical line contained in $P_1$, with the column to the right of the vertical line contained in $P_0$, then $\mathcal{M}_w(\Phi)\nsubseteq\mathcal{EX}(\mathrm{PG}(3,2)\backslash e)$.
\[
\left[
\begin{array}{cc|c}
1&1&1\\
1&1&0\\
1&0&1\\
0&1&1\\
\end{array}
\right],
\left[
\begin{array}{cc|c}
1&1&1\\
1&1&0\\
1&0&0\\
0&1&1\\
0&0&1\\
\end{array}
\right],
\left[
\begin{array}{cc|c}
1&1&1\\
1&1&0\\
1&0&0\\
0&1&0\\
0&0&1\\
0&0&1\\
\end{array}
\right]
\] 
\end{lemma}
\begin{proof} See the Computations 33-35 in the Appendix.\end{proof}

\begin{lemma}\label{A.17}
If $[P_1|P_0]$ contains any of the following matrices as a submatrix, with the portion to the left of the vertical line contained in $P_1$ and the portion to the right contained in $P_0$, then $\mathcal{M}_w(\Phi)\nsubseteq\mathcal{EX}(\mathrm{PG}(3,2)\backslash e)$:
\[
\left[
\begin{array}{c|cc}
1&1&0\\
1&0&1\\
0&1&1\\
0&1&1\\
\end{array}
\right],
\left[
\begin{array}{c|cc}
1&1&0\\
1&0&1\\
0&1&1\\
0&1&0\\
0&0&1\\
\end{array}
\right],
\left[
\begin{array}{c|cc}
1&1&1\\
1&1&0\\
0&1&1\\
0&1&1\\
\end{array}
\right],
\left[
\begin{array}{c|cc}
1&1&1\\
1&1&0\\
0&1&0\\
0&1&1\\
0&0&1\\
\end{array}
\right]
\]\end{lemma}
\begin{proof} See Computations 36-39 in the Appendix.\end{proof}

\begin{lemma}\label{A.18}
If $[P_1|P_0]$ contains either of the matrices below as a submatrix, with the portion to the left of the vertical line contained in $P_1$ and the portion to the right contained in $P_0$, then $\mathcal{M}_w(\Phi)\nsubseteq\mathcal{EX}(\mathrm{PG}(3,2)\backslash e)$:
\[
\left[
\begin{array}{c|cc}
1&1&1\\
1&1&0\\
1&0&1\\
0&1&1\\
\end{array}
\right],
\left[
\begin{array}{c|cc}
1&1&1\\
1&1&0\\
1&0&1\\
0&1&0\\
0&0&1\\
\end{array}
\right]
\]\end{lemma}
\begin{proof} See Computations 40 and 41 in the Appendix.\end{proof}

\begin{lemma}\label{A.19}
If $P_1$ contains either of the following submatrices, then $\mathcal{M}_w(\Phi)\nsubseteq\mathcal{EX}(\mathrm{PG}(3,2)_{-2})$:
\[
\left[
\begin{array}{cc}
1&1\\
1&0\\
0&1\\
\end{array}
\right],
\left[
\begin{array}{cc}
1&1\\
1&1\\
1&0\\
0&1\\
\end{array}
\right]
\]\end{lemma}
\begin{proof} See Computations 42 and 43 in the Appendix.\end{proof}

\begin{lemma}\label{A.20}
If $[P_1|P_0]$ contains either of the following submatrices, with the portion to the left of the vertical line contained in $P_1$ and the portion to the right contained in $P_0$, then $\mathcal{M}_w(\Phi)\nsubseteq\mathcal{EX}(\mathrm{PG}(3,2)_{-2})$.
\[
\left[
\begin{array}{c|c}
1&1\\
1&1\\
0&1\\
0&1\\
\end{array}
\right],
\left[
\begin{array}{c|cc}
1&1\\
1&0\\
0&1\\
0&1\\
\end{array}
\right]
\]\end{lemma}
\begin{proof} See Computations 44 and 45 in the Appendix.\end{proof}

\begin{lemma}\label{A.21}
If $[P_1|P_0]$ contains the submatrix $\left[
\begin{array}{c|c}
1&1\\
1&1\\
1&0\\
0&1\\
\end{array}
\right]$, with the portion to the left of the vertical line contained in $P_1$ and the portion to the right contained in $P_0$, then $\mathcal{M}_w(\Phi)\nsubseteq \mathcal{EX}(\mathrm{PG}(3,2)_{-2})$.\end{lemma}
\begin{proof} See Computation 46 in the Appendix.\end{proof}

\begin{lemma}\label{A.22}
If $P_0$ contains the submatrix $\left[
\begin{array}{ccc}
1&1&1\\
1&0&1\\
1&1&0\\
0&1&1\\
\end{array}
\right]$, then $\mathcal{M}_w(\Phi)\nsubseteq\mathcal{EX}(\mathrm{PG}(3,2)_{-2})$.
\end{lemma}
\begin{proof} See Computation 47 in the Appendix.\end{proof}

\section{Even-Cycle Matroids}
\label{Even-Cycle Matroids}

Theorem ~\ref{proveexminors} shows that the class of even-cycle matroids is contained in $\mathcal{EX}(\mathrm{PG}(3,2)\backslash e, L_{19}, L_{11})$. We will prove Theorem ~\ref{connectedevencycle}, which shows that for sufficiently highly connected binary matroids, the reverse inclusion holds. First, we prove several lemmas. The conclusions of the lemmas in this section are stated up to reordering of the rows and columns of $A_1[X,Y_1]$.

\begin{lemma}
 \label{cecconnected}
Let $\Phi$ be a template with $C=\emptyset$ and with $\Lambda$ trivial. Then at least one of the following holds:
\begin{enumerate}
\item There exists $k\in\mathbb{Z}_+$ such that no simple $k$-connected matroid with at least $2k$ elements virtually conforms to $\Phi$, or
\item $A_1$ is of the following form, with $Y_1=Y_1'\cup Y_1''$ and each $P_i$ an arbitrary binary matrix:
\begin{center}
\begin{tabular}{ |c|c|c| }
\multicolumn{1}{c}{$Y_1'$}&\multicolumn{1}{c}{$Y_1''$}&\multicolumn{1}{c}{$Y_0$}\\
\hline
$I$&$P_1$&$P_0$\\
\hline
\end{tabular}
\end{center}
\end{enumerate}
\end{lemma}

\begin{proof}
By operation (4), we may assume that $A_1$ is of the following form, with $Y_0=V_0\cup V_1$, with $Y_1=Y_1'\cup Y_1''$ and with each $P_i$ an arbitrary binary matrix:

\begin{center}
\begin{tabular}{ |c|c|c|c| }
\multicolumn{1}{c}{$Y_1'$}&\multicolumn{1}{c}{$Y_1''$}&\multicolumn{1}{c}{$V_0$}&\multicolumn{1}{c}{$V_1$}\\
\hline
$I$&$P_1$&0&$P_0$\\
\hline
0&0&$I$&$P_2$\\
\hline
\end{tabular}
.
\end{center}

If $V_0\neq\emptyset$, let $k\geq |Y_0|$, and let $M$ be a matroid virtually conforming to $\Phi$. Then $\lambda(Y_0)=r(Y_0)+r(E(M)-Y_0)-r(M)<r(Y_0)\leq|Y_0|\leq k$. Therefore, $k$-connectivity implies that $|E(M)-Y_0|<|Y_0|$. But then $|E(M)|<2|Y_0|\leq2k$. Thus, (1) holds.

Therefore, we may assume that $V_0=\emptyset$. In this case, (2) holds.
\end{proof}

\begin{lemma}
 \label{cecY1}
Let $\Phi$ be a template with $C=\emptyset$, with $\Delta$ and $\Lambda$ both trivial, and with $A_1$ of the form given in conclusion (2) of Lemma ~\ref{cecconnected}. Then at least one of the following holds:
\begin{enumerate}
 \item $\mathcal{M}_w(\Phi)\nsubseteq\mathcal{EX}(\mathrm{PG}(3,2)\backslash e)$,
\item $P_1$ is a restriction of a matrix of the form $\left[\begin{matrix}
1\cdots1\\
 I\\
\end{matrix}\right]$ and $P_0$ is of the following form, where each column of $[H_1|H_0]$ has at most two nonzero entries:
\[
\left[
\begin{array}{c|c}
1\cdots1&0\cdots0\\
\hline
H_1&H_0\\
\end{array}
\right],
\]
\item $P_1$ is a restriction of a matrix of the form $\begin{bmatrix}
1\cdots1\\
1\cdots1\\
 I\\
\end{bmatrix}$ and $P_0$ is of the following form, where each column of $H_{1,1}$ and $H_{0,0}$ has at most two nonzero entries and where each column of $H_{0,1}$ and $H_{1,0}$ is a unit column or a zero column:
\[
\left[
\begin{array}{c|c|c|c}
1\cdots1&1\cdots1&0\cdots0&0\cdots0\\
1\cdots1&0\cdots0&1\cdots1&0\cdots0\\
\hline
H_{1,1}&H_{1,0}&H_{0,1}&H_{0,0}\\
\end{array}
\right],
\]
\item $P_1=[1,1,0,\dots,0]^T$ or $P_1=[1,1,1,0,\dots,0]^T$, or
\item $Y_1''=\emptyset$.
\end{enumerate}
\end{lemma}

\begin{proof}
 Suppose (1) does not hold. By Lemma ~\ref{simpleY1}, $M(A_1[X,Y_1])$ is simple. By Lemma \ref{A.3}, $M(A_1[X,Y_1])$ has no circuit of size at least 5; thus each column of $P_1$ has at most three nonzero entries. Therefore, every column of $P_1$ has exactly two or three nonzero entries. Let $P_1=[P_{1,2}|P_{1,3}]$, where $P_{1,i}$ consists entirely of columns with exactly $i$ nonzero entries. By Lemmas \ref{A.4} and \ref{A.5}, along with the fact that $M(A_1[X,Y_1])$ is simple, $P_{1,2}$ must be a restriction of a matrix of the form: $\begin{bmatrix}
1\cdots1\\
 I\\
\end{bmatrix}$.

By Lemma \ref{A.4}, every pair of columns in $P_{1,3}$ must both have nonzero entries in the same two rows. Suppose a third column in $P_{1,3}$ does not have nonzero entries in the same two rows as the other two columns. It follows that $P_{1,3}$ must contain the following submatrix:
\[
\left[
\begin{array}{ccc}
1&1&1\\
1&1&0\\
1&0&1\\
0&1&1\\
\end{array}
\right]
\]

However, this contains the submatrix forbidden by Lemma \ref{A.5}. Therefore, $P_{1,3}$ is of the following form: $\begin{bmatrix}
1\cdots1\\
1\cdots1\\
 I\\
\end{bmatrix}$.

We will now show that either $P_{1,2}$ or $P_{1,3}$ is an empty matrix. Consider the matrix
\[
P_1'=\left[
\begin{array}{cc}
1&0\\
1&0\\
1&1\\
0&1\\
\end{array}
\right].
\]
As observed in Lemma \ref{A.7}, this is not a submatrix of $P_1$ because in the matrix $[I_4|P_1']$, columns 1, 2, 4, 5, and 6 form a circuit of size 5. Thus, if $v$ is a column of $P_{1,3}$, then every column $w$ of $P_{1,2}$ must have its nonzero entries in two of the rows that contain the nonzero entries of $v$. This fact, with Lemma \ref{A.5}, implies that if $P_{1,3}$ contains exactly one column, then $P_1$ is a restriction of the following matrix: $\begin{bmatrix}
1&1&1\\
1&0&1\\
0&1&1\\
\end{bmatrix}$, but this is the submatrix forbidden by Lemma \ref{A.6}. Therefore, we see that if $P_{1,3}$ is a nonempty matrix, then $P_1$ must be a restriction of a matrix of the following form:

\[
\left[
\begin{array}{c|c}
1&1\cdots1\\
1&1\cdots1\\
\hline
0&\\
\vdots&I\\
0&
\end{array}
\right]
\]

However, if $P_1$ is of this form and contains the column with two nonzero entries, then by adding the first row to the second and swapping the resulting unit column with the appropriate column of the original identity matrix, we obtain a matrix of the form: $\begin{bmatrix}
1\cdots1\\
 I\\
\end{bmatrix}$. Thus, $P_1$ is either of the form  $\begin{bmatrix}
1\cdots1\\
 I\\
\end{bmatrix}$ or of the form $\begin{bmatrix}
1\cdots1\\
1\cdots1\\
 I\\
\end{bmatrix}$.

Suppose $P_1$ is of the form $\begin{bmatrix}
1\cdots1\\
 I\\
\end{bmatrix}$. If (2) does not hold, then either (4) or (5) holds or $[P_1|P_0]$ contains one of the following submatrices, with the submatrix to the left of the vertical line contained in $P_1$, with the column to the right of the vertical line contained in $P_0$, and with either $x=0$ or $x=1$:
\[
\left[
\begin{array}{cc|c}
1&1&x\\
1&0&1\\
0&1&1\\
0&0&1\\
\end{array}
\right],
\left[
\begin{array}{cc|c}
1&1&x\\
1&0&0\\
0&1&1\\
0&0&1\\
0&0&1\\
\end{array}
\right],
\left[
\begin{array}{cc|c}
1&1&x\\
1&0&0\\
0&1&0\\
0&0&1\\
0&0&1\\
0&0&1\\
\end{array}
\right].
\]
By Lemma \ref{A.15}, (1) holds.

Now, suppose $P_1$ is of the form $\begin{bmatrix}
1\cdots1\\
1\cdots1\\
 I\\
\end{bmatrix}$. If (3) does not hold, then either (4) or (5) holds or $[P_1|P_0]$ contains one of the following submatrices, with the submatrix to the left of the vertical line contained in $P_1$, with the column to the right of the vertical line contained in $P_0$, and with either $x=0$ or $x=1$:
\[
\left[
\begin{array}{cc|c}
1&1&1\\
1&1&0\\
1&0&1\\
0&1&1\\
\end{array}
\right],
\left[
\begin{array}{cc|c}
1&1&1\\
1&1&0\\
1&0&0\\
0&1&1\\
0&0&1\\
\end{array}
\right],
\left[
\begin{array}{cc|c}
1&1&1\\
1&1&0\\
1&0&0\\
0&1&0\\
0&0&1\\
0&0&1\\
\end{array}
\right],
\]

\[
\left[
\begin{array}{cc|c}
1&1&x\\
1&1&x\\
1&0&1\\
0&1&1\\
0&0&1\\
\end{array}
\right],
\left[
\begin{array}{cc|c}
1&1&x\\
1&1&x\\
1&0&0\\
0&1&1\\
0&0&1\\
0&0&1\\
\end{array}
\right],
\left[
\begin{array}{cc|c}
1&1&x\\
1&1&x\\
1&0&0\\
0&1&0\\
0&0&1\\
0&0&1\\
0&0&1\\
\end{array}
\right].
\]
If $[P_1|P_0]$ contains any of the first three of these submatrices, then Lemma \ref{A.16} implies that (1) holds. The last three contain submatrices forbidden by Lemma \ref{A.15}; therefore (1) holds in that case as well.
\end{proof}

\begin{lemma}
 \label{cec1}
Let $\Phi$ be a template with $C=\emptyset$, with $\Delta$ and $\Lambda$ both trivial, with $A_1$ of the form given in conclusion (2) of Lemma ~\ref{cecconnected}, and with $|Y_1''|=1$.  Then at least one of the following holds:
\begin{enumerate}
 \item $\mathcal{M}_w(\Phi)\nsubseteq \mathcal{EX}(\mathrm{PG}(3,2)\backslash e)$,
\item $P_1=[1,1,0\dots,0]^T$ and $P_0$ is of the following form, where no column of $[H_{1,1}|H_{1,0}|H_{0,1}|H_{0,0}]$ has three or more nonzero entries and where at most one of $H_{1,1}$, $H_{1,0}$ and $H_{0,1}$ has columns with two nonzero entries:
\[
\left[
\begin{array}{c|c|c|c}
1\cdots1&1\cdots1&0\cdots0&0\cdots0\\
1\cdots1&0\cdots0&1\cdots1&0\cdots0\\
\hline
H_{1,1}&H_{1,0}&H_{0,1}&H_{0,0}\\
\end{array}
\right], or
\]
\item $P_1=[1,1,1,0\dots,0]^T$ and $P_0$ is a restriction of a matrix of the following form, where each column of $H_{1,1,1}$, $H_{1,0,0}$, $H_{0,1,0}$, and $H_{0,0,1}$ is either a unit column or a zero column, and where each column of $H_{1,1,0}$ and $H_{0,0,0}$ has at most two nonzero entries:
\[
\left[
\begin{array}{c|c|cc|c|c|c|c}
1\cdots1&1\cdots1&1&0&1\cdots1&0\cdots0&0\cdots0&0\cdots0\\
1\cdots1&1\cdots1&0&1&0\cdots0&1\cdots1&0\cdots0&0\cdots0\\
1\cdots1&0\cdots0&1&1&0\cdots0&0\cdots0&1\cdots1&0\cdots0\\
\hline
H_{1,1,1}&H_{1,1,0}&0&0&H_{1,0,0}&H_{0,1,0}&H_{0,0,1}&H_{0,0,0}\\
\end{array}
\right].
\]
\end{enumerate}
\end{lemma}

\begin{proof}
 Suppose (1) does not hold. Since $|Y_1''|=1$, Lemma ~\ref{cecY1} implies that either $P_1=[1,1,0\dots,0]^T$ or $P_1=[1,1,1,0\dots,0]^T$. First, consider the case where $P_1=[1,1,0\dots,0]^T$. If $P_0$ contains either $[0,0,1,1,1]^T$ or $[1,0,1,1,1]^T$ as a submatrix, then by contracting that column of $Y_0$, we obtain in $P_1$ a submatrix forbidden by Lemmas \ref{A.4} or \ref{A.7}, respectively. This, with the fact that no column of $P_0$ can contain five nonzero entries, shows that $[P_1|P_0]$ is of the following form, where each column of $[H_{1,1}|H_{1,0}|H_{0,1}|H_{0,0}]$ contains at most two nonzero entries:
\[
\left[
\begin{array}{c|c|c|c|c}
1&1\cdots1&1\cdots1&0\cdots0&0\cdots0\\
1&1\cdots1&0\cdots0&1\cdots1&0\cdots0\\
\hline
0&H_{1,1}&H_{1,0}&H_{0,1}&H_{0,0}\\
\end{array}
\right].
\]
By Lemmas \ref{A.8} and \ref{A.17}, at most one of $H_{1,1}$, $H_{1,0}$ and $H_{0,1}$ contains a column with two nonzero entries. Therefore, (2) holds.

Now consider the case where $P_1=[1,1,1,0\dots,0]^T$. If $P_0$ contains either $[0,0,0,1,1,1]^T$ or $[1,0,0,1,1]^T$ as a submatrix, then by contracting that column of $Y_0$, we obtain in $P_1$ a submatrix forbidden by Lemmas \ref{A.4} or \ref{A.7}, respectively. This, with the fact that no column of $P_0$ can contain five nonzero entries, shows that $[P_1|P_0]$ is of the following form, where each column of $H_{1,1,1}$, $H_{1,0,0}$, $H_{0,1,0}$, and $H_{0,0,1}$ is either a unit column or a zero column, and where each column of $H_{1,1,0}$, $H_{1,0,1}$, $H_{0,1,1}$, and $H_{0,0,0}$ has at most two nonzero entries
\[
\left[
\begin{array}{c|c|c|c|c|c|c|c}
1\cdots1&1\cdots1&1\cdots1&0\cdots0&1\cdots1&0\cdots0&0\cdots0&0\cdots0\\
1\cdots1&1\cdots1&0\cdots0&1\cdots1&0\cdots0&1\cdots1&0\cdots0&0\cdots0\\
1\cdots1&0\cdots0&1\cdots1&1\cdots1&0\cdots0&0\cdots0&1\cdots1&0\cdots0\\
\hline
H_{1,1,1}&H_{1,1,0}&H_{1,0,1}&H_{0,1,1}&H_{1,0,0}&H_{0,1,0}&H_{0,0,1}&H_{0,0,0}\\
\end{array}
\right].
\]

By Lemma \ref{A.18}, at most one of $H_{1,1,0}$, $H_{1,0,1}$, and $H_{0,1,1}$ contains any nonzero entries. Therefore, if any simple matroids virtually conform to $\Phi$, then (3) holds. Otherwise, we discard $\Phi$ since the results in this paper deal with simple matroids.
\end{proof}

\begin{lemma}
\label{cecY0}Let $\Phi$ be a template with $C=\emptyset$, with $\Delta$ and $\Lambda$ both trivial, and with $A_1$ of the form given in conclusion (2) of Lemma ~\ref{cecconnected}. Then at least one of the following holds:
\begin{enumerate}
 \item $\mathcal{M}_w(\Phi)\nsubseteq \mathcal{EX}(\mathrm{PG}(3,2)\backslash e, L_{11})$,
\item $P_0$ is of the following form, where each column of $[H_1|H_0]$ has at most two nonzero entries:
\[
\left[
\begin{array}{c|c}
1\cdots1&0\cdots0\\
\hline
H_1&H_0\\
\end{array}
\right], or
\]

\item $P_0$ is of the following form, where each column of $H_{1,1}$ and $H_{0,0}$ has at most two nonzero entries and where each column of $H_{0,1}$ and $H_{1,0}$ is a unit column or a zero column:
\[
\left[
\begin{array}{c|c|c|c}
1\cdots1&1\cdots1&0\cdots0&0\cdots0\\
1\cdots1&0\cdots0&1\cdots1&0\cdots0\\
\hline
H_{1,1}&H_{1,0}&H_{0,1}&H_{0,0}\\
\end{array}
\right]
\]
\end{enumerate}
\end{lemma}

\begin{proof}
If $P_0$ satisfies the hypotheses of any one of Lemmas \ref{A.8}-\ref{A.14}, then the result holds. Now suppose $P_0$ satisfies none of the hypotheses of Lemmas \ref{A.8}-\ref{A.14}. This fact, along with the fact that no column of $P_0$ can contain five nonzero entries, implies that $P_0$ has at most one column, other than duplicate columns, with more than two nonzero entries. In this case, clearly (2) holds.
\end{proof}

We now prove Theorem ~\ref{connectedevencycle}.

\begin{proof}[Proof of Theorem ~\ref{connectedevencycle}]
 Let $\mathcal{M}$ denote the class of binary matroids in $\mathcal{EX}(\mathrm{PG}(3,2)\backslash e, L_{19}, L_{11})$ and let $\mathcal{T}=\{\Phi_1,\dots\Phi_s,\Psi_1,\dots,\Psi_t\}$ be the set of templates given by Hypothesis \ref{hyp:connected-template}. Consider a template $\Psi\in\{\Psi_1,\dots,\Psi_t\}$. Recall that every matroid coconforming to $\Psi$ must be contained in the minor-closed class $\mathcal{M}$. Every cographic matroid is a minor of a matroid that coconforms to $\Psi$. Therefore,  $\Psi$ does not exist since $\mathcal{M}$ does not contain $L_{19}$, which is cographic. Thus, $t=0$ and $\mathcal{T}=\{\Phi_1,\dots\Phi_s\}$. Because $\mathrm{PG}(3,2)\backslash e$, $L_{19}$, and $L_{11}$ are simple matroids, it suffices to consider the simple matroids conforming to these templates.

Any matroid containing $\mathrm{PG}(3,2)$ as a minor of course also contains $\mathrm{PG}(3,2)\backslash e$. Therefore, Lemma ~\ref{PG32Phi} implies that, for any template $\Phi\in\{\Phi_1,\dots\Phi_s\}$, either $\Phi\preceq\Phi_X$ or $\Phi$ is a template with $C=\emptyset$ and with $\Lambda$ and $\Delta$ trivial. We will show that in fact $\Phi\preceq\Phi_X$. In this case, we will be able to assume that $\mathcal{T}=\{\Phi_X\}$, since $\mathcal{M}(\Phi_X)$ is the class of even-cycle matroids and is therefore minor-closed.

Suppose, for contradiction, that $\Phi$ is a template with $C=\emptyset$ and with $\Lambda$ and $\Delta$ trivial. Since $\mathcal{M}_w(\Phi)\subseteq \mathcal{EX}(\mathrm{PG}(3,2)\backslash e)$, conclusion (2) of Lemma ~\ref{cecconnected} holds, and one of conclusions (2)-(5) of Lemma ~\ref{cecY1} holds. If conclusion (2) of Lemma ~\ref{cecY1} holds, then any matroid virtually conforming to $\Phi$ is clearly an even-cycle matroid. Similarly, if conclusion (3) holds, then by adding the first row to the second we see that any matroid virtually conforming to $\Phi$ is an even-cycle matroid. Therefore, if either (2) or (3) of Lemma ~\ref{cecY1} holds, then $\Phi\preceq\Phi_X$. Since we already know that $\mathcal{M}(\Phi_X)\subseteq\mathcal{M}$, we may discard $\Phi$ as a template that describes $\mathcal{M}$.

Now suppose conclusion (4) of Lemma ~\ref{cecY1} holds. Then, in Lemma ~\ref{cec1}, either (2) or (3) holds. By adding the first row to the second we see that any matroid virtually conforming to $\Phi$ is an even-cycle matroid. Therefore, we may again discard $\Phi$ as a template that describes $\mathcal{M}$.

Now suppose conclusion (5) of Lemma ~\ref{cecY1} holds; so $Y_1''=\emptyset$. In Lemma ~\ref{cecY0}, either conclusion (2) or (3) holds. By adding a row to another if necessary, we again see that $\Phi\preceq\Phi_X$. Thus, $\Phi$ may be discarded.
\end{proof}

We now prove Theorem ~\ref{blockingpair}.

\begin{proof}[Proof of Theorem ~\ref{blockingpair}]
 Let $\mathcal{M}$ denote the class of binary matroids in $\mathcal{EX}(\mathrm{PG}(3,2)\backslash L,M^*(K_6))$ and let $\mathcal{T}=\{\Phi_1,\dots\Phi_s,\Psi_1,\dots,\Psi_t\}$ be the set of templates for $\mathcal{M}$ given by Hypothesis \ref{hyp:connected-template}. As above, since $M^*(K_6)$ is cographic, $t=0$ and $\mathcal{T}=\{\Phi_1,\dots\Phi_s\}$. Recall that $\mathcal{M}_w(\Phi_{Y_1})$ is the class of even-cycle matroids with a blocking pair. We will show that $\mathcal{T}=\{\Phi_{Y_1}\}$. Because $\mathrm{PG}(3,2)\backslash L$, and $M^*(K_6)$ are simple matroids, it suffices to consider the simple matroids conforming to these templates.

Let $\Phi$ be a template such that $\mathcal{M}_w(\Phi)\subseteq\mathcal{EX}(\mathrm{PG}(3,2)\backslash L, M^*(K_6))$. To show that $\mathrm{PG}(3,2)\backslash L$ is a minor of some matroid virtually conforming to a template, it suffices to show that $\mathrm{PG}(3,2)_{-2}$ is a minor of some matroid virtually conforming to that template. Moreover, for each computation in the Appendix where we showed that $L_{11}$ was a minor of a matroid virtually conforming to some template, we also showed that $\mathrm{PG}(3,2)_{-2}$ is a minor of some matroid virtually conforming to that template. Therefore, $\mathcal{M}_w(\Phi)\subseteq\mathcal{EX}(\mathrm{PG}(3,2)\backslash e, L_{19}, L_{11})\subseteq\mathcal{EX}(\mathrm{PG}(3,2))$.

Thus, by Lemma ~\ref{PG32Phi}, since $\mathrm{PG}(3,2)\backslash L$ conforms to $\Phi_X$, we have that $\Phi$ is a template with $C=\emptyset$ and with $\Lambda$ and $\Delta$ trivial. By Lemma ~\ref{cecconnected}, we may assume that $A_1$ is of the following form, with $Y_1=Y_1'\cup Y_1''$ and each $P_i$ an arbitrary binary matrix:
\begin{center}
\begin{tabular}{ |c|c|c| }
\multicolumn{1}{c}{$Y_1'$}&\multicolumn{1}{c}{$Y_1''$}&\multicolumn{1}{c}{$Y_0$}\\
\hline
$I$&$P_1$&$P_0$\\
\hline
\end{tabular}
\end{center}
By Lemma ~\ref{cecY1} and Lemma \ref{A.19}, either $P_1=[1,1,0,\dots,0]^T$, or $P_1=[1,1,1,0,\dots,0]^T$, or $Y_1''=\emptyset$.

\textit{Case 1:}
Suppose $P_1=[1,1,0,\dots,0]^T$. By Lemma \ref{cec1}, $P_0$ is of the following form, where no column of $[H_{1,1}|H_{1,0}|H_{0,1}|H_{0,0}]$ has three or more nonzero entries:
\[
\left[
\begin{array}{c|c|c|c}
1\cdots1&1\cdots1&0\cdots0&0\cdots0\\
1\cdots1&0\cdots0&1\cdots1&0\cdots0\\
\hline
H_{1,1}&H_{1,0}&H_{0,1}&H_{0,0}\\
\end{array}
\right]\]

By Lemma \ref{A.20}, no column of $[H_{1,1}|H_{1,0}|H_{0,1}]$ has two or more nonzero entries. Let $A$ be any matrix with $r$ rows that virtually conforms to $\Phi$. Add a row $r+1$ to the matrix, where row $r+1$ is the sum of rows $2,\dots, r$. Then one can see that $M(A)$ is an even-cycle matroid with row 1 being the sign row and with a blocking pair represented by rows $2$ and $r+1$. Therefore, $\Phi\preceq\Phi_{Y_1}$, and we may discard $\Phi$.

\textit{Case 2:}
Suppose $P_1=[1,1,1,0,\dots,0]^T$. By Lemma \ref{cec1}, $P_0$ is of the following form, where each column of $H_{1,1,1}$, $H_{1,0,0}$, $H_{0,1,0}$, and $H_{0,0,1}$ is either a unit column or a zero column, and where each column of $H_{1,1,0}$ and $H_{0,0,0}$ has at most two nonzero entries:
\[
\left[
\begin{array}{c|c|cc|c|c|c|c}
1\cdots1&1\cdots1&1&0&1\cdots1&0\cdots0&0\cdots0&0\cdots0\\
1\cdots1&1\cdots1&0&1&0\cdots0&1\cdots1&0\cdots0&0\cdots0\\
1\cdots1&0\cdots0&1&1&0\cdots0&0\cdots0&1\cdots1&0\cdots0\\
\hline
H_{1,1,1}&H_{1,1,0}&0&0&H_{1,0,0}&H_{0,1,0}&H_{0,0,1}&H_{0,0,0}\\
\end{array}
\right]
\]

By Lemma \ref{A.21}, $H_{1,1,0}$ must be a zero matrix. Therefore, by adding the first row to the second, we see that any matroid virtually conforming to $\Phi$ is even-cycle with a blocking pair by taking the first row to be the sign row and taking the blocking pair to be represented by the second and third rows. Therefore, $\Phi\preceq\Phi_{Y_1}$, and we may discard $\Phi$.

\textit{Case 3:}
Suppose $Y_1''=\emptyset$. Either conclusion (2) or conclusion (3) of Lemma ~\ref{cecY0} must hold. First, suppose (2) of ~\ref{cecY0} holds. Then $P_0$ is of the following form, where each column of $[H_1|H_0]$ has at most two nonzero entries:
\[
\left[
\begin{array}{c|c}
1\cdots1&0\cdots0\\
\hline
H_1&H_0\\
\end{array}
\right]
\]
We may assume that no column of $H_1$ is a zero column because any such columns can be obtained from $Y_1$. Therefore, by Lemma \ref{A.22}, $H_1$ is a restriction of a matrix of the form $\begin{bmatrix}
1\cdots1\\
 I\\
\end{bmatrix}$. By a similar argument as was used in Case 1, we see that every matroid virtually conforming to $\Phi$ is even-cycle with a blocking pair. Therefore, $\Phi\preceq\Phi_{Y_1}$, and we may discard $\Phi$.

Now, suppose (2) of ~\ref{cecY0} holds. Then $P_0$ is of the following form, where each column of $H_{1,1}$ and $H_{0,0}$ has at most two nonzero entries and where each column of $H_{0,1}$ and $H_{1,0}$ is a unit column or a zero column:
\[
\left[
\begin{array}{c|c|c|c}
1\cdots1&1\cdots1&0\cdots0&0\cdots0\\
1\cdots1&0\cdots0&1\cdots1&0\cdots0\\
\hline
H_{1,1}&H_{1,0}&H_{0,1}&H_{0,0}\\
\end{array}
\right]
\]
By Lemma \ref{A.22}, $H_{1,1}$ is a restriction of a matrix of the form $\begin{bmatrix}
1\cdots1\\
 I\\
\end{bmatrix}$, with possibly a zero column as well. Therefore, by adding the first row to the second, we see that any matroid virtually conforming to $\Phi$ is even-cycle with a blocking pair by taking the first row to be the sign row and taking the blocking pair to be represented by the second and third rows. Therefore, $\Phi\preceq\Phi_{Y_1}$, and we may discard $\Phi$. This completes the proof.
\end{proof}

\section{Some Technical Lemmas Proved with SageMath: Even-Cut Matroids}
\label{sec:technical-cut}
In this section, we will list several technical lemmas that we will need to prove Theorems \ref{connectedevencut}. As was the case with Section \ref{sec:technical-cycle}, many of the proofs will merely refer the reader to a computation in the Appendix; the computations will use the SageMath software system. The reader may prefer to move on to Section \ref{Even-Cut Matroids}, referring to Section \ref{sec:technical-cut} as necessary. In Lemmas \ref{B.3}-\ref{B.10}, $\Psi$ is a template with $C=\emptyset$, with $\Lambda$ trivial, and with $\Delta=\{0,\bar{x}\}$ for some row vector $\bar{x}$. Moreover, there are partitions $Y_1=Y'_1\cup Y_{1,0}\cup Y_{1,1}$ and $Y_0=Y_{0,0}\cup Y_{0,1}$ such that $\bar{x}_y=1$ if and only if $y\in Y_{1,1}\cup Y_{0,1}$ and such that $A_1$ is of the following form:
\begin{center}
\begin{tabular}{ |c|c|c|c|c| }
\multicolumn{1}{c}{$Y'_1$}&\multicolumn{1}{c}{$Y_{1,0}$}&\multicolumn{1}{c}{$Y_{1,1}$}&\multicolumn{1}{c}{$Y_{0,0}$}&\multicolumn{1}{c}{$Y_{0,1}$}\\
\hline
$I$&$A_{Y_1}$&$B_{Y_1}$&$A_{Y_0}$&$B_{Y_0}$\\
\hline
\end{tabular}
.
\end{center}

Computation 48 in the Appendix gives a SageMath function which builds the largest possible matrix $A$ conforming to such a template whose vector matroid is a simple matroid of rank $r+|X|$. The variable \texttt{B\_Jrows} specifies the number of row indices $b\in B$ for which $A[b,Y_0\cup Y_1]=\bar{x}$.

\begin{lemma}\label{B.3}
 If $B_{Y_1}$ contains the submatrix $[1,0]$, then $\mathcal{M}_w(\Psi)\nsubseteq\mathcal{EX}(H_{12})$.\end{lemma}
\begin{proof} See Computation 49 in the Appendix.\end{proof}

\begin{lemma}\label{B.4}
If $[A_{Y_1}|B_{Y_1}]$ contains either of the following submatrices, with the column to the left of the vertical line contained in $A_{Y_1}$, and the column to the right of the vertical line contained in $B_{Y_1}$, then $\mathcal{M}_w(\Psi)\nsubseteq\mathcal{EX}(H_{12})$.
\[
\left[
\begin{array}{c|c}
1&1\\
1&1\\
\end{array}
\right],
\left[
\begin{array}{c|c}
1&0\\
1&0\\
\end{array}
\right],
\left[
\begin{array}{c|c}
1&1\\
1&0\\
\end{array}
\right]
\]
\end{lemma}
\begin{proof} See Computations 50-52 in the Appendix.\end{proof}

\begin{lemma}\label{B.5}
If $A_{Y_1}$ contains either of the following submatrices, then $\mathcal{M}_w(\Psi)\nsubseteq\mathcal{EX}(H_{12})$.
\[
\left[
\begin{array}{cc}
1&1\\
1&1\\
0&1\\
\end{array}
\right],
\left[
\begin{array}{cc}
1&0\\
1&1\\
0&1\\
\end{array}
\right],
\left[
\begin{array}{cc}
1&0\\
1&0\\
0&1\\
0&1\\
\end{array}
\right]
\]
\end{lemma}
\begin{proof} See Computations 53-55 in the Appendix.\end{proof}

\begin{lemma}\label{old4.6.1}
If $A_1$ contains either of the following matrices, with the column on the left indexed by an element of $Y_1$ and the column on the right is indexed by an element of $Y_0$, then $\mathcal{M}_w(\Psi)\nsubseteq\mathcal{EX}(H_{12})$.
\[
\left[
\begin{array}{c|c}
1&1\\
1&1\\
0&1\\
0&1\\
\end{array}
\right],
\left[
\begin{array}{c|c}
1&0\\
1&1\\
0&1\\
0&1\\
\end{array}
\right],
\left[
\begin{array}{c|c}
1&0\\
1&0\\
0&1\\
0&1\\
0&1\\
\end{array}
\right],
\left[
\begin{array}{c|c}
1&1\\
1&1\\
1&0\\
0&1\\
\end{array}
\right],
\left[
\begin{array}{c|c}
1&0\\
1&0\\
1&1\\
1&1\\
1&1\\
\end{array}
\right]
\]
\end{lemma}
\begin{proof} If the column on the left is contained in $A_{Y_1}$ and the column on the right is contained in $A_{Y_0}$, then these matrices are forbidden because contraction of the element indexing this column of $A_{Y_0}$ produces a new $A_{Y_1}$, containing a column originally in the identity matrix, that contains one of the  submatrices listed in Lemma \ref{B.5}. Since we may choose the zero vector for every element of $\Delta$, we also have $\mathcal{M}_w(\Psi)\nsubseteq\mathcal{EX}(H_{12})$ if these submatrices are contained in $[A_{Y_1}|B_{Y_0}]$, in $[B_{Y_1}|B_{Y_0}]$, or in $[B_{Y_1}|A_{Y_0}]$.\end{proof}

\begin{lemma}\label{B.6}
If $[A_{Y_1}|B_{Y_0}]$ contains either of the following submatrices, with the column to the left of the vertical line contained in $A_{Y_1}$, and the column to the right of the vertical line contained in $B_{Y_0}$, then $\mathcal{M}_w(\Psi)\nsubseteq\mathcal{EX}(H_{12})$.
\[
\left[
\begin{array}{c|c}
1&0\\
1&0\\
\end{array}
\right],
\left[
\begin{array}{c|c}
1&1\\
1&0\\
\end{array}
\right],
\left[
\begin{array}{c|c}
1&1\\
1&1\\
\end{array}
\right]
\]\end{lemma}
\begin{proof} See Computations 56-58 in the Appendix.\end{proof}

\begin{lemma}\label{B.7}
If $[B_{Y_1}|A_{Y_0}]$ contains either of the following submatrices, with the column to the left of the vertical line contained in $B_{Y_1}$, and the column to the right of the vertical line contained in $A_{Y_0}$, then $\mathcal{M}_w(\Psi)\nsubseteq\mathcal{EX}(H_{12})$.
\[
\left[
\begin{array}{c|c}
0&1\\
0&1\\
0&1
\end{array}
\right],
\left[
\begin{array}{c|c}
1&1\\
0&1\\
0&1\\
\end{array}
\right],
\left[
\begin{array}{c|c}
1&1\\
1&1\\
0&1
\end{array}
\right],
\left[
\begin{array}{c|c}
1&1\\
1&1\\
1&1
\end{array}
\right]
\]\end{lemma}
\begin{proof} See Computations 59-62 in the Appendix.\end{proof}

\begin{lemma}\label{B.8}
 If $[B_{Y_1}|B_{Y_0}]$ contains  any of the following submatrices, with the column to the left of the vertical line contained in $B_{Y_1}$, and the column to the right of the vertical line contained in $B_{Y_0}$, then $\mathcal{M}_w(\Psi)\nsubseteq\mathcal{EX}(H_{12})$.
\[
\left[
\begin{array}{c|c}
0&1\\
0&1\\
\end{array}
\right],
\left[
\begin{array}{c|c}
1&0\\
0&1\\
\end{array}
\right],
\left[
\begin{array}{c|c}
1&0\\
1&0\\
\end{array}
\right]
\]
\end{lemma}
\begin{proof} See Computations 63-65 in the Appendix.\end{proof}

\begin{lemma}\label{B.9}
 If $A_{Y_0}$ contains the following submatrix, then $\mathcal{M}_w(\Psi)\nsubseteq\mathcal{EX}(H_{12})$.
\[
\left[
\begin{array}{cccc}
1&0&1&0\\
1&1&0&0\\
1&1&1&1\\
1&1&1&1\\
1&1&1&1\\
0&1&1&0\\
\end{array}
\right]
\]
\end{lemma}
\begin{proof} See Computation 66 in the Appendix.\end{proof}

\begin{lemma}\label{B.10}
 If $B_{Y_0}$ contains the following submatrix, then $\mathcal{M}_w(\Psi)\nsubseteq\mathcal{EX}(H_{12})$.
\[
\left[
\begin{array}{cccc}
1&1&0&0\\
1&0&1&0\\
\end{array}
\right]
\]
\end{lemma}
\begin{proof} See Computation 67 in the Appendix.\end{proof}

\begin{lemma}\label{new}
 If $A_{Y_0}$ contains the following submatrix, then $\mathcal{M}_w(\Psi)\nsubseteq\mathcal{EX}(H_{12})$.
\[
\left[
\begin{array}{cccc}
1&1&1&1\\
1&1&1&1\\
1&1&0&0\\
1&0&1&0\\
1&0&0&1\\
\end{array}
\right]
\]
\end{lemma}
\begin{proof} See Computation 68 in the Appendix.\end{proof}

\section{Even-Cut Matroids}
\label{Even-Cut Matroids}
In this section, we prove Theorem ~\ref{connectedevencut}. Recall that we use the following definition: An \textit{even-cut matroid} is a matroid $M$ that can be represented by a binary matrix with a row whose removal results in a matrix representing a cographic matroid. Thus, there is some binary extension $N$ of $M$ on ground set $E(M)\cup\{e\}$ such that $N/e$ is cographic. Thus, to check if a binary matroid $M$ is even-cut, it suffices to check if $N/e$ is cographic for some binary extension $N$ of $M$. It will be useful to consider the dual situation. Therefore, it suffices to check if there is a binary coextension $N^*$ of $M^*$ such that $N^*\backslash e$ is graphic. If this is the case, then $M^*\in\mathcal{M}(\Phi_C)$. We see then that $\mathcal{M}^*(\Phi_C)$ is exactly the class of even-cut matroids. Recall from Lemma ~\ref{YCD} that $\Phi_{Y_1}\preceq\Phi_C$. This property reflects the fact, first observed by Pivotto ~\cite{p11}, that even-cycle matroids with a blocking pair are duals of even-cut matroids.

Recall that $H_{12}$ is the matroid with the even-cycle representation given in Figure \ref{fig:H12}. Thus, $H_{12}$ is the vector matroid of the binary matrix below. In that matrix, the top row is the sign row.

\[
\left[
\begin{array}{cccccccccccc}
1&0&1&0&1&0&1&0&1&0&1&0\\
0&0&0&0&1&1&1&1&0&0&0&0\\
0&0&0&0&1&1&0&0&1&1&0&0\\
1&1&0&0&0&0&0&0&0&0&1&1\\
0&0&1&1&0&0&1&1&0&0&1&1\\
\end{array}
\right]
\]

Lemma ~\ref{evencutexcludedminors} shows that the class of even-cut matroids is contained in $\mathcal{EX}(M(K_6), H^*_{12})$. Theorem ~\ref{connectedevencut} is the claim that for sufficiently highly connected matroids, the reverse inclusion holds. We will prove Theorem ~\ref{connectedevencut} after giving a definition and proving some lemmas.

\begin{definition}
 Let $|C|=2$ and let $\Delta$ be the subgroup of $\mathrm{GF}(2)^C$ generated by $[1,0]$ and $[0,1]$. The template $\Phi^2_C$ is given by
\[\Phi^2_C=(C,\emptyset,\emptyset,\emptyset,[\emptyset],\Delta,\{0\}).\]
\end{definition}

\begin{lemma}
\label{PhiC2}
 For a template $\Phi$, either $\Phi^2_C\preceq\Phi$ or $\Phi$ is equivalent to a template with $|C_1|\leq1$, where $C_1$ is as in Figure \ref{fig:A'}.
\end{lemma}

\begin{proof}
There are three cases to consider.

\textit{Case 1:} Every element of $\Delta|C$ is in the row space of $A_1[X_0,C]$. Then contraction of $C_0$ turns the elements of $C_1$ into loops, and contraction of $C_1$ is the same as deletion of $C_1$. By deleting $C_1$ from every matrix virtually conforming to $\Phi$, we see that $\Phi$ is equivalent to a template with $C_1=\emptyset$.

\textit{Case 2:}  There is exactly one element $\bar{x}\in\Delta|C$ that is not in the row space of $A_1[X_0,C]$. Then contraction of $C_0$ turns the elements of $C_1$ into parallel elements. Thus, contraction of some element $c\in C_1$ turns the elements of $C_1-\{c\}$ into loops, and contraction of $C_1-\{c\}$ is the same as deletion of $C_1-\{c\}$. By deleting $C_1-\{c\}$ from every matrix virtually conforming to $\Phi$, we see that $\Phi$ is equivalent to a template with $|C_1|=1$.

\textit{Case 3:}  There are distinct elements $\bar{x}$ and $\bar{y}$ in $\Delta|C$ that are not in the row space of $A_1[X_0,C]$. Index the elements of $C_0$ by $\{1,2,\dots,n\}$ and the elements of $X_0$ by $\{d_1,d_2,\dots,d_n\}$. Let $S_x$ and $S_y$ be the supports of $\bar{x}|C_0$ and $\bar{y}|C_0$, respectively. Then the support of $(\bar{x}+\bar{y})|C_0$ is the symmetric difference $S_x\triangle S_y$. First, suppose that for every pair of elements $\bar{x}$ and $\bar{y}$ in $\Delta|C$ that are not in the row space of $A_1[X_0,C]$, we have that $\bar{x}+\bar{y}$ is in the row space of $A_1[X,C]$. Since the rows of $A_1[X_0,C]$ are linearly independent, it must be that the zero vector is equal to
\[\sum_{i\in S_x\triangle S_y}A_1[\{d_i\},C]+\bar{x}+\bar{y}=\sum_{i\in S_x}A_1[\{d_i\},C]+\bar{x}+\sum_{i\in S_y}A_1[\{d_i\},C]+\bar{y}\]
and therefore, since we are working in characteristic 2,
\[\sum_{i\in S_x}A_1[\{d_i\},C]+\bar{x}=\sum_{i\in S_y}A_1[\{d_i\},C]+\bar{y}.\]
Thus, contraction of $C_0$ projects $\bar{x}$ and $\bar{y}$ onto the same element of $\mathrm{GF}(2)^{C_1}$.
Moreover, this is true for any pair of elements of $\Delta|C$ that are not in the row space of $A_1[X_0,C]$. Therefore, the same argument used for Case 2 shows that $\Phi$ is equivalent to a template with $|C_1|=1$.

Therefore, we may assume that there are elements $\bar{x}$ and $\bar{y}$ in $\Delta|C$ that are not in the row space of $A_1[X_0,C]$ and such that $\bar{x}+\bar{y}$ is also not in the row space of $A_1[X_0,C]$. Repeatedly perform operations (3) and (9) on $\Phi$ until the following template is obtained:
\[(C,X,\emptyset,\emptyset,A_1[X,C],\Delta|C,\Lambda).\]
On this template, perform operations (1) and (2) to obtain the following template:
\[(C,X,\emptyset,\emptyset,A_1[X,C],\langle\bar{x},\bar{y}\rangle,\{0\}).\]
By performing elementary row operations, we see that every matrix virtually respecting this template is row equivalent to a matrix virtually respecting the following template, where $\bar{x}'|C_0$ and $\bar{y}'|C_0$ are zero vectors:
\[\Phi'=(C,X,\emptyset,\emptyset,A_1[X,C],\langle\bar{x}',\bar{y}'\rangle,\{0\}).\]
Note that $\bar{x}'|C_1$, $\bar{y}'|C_1$, and $(\bar{x}'+\bar{y}')|C_1$ are nonzero since $\bar{x}$, $\bar{y}$, and $\bar{x}+\bar{y}$ were not in the row space of $A_1[X_0,C]$ in the original template $\Phi$. Also, we must have $\bar{x}'\neq\bar{y}'$ because otherwise, $\bar{x}'+\bar{y}'=\bf{0}$, contradicting the assumption that $\bar{x}+\bar{y}$ was not in the row space of $A_1[X_0,C]$ in $\Phi$. Now, on $\Phi'$, repeatedly perform operation (6) and then operation (5) to obtain the following template:
\[\Phi''=(C_1,\emptyset,\emptyset,\emptyset,[\emptyset],\langle\bar{x}'|C_1,\bar{y}'|C_1\rangle,\{0\}).\]

Now, every matroid $M$ conforming to $\Phi''$ is obtained by contracting $C_1$ from $M(A)$, where $A$ is a matrix conforming to $\Phi''$. Thus, if there are any elements of $C_1$ that are parallel elements in $M(A)$, contracting one of these elements turns the rest of the parallel class into loops. So these elements are deleted to obtain $M$. Thus, $\Phi''$ is equivalent to a template where these elements have been deleted from $C$. There are two cases to consider. First, if it is the case that either the supports of $\bar{x}'$ and $\bar{y}'$ are disjoint or that one support is contained in the other, then in the resulting template, $|C|=2$ and $\Delta=\langle[1,0],[0,1]\rangle$. So this resulting template is $\Phi^2_C$. In the other case, $\bar{x}'$ and $\bar{y}'$ have intersecting supports but neither is contained in the other. In this case, $\Phi''$ is equivalent to the following template with $|C_1|=3$:
\[\Phi'''=(C_1,\emptyset,\emptyset,\emptyset,[\emptyset],\langle[1,1,0],[1,0,1]\rangle,\{0\}).\]
However, by contracting any element of $C$, the other two become parallel. Thus, by contracting a second element, the third becomes a loop. Therefore, the third element is deleted to obtain a matroid conforming to $\Phi'''$. Thus, $\Phi^2_C\sim\Phi'''\preceq\Phi$.
\end{proof}

\begin{lemma}
\label{PhiCD2}
 If $\Phi$ is a template with $|C_1|=1$ and with $\Lambda|X_1$ trivial, then $\Phi_{CX}\preceq\Phi$ or $\Phi$ is equivalent to a template with $C=\emptyset$.
\end{lemma}

\begin{proof}
We consider two cases, depending on whether $\Delta|C$ contains an element that is not in the row space of $A_1[X_0,C]$.

\textit{Case 1:} Every element of $\Delta|C$ is in the row space of $A_1[X_0,C]$. Let $A$ be a matrix that conforms to $\Phi$. When $C_0$ is contracted from $M(A)$, each element of $C_1$ becomes a loop and can therefore be deleted rather than contracted. Thus, $\Phi$ is equivalent to a template $\Phi'$ with $C_1=\emptyset$. Suppose there exist elements $\bar{x}\in\Delta|C_0$ and $\bar{y}\in\Lambda|X_0$ such that there are an odd number of natural numbers $i$ with $\bar{x}_i=\bar{y}_i=1$. Repeatedly perform operations (3), (9), and (5) on $\Phi'$ to obtain the following template:
\[(C_0,X_0,\emptyset,\emptyset,A_1[X_0,C_0],\Delta|C_0,\Lambda).\]
Then perform operations (1) and (2) to obtain the following template:
\[(C_0,X_0,\emptyset,\emptyset,A_1[X_0,C_0],\{\mathbf{0}, \bar{x}\},\{\mathbf{0}, \bar{y}\}).\]
Any matroid conforming to this template is obtained by contracting $C_0$ from $M(A)$, where $A$ is a matrix conforming to $\Phi$. Recall that $A[B-X_0,E-C_0]$ is a frame matrix. If $\bar{x}$ is in the row labeled by $r$ and $\bar{y}$ is in the column labeled by $c$, then when $C_0$ is contracted, 1 is added to the entry of the frame matrix in row $r$ and column $c$. Otherwise, the entry remains unchanged when $C$ is contracted. We see then that this template is equivalent to $\Phi_{CX}$, where 1s are used to replace $\bar{x}$ and $\bar{y}$.

Thus, we may assume that for every element  $\bar{x}\in\Delta|C_0$ and $\bar{y}\in\Lambda|X_0$, there are an even number of natural numbers $i$ such that $\bar{x}_i=\bar{y}_i=1$. This implies that contraction of $C$ has no effect on the frame matrix. So $\Phi'$, and therefore $\Phi$ are equivalent to a template with $\Lambda|X_0$ trivial. In this case, we see that repeated use of operation (6) produces a template equivalent to $\Phi$ with $C=\emptyset$.

\textit{Case 2:} There is an element $\delta\in\Delta|C$ that is not in the row space of $A_1[X_0,C]$. Since $|C_1|=1$, every element of $\Delta|C$ not in the row space of $A_1[X_0,C]$ becomes a 1 after $C_0$ is contracted, and every element that is in the row space becomes a 0. Therefore, we may assume that the column vector $A_1[X,C_1]$ is a zero vector and that an element of $\Delta|C$ has a 0 as its final entry if it is in the row space and a 1 otherwise.

If $\bar{x}\in\Delta|C$ and $\bar{y}\in\Lambda|X_0$ are such that there are an odd number of natural numbers $i$ such that $\bar{x}_i=\bar{y}_i=1$, then we call the ordered pair $(\bar{x},\bar{y})$ a \textit{pair of odd type}. Otherwise, $(\bar{x},\bar{y})$ is a \textit{pair of even type}.  Consider a matrix $A$ virtually conforming to $\Phi$ and contract $C$ from $M(A)$. The effect on the elements of $\Delta$ is a change of basis followed by a projection into a lower dimension. Therefore, a group structure is maintained. Let us call the resulting group $\Delta'$. There are two subcases to check.

\textit{Subcase a:} Suppose there exists a pair $(\bar{x},\bar{y})$ of odd type. If $\bar{x}$ is in the row space of $A_1[X_0,C]$, or if $\bar{x}$ is not in the row space of $A_1[X_0,C]$ but $(\delta,\bar{y})$ is a pair of even type, then we will show that $\Phi_{CX}\preceq\Phi$. On $\Phi$, repeatedly perform operations (3), (9), and (5) as needed to obtain the following template:
\[(C,X_0,\emptyset,\emptyset,A_1[X_0,C],\Delta|C,\Lambda).\]
Then perform operations (1) and (2) to obtain the following template:
\[\Phi'=(C,X_0,\emptyset,\emptyset,A_1[X_0,C],\langle\bar{x},\delta\rangle,\{\mathbf{0}, \bar{y}\}).\]
Consider the following matrix conforming to $\Phi_{CX}$:

\begin{center}
\begin{tabular}{|c|c|c|}
\hline
$0\cdots0$&$1\cdots1$&1\\
\hline
&&1\\
&&$\vdots$\\
frame&frame&1\\
matrix&matrix&0\\
&&$\vdots$\\
&&0\\
\hline
\end{tabular}
\end{center}
The matrix below conforms to $\Phi'$ and results in the same matroid when $C$ is contracted:
\begin{center}
\begin{tabular}{|c|c|c|c|}
\multicolumn{2}{c}{}&\multicolumn{1}{c}{$C_0$}&\multicolumn{1}{c}{$C_1$}\\
\hline
&&&0\\
0&$\bar{y}\cdots\bar{y}$&$I$&$\vdots$\\
&&&0\\
\hline
$0\cdots0$&$0\cdots0$&\multicolumn{2}{c|}{$\delta$}\\
\hline
&&\multicolumn{2}{c|}{$\bar{x}$}\\
&&\multicolumn{2}{c|}{$\vdots$}\\
frame&frame&\multicolumn{2}{c|}{$\bar{x}$}\\
\cline{3-4}
matrix&matrix&\multicolumn{2}{c|}{}\\
&&\multicolumn{2}{c|}{0}\\
&&\multicolumn{2}{c|}{}\\
\hline
\end{tabular}
\end{center}

Therefore, we may assume that an element $\delta'\in\Delta|C$ is in the row space of $A_1[X_0,C]$ if and only if $(\delta',\bar{y})$ is a pair of even type. Moreover, this is true for every nonzero element of $\Lambda|X_0$. Thus, if $\bar{y_1}$ and $\bar{y_2}$ are nonzero elements of $\Lambda|X_0$, then both $(\delta,\bar{y_1})$ and $(\delta,\bar{y_2})$ are pairs of odd type, since $\delta$ is not in the row space of $A_1[X_0,C]$. This implies that $(\delta,\bar{y_1}+\bar{y_2})$ is a pair of even type. But we have just shown that this implies that $\bar{y_1}+\bar{y_2}$ is the zero vector. Thus $\bar{y_1}=\bar{y_2}$ and $\Lambda|X_0=\{\mathbf{0}, \bar{y}\}$ in the original template $\Phi$. By a similar argument, $\bar{x}=\delta$ and $\Delta|C=\{\mathbf{0}, \delta\}$ in the original template $\Phi$. Therefore, $\Phi$ is equivalent to a template with $|C_0|=|X_0|=1$, with $A_1[X_0,C_0]=[1]$, with $A_1[X_0,C_1]=[0]$, with $\Lambda|X_0=\{[0],[1]\}$, and with $\Delta|C=\{[0,0],[1,1]\}$. We may now assume that the original template $\Phi$ was of this form and will study this template $\Phi$ in this form but before any operations have been performed on it.

We will show that $\Phi$ is equivalent to the following template $\Phi'$, with $C=\emptyset$, obtained by adjoining an element $y$ to $Y_1$ and letting $A_1[X_1,y]$ be the zero vector. We will define $\Delta''$ below.
\[\Phi'=(\emptyset,X_1,Y_0,Y_1\cup y,A_1[X_1,Y_0\cup Y_1\cup y],\Delta'',\{0\})\]
Recall that $\Delta'$ is the group obtained from $\Delta$ after $C$ is contracted. Let $\Delta''$ be the subgroup of $\mathrm{GF}(2)^{Y_0\cup Y_1\cup y}$ consisting of all the row vectors obtained by adjoining to any element of $\Delta'$ either a zero or a 1. So $|\Delta''|=2|\Delta'|$. Let $A$ be a matrix that virtually conforms to $\Phi$. Recall that columns indexed by elements of $Z$ are formed by adding a column indexed by $Y_1$ to a column indexed by $Z$ in a matrix that respects $\Phi$. If $A[B-X,C]$ is the zero matrix, then $M(A)/C$ conforms to $\Phi'$ because we may simply choose never to use $y$ to build a column indexed by $Z$.

Otherwise, choose an element $r$ of $B-X$ such that $A[r,C]=[1,1]$. Let $S$ be the subset of $B-(X\cup r)$ such that $s\in S$ if and only if $A[s,C]=[1,1]$. Let $T=B-(X\cup S\cup r)$. The effect on the frame matrix of contracting $C$ from $M(A)$ is to remove $r$ and to add a 1 to each entry $A_{s,c}$ of the frame matrix where $s\in S$ and where $c$ is an element of $E-(C\cup Y_0\cup Y_1\cup Z)$ with $A_{r,c}=1$. Let $\hat{A}$ be the matrix that results from $A$ by contracting $C$. Recall that every column of the frame matrix $A[B-X,E-(C\cup Y_0\cup Y_1\cup Z)]$ contains at most two nonzero entries. Thus, for a column $c$ with $A_{r,c}=1$, the column $A[B-(X\cup r),c]$ must be either a unit column or a zero column. Therefore, there are several possibilities for $\hat{A}[B-(X\cup r),c]$. Either $\hat{A}[S,c]=[1,\dots ,1]^T$ and $\hat{A}[T,c]=[0,\dots ,0]^T$, or $\hat{A}[S,c]=[1,\dots ,1]^T$ and $\hat{A}[T,c]$ is an identity column, or $\hat{A}[S,c]$ is the complement of an identity column and $\hat{A}[T,c]=[0,\dots ,0]^T$. This exact same situation can be obtained with $\Phi'$ using the new column $y$. Thus, $\Phi$ is equivalent to $\Phi'$, a template with $C=\emptyset$.

\textit{Subcase b:} Therefore, we may assume that every pair of elements $(\bar{x},\bar{y})\in\Delta\times(\Lambda|X_0)$ is a pair of even type. Thus, contraction of $C_0$ has no effect on the frame matrix. This implies that $\Phi$ is equivalent to a template with $\Lambda$ trivial. By repeated use of operation (6), we obtain a template equivalent to $\Phi$ with $C_0=\emptyset$, with $|C_1|=1$, and with $\Delta|C=\{[0],[1]\}$. Using an argument similar to the one used at the end of Subcase a, we see that $\Phi$ is equivalent to a template with $C=\emptyset$ by adjoining an element to $Y_1$.
\end{proof}

\begin{lemma}
 \label{YC}
Let $\Phi=(C,X,Y_0,Y_1,A_1,\Delta,\Lambda)$ be a template. If $\Phi'=(C',X,\emptyset,\emptyset,A_1,\Delta,\Lambda)$, where $C'=Y_0\cup Y_1\cup C$, then every matroid conforming to $\Phi'$ is a minor of a matroid conforming to $\Phi$.
\end{lemma}

\begin{proof}
 Let $\Phi''=(C,X,Y'_0,\emptyset,A_1,\Delta,\Lambda)$, where $Y'_0=Y_0\cup Y_1$. By Lemma \ref{yshift}, $\Phi''\preceq\Phi$. Any matroid conforming to $\Phi'$ is obtained from a matroid conforming to $\Phi''$ by contracting $Y'_0$.
\end{proof}

We now prove Theorem ~\ref{connectedevencut}.

\begin{proof}[Proof of Theorem ~\ref{connectedevencut}]
Let $\mathcal{M}=\mathcal{EX}(M(K_6),H^*_{12})$, and let $\mathcal{T}=\{\Phi_1$, $\dots$,$\Phi_s$, $\Psi_1$,$\dots$,$\Psi_t\}$ be the set of templates given by Hypothesis \ref{hyp:connected-template} for $\mathcal{M}$. Consider a template $\Phi\in\{\Phi_1,\dots,\Phi_s\}$. Recall that every matroid conforming to $\Phi$ must be contained in the minor-closed class $\mathcal{M}$. Every graphic matroid is a minor of a matroid that conforms to $\Phi$. Since $\mathcal{M}$ does not contain $M(K_6)$, it must be the case that $\Phi$ does not exist. Thus, $s=0$ and $\mathcal{T}=\{\Psi_1,\dots,\Psi_t\}$. Therefore, we will study the highly connected matroids in $\mathcal{M}$ by considering their dual matroids which virtually conform to some template $\Psi\in\{\Psi_1,\dots,\Psi_t\}$. Because $M(K_6)$, and $H^*_{12}$ are cosimple matroids, it suffices to consider cosimple matroids in  $\mathcal{M}$. Thus, it suffices to consider simple matroids that are duals of matroids in $\mathcal{M}$. Therefore, we only consider simple matroids conforming to $\Psi$.

Let $\Psi=(C,X,Y_0,Y_1,A_1,\Delta,\Lambda)$ be a template in $\mathcal{T}$. We know that $\mathcal{M}^*(\Phi_C)$ is the class of even-cut matroids. Therefore,  we may assume that $\Phi_C\in\{\Psi_1,\dots,\Psi_t\}$, and if any template $\Psi\preceq\Phi_C$, we may discard $\Psi$ form the set $\{\Psi_1,\dots,\Psi_t\}$. Since $H_{12}$ is an even-cycle matroid, $H_{12}$ conforms to $\Phi_X$. Thus, we have $\Phi_X\npreceq\Psi$. By Lemma ~\ref{PhiD}, $\Lambda|X_1$ is trivial. Moreover, by Lemma ~\ref{YCD}, we have $\Phi_{CX}\npreceq\Psi$.

The following matrix conforms to $\Phi^2_C$, with $C$ indexing the last two columns:
\[
\left[
\begin{array}{cccccccccccc|cc}
0&0&0&1&1&0&0&0&0&0&0&0&1&0\\
0&0&0&0&0&0&0&0&0&0&0&1&0&1\\
0&0&1&0&1&0&1&0&1&0&1&0&0&1\\
0&0&0&0&0&1&1&0&0&0&0&0&1&0\\
0&0&0&0&0&0&0&1&1&0&0&0&1&0\\
1&0&0&0&0&0&0&0&0&1&1&0&0&1\\
0&1&0&1&0&1&0&1&0&1&0&0&0&1\\
\end{array}
\right].
\]
By contracting $C$, we obtain the following matrix $A$ with $M(A)$ conforming to $\Phi^2_C$:
\[
A=\left[
\begin{array}{cccccccccccc}
0&0&1&0&1&0&1&0&1&0&1&1\\
0&0&0&1&1&1&1&0&0&0&0&0\\
0&0&0&1&1&0&0&1&1&0&0&0\\
1&0&0&0&0&0&0&0&0&1&1&1\\
0&1&0&1&0&1&0&1&0&1&0&1\\
\end{array}
\right].
\]
By adding the first and third rows to the fifth row, we see that this matrix represents $H_{12}$. Therefore, $\Phi^2_C\npreceq\Psi$ and by Lemma ~\ref{PhiC2}, we may assume that $|C_1|\leq1$.

Since $\Phi_{CX}\npreceq\Psi$, Lemma ~\ref{PhiCD2} implies that $\Psi$ is equivalent to a template with $C=\emptyset$. Hence we will assume from now on that $C=\emptyset$.

Since $C_0=\emptyset$, we have $X_0=\emptyset$. Also, we have seen that $\Lambda|X_1$ is trivial. Therefore, $\Lambda$ is trivial. By performing elementary row operations on every matrix respecting $\Psi$, we may assume that $A_1$ is of the following form, with $Y_0=V_0\cup V_1$, with $Y_1=Y'_1\cup Y''_1$, and with the stars representing arbitrary binary matrices:

\begin{center}
\begin{tabular}{ |c|c|c|c| }
\multicolumn{1}{c}{$Y'_1$}&\multicolumn{1}{c}{$Y''_1$}&\multicolumn{1}{c}{$V_0$}&\multicolumn{1}{c}{$V_1$}\\
\hline
$I$&$*$&0&$Q$\\
\hline
0&0&$I$&$*$\\
\hline
\end{tabular}
.
\end{center}
By Lemma \ref{cecconnected}, we must have $V_0=\emptyset$. Thus, $A_1$ is of the form 
\begin{center}
\begin{tabular}{ |c|c|c| }
\multicolumn{1}{c}{$Y'_1$}&\multicolumn{1}{c}{$Y''_1$}&\multicolumn{1}{c}{$Y_0$}\\
\hline
$I$&$P_1$&$P_0$\\
\hline
\end{tabular}
.
\end{center}
Also, by elementary row operations, we may assume that $\Delta|Y'_1$ is trivial.

We will now show that $|\Delta|\leq2$. Suppose otherwise. Then $\Delta$ contains a subgroup $\Delta'$ isomorphic to $(\mathbb{Z}/2\mathbb{Z})\times(\mathbb{Z}/2\mathbb{Z})$. Repeatedly perform $y$-shifts and operation (11) and then perform operation (2) to obtain the following template:

\[(\emptyset,\emptyset,Y''_1\cup V_1,\emptyset,[\emptyset],\Delta',\{0\}).\]
By Lemma ~\ref{YC}, if $C'=Y''_1\cup V_1$, then every matroid conforming to the following template is a minor of a matroid conforming to $\Psi$:
\[(C',\emptyset,\emptyset,\emptyset,[\emptyset],\Delta',\{0\})\]
The latter template is equivalent to $\Phi^2_C$ since $\Delta'\cong\langle[1,0],[0,1]\rangle$. By contradiction, we deduce that $|\Delta|\leq2$. Therefore, there is at most one nonzero element $\bar{x}\in\Delta$. Let $Y_{1,i}$ consist of the elements $y\in Y''_1$ such that $\bar{x}_y=i$. Similarly, let $Y_{0,i}$ consist of the elements $y\in Y_0$ such that $\bar{x}_y=i$. Thus, $A_1$ is of the following form, where $Y_1=Y'_1\cup Y_{1,0}\cup Y_{1,1}$ and where $Y_0=Y_{0,0}\cup Y_{0,1}$:
\begin{center}
\begin{tabular}{ |c|c|c|c|c| }
\multicolumn{1}{c}{$Y'_1$}&\multicolumn{1}{c}{$Y_{1,0}$}&\multicolumn{1}{c}{$Y_{1,1}$}&\multicolumn{1}{c}{$Y_{0,0}$}&\multicolumn{1}{c}{$Y_{0,1}$}\\
\hline
$I$&$A_{Y_1}$&$B_{Y_1}$&$A_{Y_0}$&$B_{Y_0}$\\
\hline
\end{tabular}
.
\end{center}

By Lemma \ref{B.3}, each row of $B_{Y_1}$ consists either entirely of 0s or entirely of 1s. Any duplicate columns in either $[I|A_{Y_1}]$ or $B_{Y_1}$ produce the same columns in a matrix virtually conforming to $\Psi$. Therefore, we may assume that $|Y_{1,1}|\leq1$, that every column of $A_{Y_1}$ contains at least two nonzero entries, and that no column of $A_{Y_1}$ is a copy of another. Since we are only considering templates to which simple matroids conform, we may assume that no column of $A_{Y_0}$ is a copy of another and also that no column of $B_{Y_0}$ is a copy of another. By Lemma \ref{B.4}, either $Y_{1,0}$ or $Y_{1,1}$ is empty. If $|Y_{1,0}|\geq2$, then $A_{Y_1}$ contains one of the submatrices below, all of which are forbidden by Lemma \ref{B.5}. Therefore, $|Y_{1,0}|\leq1$.
\[
\left[
\begin{array}{cc}
1&1\\
1&1\\
0&1\\
\end{array}
\right],
\left[
\begin{array}{cc}
1&0\\
1&1\\
0&1\\
\end{array}
\right],
\left[
\begin{array}{cc}
1&0\\
1&0\\
0&1\\
0&1\\
\end{array}
\right]
\]

By Lemma \ref{B.6}, if $|Y_{1,0}|=1$, then $Y_{0,1}=\emptyset$. By Lemma \ref{B.7}, if $|Y_{1,1}|=1$, then each column of $A_{Y_{0}}$ contains at most two nonzero entries.

Recall that a  binary matroid $M$ conforms to $\Phi_C$ if there is some binary coextension $N$ of $M$ on ground set $E(M)\cup\{e\}$ such that $N\backslash e$ is graphic. Thus, checking if a binary matroid conforms to $\Phi_C$ amounts to checking if some row can be added to the matrix to make the resulting matroid graphic. There are four cases to check:

Case I: $|Y_{1,0}|=1$

Case II: $|Y_{1,1}|=1$

Case III: $Y_{1,0}=Y_{1,1}=Y_{0,1}=\emptyset$

Case IV: $Y_{1,0}=Y_{1,1}=\emptyset$ and $Y_{0,1}\neq\emptyset$

In the diagrams of matrices below, we will use the  abbreviations \textit{n.p.c.} and \textit{z.p.c.} to stand for ``nonzero entries per column'' and ``zeros per column'' respectively.

\textit{Case I:}
Since $|Y_{1,0}|=1$, the arguments above imply that $Y_{1,1}=Y_{0,1}=\emptyset$. We study the matrix $A_{Y_0}$. Let $X=R\cup S$, where $R$ is the set of rows where $A_{Y_1}$ has its nonzero entries. We will show that $\Psi\preceq\Phi_C$ by appending to every matrix virtually conforming to $\Psi$ a row indexed by $d$ so that the resulting matroid is graphic. The row indexed by $d$ also has a 1 in the column indexed by $Y_{1,0}$, and we will add this row to every row in $R$. After the row indexed by $d$ has been added to every row in $R$, the entire resulting matrix (with rows indexed by $B\cup d$ and columns indexed by $E$) will have at most two nonzero entries per column and will therefore represent a graphic matroid. The form of $A_{Y_0}$ itself (without the row indexed by $d$) is determined by Lemma \ref{old4.6.1}.

First, we study $A_{Y_0}$ when $|R|=2$.
\begin{center}
\begin{tabular}{r|c|c|c| }
\cline{2-4}
$d$&0&0&$1\cdots1$\\
\cline{2-4}
$R$&0&1 n.p.c.&2 n.p.c.\\
\cline{2-4}
$S$&$\leq2$ n.p.c.&$\leq1$ n.p.c&$\leq1$ n.p.c.\\
\cline{2-4}
\end{tabular}
\end{center}

Next, we study $A_{Y_0}$ when $|R|=3$.
\begin{center}
\begin{tabular}{r|c|c|c|c| }
\cline{2-5}
$d$&0&0&0&$1\cdots1$\\
\cline{2-5}
$R$&0&1 n.p.c.&2 n.p.c.&3 n.p.c.\\
\cline{2-5}
$S$&$\leq2$ n.p.c.&$\leq1$ n.p.c&0&$\leq1$ n.p.c\\
\cline{2-5}
\end{tabular}
\end{center}

Next, we study $A_{Y_0}$ when $|R|\geq4$. Here, $J$ denotes a matrix where every entry is a 1.
\begin{center}
\begin{tabular}{r|c|c|c|c|c| }
\cline{2-6}
$d$&0&0&0&$1\cdots1$&$1\cdots1$\\
\cline{2-6}
$R$&0&1 n.p.c.&2 n.p.c.&1 z.p.c.&$J$\\
\cline{2-6}
$S$&$\leq2$ n.p.c.&$\leq1$ n.p.c&0&0&$\leq1$ n.p.c\\
\cline{2-6}
\end{tabular}
\end{center}

\textit{Case II:}
Since $|Y_{1,1}|=1$, the arguments above imply that $|Y_{1,0}|=\emptyset$ and that each column of $A_{Y_0}$ has at most two nonzero entries. We study the matrix $B_{Y_0}$. By Lemma \ref{B.8}, the submatrices below, with the column to the left of the vertical line contained in $B_{Y_1}$, and the column to the right of the vertical line contained in $B_{Y_0}$, are forbidden.
\[
\left[
\begin{array}{c|c}
0&1\\
0&1\\
\end{array}
\right],
\left[
\begin{array}{c|c}
1&0\\
0&1\\
\end{array}
\right],
\left[
\begin{array}{c|c}
1&0\\
1&0\\
\end{array}
\right]
\]
This fact, along with Lemma \ref{old4.6.1}, determines the form of $B_{Y_0}$.

Let $X=R\cup S$, where $R$ is the set of rows where $B_{Y_1}$ has its nonzero entries. We will show that $\Psi\preceq\Phi_C$ by appending to every matrix virtually conforming to $\Psi$ a row indexed by $d$ so that the resulting matroid is graphic. The row indexed by $d$ also has a 1 in the column indexed by $Y_{1,0}$, and we are adding this row to every row in $R$, as well as to every row indexed by an element of $B-X$ where the nonzero element $\bar{x}$ of $\Delta$ is used.

First, we study the case when $R=\emptyset$.
\begin{center}
\begin{tabular}{r|c|c| }
\multicolumn{1}{c}{}&\multicolumn{1}{c}{$Y_{0,1}$}&\multicolumn{1}{c}{$Y_{0,0}$}\\
\cline{2-3}
$d$&$1\cdots1$&$0\cdots0$\\
\cline{2-3}
$S$&$\leq1$ n.p.c.&$A_{Y_0}$ ($\leq2$ n.p.c.)\\
\cline{2-3}
\end{tabular}
\end{center}

Now we study the case when $R\neq\emptyset$. Here $J$ denotes a matrix where every entry is a 1.
\begin{center}
\begin{tabular}{r|c|c|c| }
\multicolumn{1}{c}{}&\multicolumn{2}{c}{$Y_{0,1}$}&\multicolumn{1}{c}{$Y_{0,0}$}\\
\cline{2-4}
$d$&$1\cdots1$&$1\cdots1$&$0\cdots0$\\
\cline{2-4}
$R$&1 z.p.c.&$J$&$A_{Y_0}$\\
\cline{2-3}
$S$&0&$\leq1$ n.p.c.&($\leq2$ n.p.c.)\\
\cline{2-4}
\end{tabular}
\end{center}

\textit{Case III:}
Since $Y_{1,0}=Y_{1,1}=Y_{0,1}=\emptyset$, we have $Y_1=Y'_1$ and $Y_0=Y_{0,0}$. The submatrices below are forbidden from $A_{Y_0}$ because by deleting the rest of $Y_0$ and contracting the elements of $Y_0$ indexing the two given columns, we produce one of the submatrices forbidden by Lemma \ref{B.5}.
\[
Q_1=\left[
\begin{array}{cc}
1&1\\
1&1\\
1&0\\
1&0\\
0&1\\
\end{array}
\right]
Q_2=\left[
\begin{array}{cc}
1&0\\
1&0\\
1&1\\
0&1\\
0&1\\
\end{array}
\right]
Q_3=\left[
\begin{array}{cc}
1&0\\
1&0\\
1&0\\
0&1\\
0&1\\
0&1\\
\end{array}
\right]
Q_4=\left[
\begin{array}{cc}
1&0\\
1&0\\
1&0\\
1&1\\
1&1\\
1&1\\
\end{array}
\right]
\]

If every column of $A_{Y_0}$ has at most two nonzero entries, then $\Psi\sim\Phi_0\preceq\Phi_C$ and can be discarded. Thus, we may assume that there is a column of $A_{Y_0}$ with at least three nonzero entries. Let $H$ be the submatrix of $A_{Y_0}$ consisting of all the columns with at most two nonzero entries.

Let $y$ be an element of $Y_0$ such that $A_{Y_0}[X,\{y\}]$ has a maximum number of nonzero entries among all elements of $Y_0$. Let $X=R\cup S$, where $(A_{Y_0})_{r,y}=1$ for each $r\in R$ and  $(A_{Y_0})_{s,y}=0$ for each $s\in S$. We will prove the following.
\begin{claim}
 \label{vR}Let $v$ be a column of $A_{Y_0}$ such that $v|R$ has at least two zeros. Then either $v$ has at most two nonzero entries, or $v|R$ has exactly two zeros and $v|S$ is a zero vector.
\end{claim}

\begin{subproof}
 If there is a column $v$ of $A_{Y_0}$ such that $v|R$ has exactly two zeros, then if $|R|=3$, the fact that $Q_2$ is forbidden implies that $v$ has at most two nonzero entries. If $|R|>3$, then the fact that $Q_1$ is forbidden implies that $v|S$ is a zero vector. Now, if there is a column $v$ of $A_{Y_0}$ such that $v|R$ has at least three zeros, then since $Q_4$ is forbidden, $v|R$ has at most two nonzero entries. Since $Q_1$, $Q_2$, and $Q_3$ are forbidden (corresponding to when $v|R$ has two, one, or zero nonzero entries), $v$ has at most two nonzero entries.
\end{subproof}

Suppose there are two elements other than $y$ that index columns $v_1$ and $v_2$ of $A_{Y_0}$ with $|R|$ nonzero entries. Since $Q_2$ is forbidden, $v_1|R$ and $v_2|R$ each have at most one zero. Since we are only considering simple matroids conforming to $\Psi$, we have that $v_1|R$ and $v_2|R$ each have exactly one zero.  If $v_1|R$ and $v_2|R$ have their zeros in different rows, then again since $Q_2$ is forbidden, the nonzero entries in $v_1|S$ and $v_2|S$ must be in the same row. Thus, we may divide this case into three subcases:
\begin{enumerate}
\item There is at least one column other than the one indexed by $y$ with $|R|$ nonzero entries, and all such columns have a zero in the same row $r$ of $R$.
\item There are at least two columns other than the one indexed by $y$ with $|R|$ nonzero entries, and all such columns have a nonzero entry in the same row $s$ of $S$.
\item No column other than the one indexed by $y$ has $|R|$ nonzero entries.
\end{enumerate}

In subcase (1), we need to determine the structure of the columns $v$ such that $v|R$ has at least two zeros. In fact, by Claim \ref{vR}, either $v$ is a column of $H$ or $v|R$ has exactly two zeros and $v|S$ is a zero vector. If $v$ is not a column of $H$ but is such that $v|R$ has exactly two zeros, then we must have $|R|\geq5$. Since $Q_1$ is forbidden, $v$ has a zero in the row indexed by $r$. Therefore, $A_{Y_0}[X,Y_0-y]$ has the following form, where $J$ denotes a matrix where every entry is a 1:
\begin{center}
\begin{tabular}{r|c|c|c|}
\cline{2-4}
$R-r$&$J$&1 z.p.c.&\\
\cline{2-3}
$r$&$0\cdots0$&$0\cdots0$&$H$\\
\cline{2-3}
$S$&$\leq1$ n.p.c.&0&\\
\cline{2-4}
\end{tabular}
\end{center}
Below, we append the row $d$ to the matrix, where $d$ also has a 1 in the entry in the column of $y$. By adding $d$ to every row in $R-r$, we see that the resulting matroid is graphic. Thus, in this subcase, $\Psi\preceq\Phi_C$.
\begin{center}
\begin{tabular}{r|c|c|c|}
\cline{2-4}
$d$&$1\cdots1$&$1\cdots1$&$0\cdots0$\\
\cline{2-4}
$R-r$&$J$&1 z.p.c.&\\
\cline{2-3}
$r$&$0\cdots0$&$0\cdots0$&$H$\\
\cline{2-3}
$S$&$\leq1$ n.p.c.&0&\\
\cline{2-4}
\end{tabular}
\end{center}

We now consider subcase (2). Suppose there is some column $v$ of $A_{Y_0}$ such that $v|R$ has two zeros. If $v$ has more than two nonzero entries, then Claim \ref{vR} implies that $v|S$ is a zero vector. By Lemma \ref{B.9}, along with the facts that $Q_1$ is forbidden and that we are only considering templates to which simple matroids conform, no such column $v$ can exist. Therefore, $A_{Y_0}[X,Y_0-y]$ has the following form:
\begin{center}
\begin{tabular}{r|c|c|}
\cline{2-3}
$R$&1 z.p.c.&\\
\cline{2-2}
$s$&$1\cdots1$&$H$\\
\cline{2-2}
$S-s$&0&\\
\cline{2-3}
\end{tabular}
\end{center}
Below, we append the row $d$ to the matrix, where $d$ also has a 1 in the entry in the column of $y$. By adding $d$ to every row in $R\cup \{s\}$, we see that the resulting matroid is graphic. Thus, in this subcase, $\Psi\preceq\Phi_C$.
\begin{center}
\begin{tabular}{r|c|c|c|}
\cline{2-3}
$d$&$1\cdots1$&$0\cdots0$\\
\cline{2-3}
$R$&1 z.p.c.&\\
\cline{2-2}
$s$&$1\cdots1$&$H$\\
\cline{2-2}
$S-s$&0&\\
\cline{2-3}
\end{tabular}
\end{center}

We now consider subcase (3). First, suppose that there is a column $w$ of $A_{Y_0}$ such that $w|R$ has at least two zeros and at least three nonzero entries (so $|R|\geq5$). By Claim \ref{vR}, $w|R$ has exactly two zeros and $w|S$ is a zero vector. Since $Q_2$ is forbidden, every pair of such columns must have zero entries in a common row. By Lemma \ref{new}, all such columns must have a zero in the same row $r$ of $R$. Since $Q_1$ is forbidden, a column $v$ such that $v|$ has exactly one zero must also have its zero in row $r$. Since we are only considering simple matroids, there is at most one such column $v$. Therefore, $A_{Y_0}[X,Y_0-y]$ has the following form, with or without $v$:
\begin{center}
\begin{tabular}{r|c|c|c|}
\multicolumn{1}{c}{}&\multicolumn{1}{c}{$v$}&\multicolumn{2}{c}{}\\
\cline{2-4}
&1&&\multirow{7}{*}{$H$}\\
$R-r$&$\vdots$&1 z.p.c&\\
&1&&\\
\cline{2-3}
$r$&0&$0\cdots0$&\\
\cline{2-3}
&0&&\\
$S$&$\vdots$&0&\\
&0&&\\
\cline{2-4}
\end{tabular}
\end{center}
We append the row $d$ to the matrix, where $d$ also has a 1 in the entry in the column of $y$. By adding $d$ to every row in $R-r$, we see that the resulting matroid is graphic.
\begin{center}
\begin{tabular}{r|c|c|c|}
\multicolumn{1}{c}{}&\multicolumn{1}{c}{$v$}&\multicolumn{2}{c}{}\\
\cline{2-4}
$d$&1&$1\cdots1$&$0\cdots0$\\
\cline{2-4}
&1&&\multirow{7}{*}{$H$}\\
$R-r$&$\vdots$&1 z.p.c&\\
&1&&\\
\cline{2-3}
$r$&0&$0\cdots0$&\\
\cline{2-3}
&0&&\\
$S$&$\vdots$&0&\\
&0&&\\
\cline{2-4}
\end{tabular}
\end{center}

Therefore, we may assume that no column $w$ of $A_{Y_0}$ exists such that $w|R$ has two zeros and at least three nonzero entries. Thus, Claim \ref{vR} implies that $A_{Y_0}[X,Y_0-y]$ is of the following form:
\begin{center}
\begin{tabular}{r|c|c|}
\cline{2-3}
$R$&1 z.p.c.&\multirow{2}{*}{$H$}\\
\cline{2-2}
$S$&0&\\
\cline{2-3}
\end{tabular}
\end{center}
We add the row $d$ to the matrix, where $d$ also has a 1 in the entry in the column of $y$. By adding $d$ to every row in $R$, we see that the resulting matroid is graphic. Thus, in this subcase, $\Psi\preceq\Phi_C$.
\begin{center}
\begin{tabular}{r|c|c|}
\cline{2-3}
$d$&$1\cdots1$&$0\cdots0$\\
\cline{2-3}
$R$&1 z.p.c.&\multirow{2}{*}{$H$}\\
\cline{2-2}
$S$&0&\\
\cline{2-3}
\end{tabular}
\end{center}

\textit{Case IV:}
If any column of $A_{Y_0}$ has three nonzero entries, then by contracting that element of $Y_{0,0}$ we make a column of the identity matrix into a column of $A_{Y_1}$ with two nonzero entries. Since $Y_{0,1}$ is nonempty, this is forbidden by Lemma \ref{B.6}. Therefore, each column of $A_{Y_0}$ contains at most two nonzero entries.

The matrices $Q_5$ and $Q_6$ below are forbidden from $B_{Y_0}$ because by contracting one of the corresponding elements of $Y_0$, we obtain, using a column of the identity matrix, a submatrix forbidden by Lemma \ref{B.8}. The matrix $Q_7$ below is forbidden by Lemma \ref{B.10}.
\[
Q_5=\left[
\begin{array}{cc}
0&1\\
0&1\\
0&1\\
\end{array}
\right]
Q_6=\left[
\begin{array}{cc}
1&0\\
1&0\\
0&1\\
\end{array}
\right]
Q_7=\left[
\begin{array}{cccc}
1&1&0&0\\
1&0&1&0\\
\end{array}
\right]
\]

Let $y$ be an element of $Y_0$ such that $B_{Y_0}[X,\{y\}]$ has a maximum number of nonzero entries among all elements of $Y_{0,1}$. Let $X=R\cup S$, where $(B_{Y_0})_{r,y}=1$ for each $r\in R$ and  $(B_{Y_0})_{s,y}=0$ for each $s\in S$. If $|R|$ is 0 or 1, then we append to each matrix $A$ virtually conforming to $\Psi$ a row which is the characteristic vector of $Y_{0,1}$. If we add this row to each row of $A$ where $\bar{x}$ has been used as the element of $\Delta$, we see that the resulting matroid is graphic. Therefore, we may assume that $|R|\geq2$.

Since $Q_5$ is forbidden, each column of $B_{Y_0}[R,Y_{0,1}]$ contains at most two zeros. Since $Q_6$ is forbidden, every column $w$ such that $w|R$ has two zeros must be such that $w|S$ is a zero vector.

Suppose there are two columns $v_1,v_2$ of $B_{Y_0}$, in addition to the column indexed by $y$, with $|R|$ nonzero entries. Since $Q_6$ is forbidden, $v_1|R$ and $v_2|R$ each have at most one zero. Since we are only considering simple matroids conforming to $\Psi$, $v_1|R$ and $v_2|R$ each have exactly one zero. If $v_1|R$ and $v_2|R$ have their zeros in different rows, then again since $Q_6$ is forbidden, the nonzero entries in $v_1|S$ and $v_2|S$ must be in the same row. Thus, we have the same three subcases as we did in Case III. In each subcase, we will determine the structure of $B_{Y_0}[X,Y_{0,1}-y]$. Since we are only considering simple matroids conforming to $\Psi$, we may assume that no column of $B_{Y_0}$ is a copy of another.

Let us consider subcase (1). If there is a column $v$ of $B_{Y_0}[R,Y_{0,1}]$ with two zeros, then since $Q_6$ is forbidden one of the zeros of $v$ must be in row $r$. Therefore, $B_{Y_0}[X,Y_{0,1}-y]$ is of the following form, where $J$ denotes a matrix where every entry is a 1:
\begin{center}
\begin{tabular}{r|c|c|}
\cline{2-3}
$R-r$&1 z.p.c.&$J$\\
\cline{2-3}
$r$&$0\cdots0$&$0\cdots0$\\
\cline{2-3}
$S$&0&$\leq1$ n.p.c.\\
\cline{2-3}
\end{tabular}
\end{center}
Append to every matrix $A$ conforming to $\Psi$ an additional row that is the characteristic vector of $Y_{0,1}$. By adding this characteristic vector to each row of $A$ where $\bar{x}$ has been used as the element of $\Delta$ as well as to each row of $R-r$, we see that the resulting matroid is graphic. Therefore, $\Psi\preceq\Phi_C$.

Now, we consider subcase (2). Then there are columns $v_1$ and $v_2$ of $B_{Y_0}$, other than the column indexed by $y$, with $|R|$ nonzero entries. Suppose $w$ is a column of $B_{Y_0}$ such that $w|R$ has two zeros. Since $Q_6$ is forbidden, $w$ must have a zero in each of the rows where $v_1$ and $v_2$ have their zeros. But then $B_{Y_0}$ contains $Q_7$. Therefore, no column of $B_{Y_0}[R,Y_{0,1}]$ has two zeros. Therefore, recalling that each column of $B_{Y_0}[R,Y_{0,1}]$ has at most two zeros, we have that $B_{Y_0}[X,Y_{0,1}-y]$ is of the following form:
\begin{center}
\begin{tabular}{r|c|}
\cline{2-2}
$R$&1 z.p.c.\\
\cline{2-2}
$s$&$1\cdots1$\\
\cline{2-2}
$S-s$&0\\
\cline{2-2}
\end{tabular}
\end{center}
Append to every matrix $A$ conforming to $\Psi$ an additional row that is the characteristic vector of $Y_{0,1}$. By adding this characteristic vector to each row of $A$ where $\bar{x}$ has been used as the element of $\Delta$ as well as to each row of $R\cup \{s\}$, we see that the resulting matroid is graphic. Therefore, $\Psi\preceq\Phi_C$.

Now, we consider subcase (3). If $|R|=2$, then since $Q_7$ is forbidden $B_{Y_0}$ (including the column indexed by $y$) is a submatrix, obtained by deleting columns or any rows but the first two, of one of the following matrices :
\[
T_1=\left[
\begin{array}{ccc}
1&1&0\\
1&0&0\\
0&0&0\\
\vdots&\vdots&\vdots\\
0&0&0\\
\end{array}
\right]
T_2=\left[
\begin{array}{ccc}
1&1&0\\
1&0&1\\
0&0&0\\
\vdots&\vdots&\vdots\\
0&0&0\\
\end{array}
\right]
\]
Append to every matrix $A$ conforming to $\Psi$ an additional row that is the characteristic vector of $Y_{0,1}$. If $B_{Y_0}$ is a submatrix of $T_i$, then add this characteristic vector to the first $i$ rows of $A$ as well as to every row of $A$ where $\bar{x}$ has been used as the element of $\Delta$. We see that the resulting matroid is graphic. Therefore, $\Psi\preceq\Phi_C$. Thus, we may assume that $|R|>2$.

Recall that each column of $B_{Y_0}[R,Y_{0,1}]$ has at most two zeros. Suppose there is a column $w$ of $B_{Y_0}$ such that $w|R$ has exactly two zeros. Since $Q_6$ and $Q_7$ are forbidden, all columns of $B_{Y_0}[X,Y_{0,1}-y]$ must have a zero in the same row $r$ of $R$. Since we are only considering simple matroids, there is at most one column $v$ where $v|R$ has exactly one zero. Therefore, $B_{Y_0}[X,Y_{0,1}-y]$ has the following form, with or without $v$:
\begin{center}
\begin{tabular}{r|c|c|}
\multicolumn{1}{c}{}&\multicolumn{1}{c}{$v$}&\multicolumn{1}{c}{}\\
\cline{2-3}
&1&\\
$R-r$&$\vdots$&1 z.p.c\\
&1&\\
\cline{2-3}
$r$&0&$0\cdots0$\\
\cline{2-3}
&0&\\
$S$&$\vdots$&0\\
&0&\\
\cline{2-3}
\end{tabular}
\end{center}
Append to every matrix $A$ conforming to $\Psi$ an additional row that is the characteristic vector of $Y_{0,1}$ and add this characteristic vector to every row of $R-r$ as well as to every row of $A$ where $\bar{x}$ has been used as the element of $\Delta$. We see that the resulting matroid is graphic.

Therefore, we may assume that every column of $B_{Y_0}[R,Y_{0,1}-y]$ has exactly one zero. Thus, $B_{Y_0}[X,Y_{0,1}-y]$ is of the following form:
\begin{center}
\begin{tabular}{r|c|}
\cline{2-2}
$R$&1 z.p.c.\\
\cline{2-2}
$S$&0\\
\cline{2-2}
\end{tabular}
\end{center}
Append to every matrix $A$ conforming to $\Psi$ an additional row that is the characteristic vector of $Y_{0,1}$. By adding this characteristic vector to each row of $A$ where $\bar{x}$ has been used as the element of $\Delta$ as well as to each row of $R$, we see that the resulting matroid is graphic. Therefore, $\Psi\preceq\Phi_C$. This completes the proof.
\end{proof}

\section{Vertical and Cyclic Connectivity}
\label{sec:vertical}
A matroid $M$ is \textit{vertically $k$-connected} if, for each partition $(X,Y)$ of the ground set of $M$ with $r(X)+r(Y)-r(M)<k-1$, either $X$ or $Y$ is spanning. If a matroid is vertically $k$-connected, then its dual is said to be \textit{cyclically $k$-connected}. The next hypothesis is similar to Hypothesis \ref{hyp:connected-template}. It is also found in \cite{gvz-prob} and is a modification of the same hypothesis of Geelen, Gerards, and Whittle \cite{ggw15}. It is also believed to be true, but its proof is also forthcoming in future papers by Geelen, Gerards, and Whittle.

\begin{hypothesis}[{\cite[Hypothesis 4.6, binary case]{gvz-prob}}]
 \label{hyp:cliquetemplate}
Let $\mathcal M$ be a proper minor-closed class of binary matroids. Then there exist $k,n\in\mathbb{Z}_+$ and frame templates $\Phi_1,\ldots,\Phi_s,\Psi_1,\ldots,\Psi_t$ such that
\begin{enumerate}
\item $\mathcal{M}$ contains each of the classes $\mathcal{M}(\Phi_1),\ldots,\mathcal{M}(\Phi_s)$,
\item $\mathcal{M}$ contains the duals of the matroids in each of the classes $\mathcal{M}(\Psi_1),\ldots,\mathcal{M}(\Psi_t)$,
\item if $M$ is a simple vertically $k$-connected member of $\mathcal M$ with an $M(K_n)$-minor, then $M$ is a member of at least one of the classes $\mathcal{M}(\Phi_1),\ldots,\mathcal{M}(\Phi_s)$, and
\item if $M$ is a cosimple cyclically $k$-connected member of $\mathcal M$ with an $M^*(K_n)$-minor, then $M^*$ is a member of at least one of the classes $\mathcal{M}(\Psi_1),\ldots,\mathcal{M}(\Psi_t)$.
\end{enumerate}
\end{hypothesis}

We claim that Theorems \ref{vert-connectedevencycle}-\ref{cyc-connectedevencut} and Corollary \ref{cor:cyc-four-terminals} follow from the proofs of Theorems \ref{connectedevencycle}-\ref{connectedevencut} and Corollary \ref{cor:four-terminals}. Indeed, the only difference is that Lemma \ref{vert-cecconnected} must be used instead of Lemma \ref{cecconnected}.

\begin{theorem}
 \label{vert-connectedevencycle}
Suppose Hypothesis \ref{hyp:cliquetemplate} holds. Then there exists $k,n\in\mathbb{Z}_+$ such that a vertically $k$-connected binary matroid with an $M(K_n)$-minor is contained in $\mathcal{EX}(\mathrm{PG}(3,2)\backslash e, L_{19}, L_{11})$  if and only if it is an even-cycle matroid.
\end{theorem}
\begin{theorem}
 \label{vert-blockingpair}
Suppose Hypothesis \ref{hyp:cliquetemplate} holds. Then there exists $k,n\in\mathbb{Z}_+$ such that a vertically $k$-connected binary matroid with an $M(K_n)$-minor is contained in $\mathcal{EX}(\mathrm{PG}(3,2)\backslash L, M^*(K_6))$ if and only if it is an even-cycle matroid with a blocking pair.
\end{theorem}
\begin{theorem}
 \label{cyc-connectedevencut}
Suppose Hypothesis \ref{hyp:cliquetemplate} holds. Then there exists $k,n\in\mathbb{Z}_+$ such that  a cyclically $k$-connected binary matroid with an $M^*(K_n)$-minor is contained in $\mathcal{EX}(M(K_6), H^*_{12})$ if and only if it is an even-cut matroid.
\end{theorem}
\begin{corollary}
 \label{cor:cyc-four-terminals}
Suppose Hypothesis \ref{hyp:cliquetemplate} holds. Then there exists $k,n\in\mathbb{Z}_+$ such that a cyclically $k$-connected binary matroid with an $M^*(K_n)$-minor is contained in $\mathcal{EX}((\mathrm{PG}(3,2)\backslash L)^*, M(K_6))$ if and only if it has an even-cut representation with at most four terminals.
\end{corollary}

It remains to prove the following lemma.
\begin{lemma}
 \label{vert-cecconnected}
Let $\Phi$ be a template with $C=\emptyset$ and with $\Lambda$ trivial. Then at least one of the following holds:
\begin{enumerate}
\item There exists $k\in\mathbb{Z}_+$ such that no simple vertically $k$-connected matroid with an $M(K_n)$-minor virtually conforms to $\Phi$, or
\item $A_1$ is of the following form, with $Y_1=Y_1'\cup Y_1''$ and each $P_i$ an arbitrary binary matrix:
\begin{center}
\begin{tabular}{ |c|c|c| }
\multicolumn{1}{c}{$Y_1'$}&\multicolumn{1}{c}{$Y_1''$}&\multicolumn{1}{c}{$Y_0$}\\
\hline
$I$&$P_1$&$P_0$\\
\hline
\end{tabular}
\end{center}
\end{enumerate}
\end{lemma}

\begin{proof}
By operation (4), we may assume that $A_1$ is of the following form, with $Y_0=V_0\cup V_1$, with $Y_1=Y_1'\cup Y_1''$ and with each $P_i$ an arbitrary binary matrix:

\begin{center}
\begin{tabular}{ |c|c|c|c| }
\multicolumn{1}{c}{$Y_1'$}&\multicolumn{1}{c}{$Y_1''$}&\multicolumn{1}{c}{$V_0$}&\multicolumn{1}{c}{$V_1$}\\
\hline
$I$&$P_1$&0&$P_0$\\
\hline
0&0&$I$&$P_2$\\
\hline
\end{tabular}
.
\end{center}

Suppose $V_0\neq\emptyset$, and let $M$ be a simple binary matroid virtually conforming to $\Phi$. Set $k>|Y_0|+1$. Thus, we have $\lambda(Y_0)=r(Y_0)+r(E(M)-Y_0)-r(M)<r(Y_0)\leq|Y_0|<k-1$. Since $\Lambda$ is trivial and $Y_1$ does not span $M(A_1)$, it is not possible for $E(M)-Y_0$ to be spanning in $M$. Thus, for $M$ to be vertically $k$-connected, we must have $Y_0$ spanning. Thus, $r(M)=r(Y_0)\leq|Y_0|$. Since $M$ is simple and binary, we have $|E(M)|\leq2^{|Y_0|}-1$. Choose $n$ sufficiently large so that $|E(M(K_n))|=\binom{n}{2}>2^{|Y_0|}-1$. Then (1) holds.

Therefore, we may assume that $V_0=\emptyset$. In this case, (2) holds.
\end{proof}

\section*{Acknowledgements}
The authors thank Irene Pivotto and Gordon Royle for providing a list of known excluded minors for the class of even-cycle matroids. The authors also thank the anonymous referees for carefully reading the manuscript and giving many helpful suggestions to improve it.



\section*{Supplementary material}

# Appendix

Computation 1

```
In [1]:  from sage.matroids.advanced import *
         def is_even_cycle(M, certificate=False):
             """
             Check if matroid M is an even-cycle matroid. If
             certificate=True, also return the matroid N such
             that M = N \ e, and N / e is graphic.
             """
             if not isinstance(M, BinaryMatroid):
                 raise TypeError("This function only works on binary matroid
         s")
             e = newlabel(M.groundset())
             for N in M.linear_extensions(e):
                 if (N / e).is_graphic():
                     if certificate:
                         return (True, N, e)
                     else:
                         return True
             return False


         def is_even_cycle_excluded_minor(M):
             """
             Check whether M is an excluded minor for the
             class of even-cycle matroids.
             """
             if is_even_cycle(M):
                 return False
             for e in M.groundset():
                 if not is_even_cycle(M \ e) or not is_even_cycle(M / e):
                     # not minimal
                     return False
             return True
```

Computation 2

```python
def is_even_cut(M, certificate=False):
    """
    Check if matroid M is an even-cut matroid. If
    certificate=True, also return the matroid N
    such that M = N \ e, and N / e is cographic.
    """
    if not isinstance(M, BinaryMatroid):
        raise TypeError("This function only works on binary matroids")
    e = newlabel(M.groundset())
    for N in M.linear_extensions(e):
        if (N / e).dual().is_graphic():
            if certificate:
                return (True, N, e)
            else:
                return True
    return False

def is_even_cut_excluded_minor(M):
    """
    Check whether M is an excluded minor for
    the class of even-cut matroids.
    """
    if is_even_cut(M):
        return False
    for e in M.groundset():
        if not is_even_cut(M \ e) or not is_even_cut(M / e):
            # not minimal
            return False
    return True
```

Computation 3

In [3]:
```python
A = Matrix(GF(2), [[1,0,0,0,0,0,0,0,0,0,0,0,0,1,1,0,0,0,0],
                   [0,1,0,0,0,0,0,0,0,0,0,0,0,1,0,1,0,0,0],
                   [0,0,1,0,0,0,0,0,0,0,0,0,0,1,0,0,1,0,0],
                   [0,0,0,1,0,0,0,0,0,0,0,0,0,1,0,0,0,1,0],
                   [0,0,0,0,1,0,0,0,0,0,0,0,0,1,0,0,0,0,1],
                   [0,0,0,0,0,1,0,0,0,0,0,0,0,0,1,1,0,0,0],
                   [0,0,0,0,0,0,1,0,0,0,0,0,0,0,1,0,1,0,0],
                   [0,0,0,0,0,0,0,1,0,0,0,0,0,0,1,0,0,1,0],
                   [0,0,0,0,0,0,0,0,1,0,0,0,0,0,1,0,0,0,1],
                   [0,0,0,0,0,0,0,0,0,1,0,0,0,0,0,1,1,0,0],
                   [0,0,0,0,0,0,0,0,0,0,1,0,0,0,0,1,0,1,0],
                   [0,0,0,0,0,0,0,0,0,0,0,1,0,0,0,1,0,0,1],
                   [0,0,0,0,0,0,0,0,0,0,0,0,1,0,0,0,1,1,0]])
L19=Matroid(field=GF(2), matrix=A)
is_even_cycle_excluded_minor(L19)
```

Out[3]: True

Computation 4

```
In [4]:  A = Matrix(GF(2), [[1,0,0,0,0,0,1,0,1,0,1],
                            [0,1,0,0,0,0,1,0,0,1,1],
                            [0,0,1,0,0,0,0,1,1,1,0],
                            [0,0,0,1,0,0,0,1,1,0,1],
                            [0,0,0,0,1,0,0,1,0,1,1],
                            [0,0,0,0,0,1,0,0,1,1,1]])
         L11=Matroid(field=GF(2), matrix=A)
         is_even_cycle_excluded_minor(L11)

Out[4]:  True
```

Computation 5

```
In [5]:  A = Matrix(GF(2), [[1,0,0,0,0,1,1,1,1,0,0,0,0,0,0],
                            [0,1,0,0,0,1,0,0,0,1,1,1,0,0,0],
                            [0,0,1,0,0,0,1,0,0,1,0,0,1,1,0],
                            [0,0,0,1,0,0,0,1,0,0,1,0,1,0,1],
                            [0,0,0,0,1,0,0,0,1,0,0,1,0,1,1]])
         MK6=Matroid(field=GF(2), matrix=A)
         is_even_cut_excluded_minor(MK6)

Out[5]:  True
```

Computation 6

```
In [6]:  A = Matrix(GF(2), [[1,0,1,0,1,0,1,0,1,0,1,0],
                            [0,0,0,0,1,1,1,1,0,0,0,0],
                            [0,0,0,0,1,1,0,0,1,1,0,0],
                            [1,1,0,0,0,0,0,0,0,0,1,1],
                            [0,0,1,1,0,0,1,1,0,0,1,1]])
         H12=Matroid(field=GF(2), matrix=A)
         is_even_cut_excluded_minor(H12.dual())

Out[6]:  True
```

Computation 7

```
def template_test(n, AY1, AY0):
    """
    Generate the universal matroid conforming to the relevant
    template, with a n-vertex complete graph, for a total rank
    of n-1+k.
    """
    if not AY1.nrows() == AY0.nrows():
        raise ValueError, "AY1 and AY0 must have the same number of
rows"
    num_elts = n * (n-1) / 2 + n * AY1.ncols() + AY0.ncols()
    A = Matrix(GF(2), n-1 + AY1.nrows(), num_elts)
    i = 0
    k = AY1.nrows()
    # Identity in front
    for j in range(n-1):
        A[j+k,j] = 1
    i = n-1
    # All columns with two nonzero entries in the frame matrix
    for S in Subsets(range(k,n+k-1),2):
        for j in S:
            A[j,i] = 1
        i = i + 1
    # Columns from Y1 with unit columns and all-zero columns
    # below it
    for l in range(AY1.ncols()):
        for j in range(k):
            A[j,i] = AY1[j,l]
        i = i + 1
        for m in range(n-1):
            A[k+m, i] = 1
            for j in range(k):
                A[j,i] = AY1[j,l]
            i = i + 1
    # Columns from Y0
    for l in range(AY0.ncols()):
        for j in range(k):
            A[j, i] = AY0[j, l]
        i=i+1
    return A

# Below, we give an example of the use of this function.
def print_matrix(A):
    s = ""
    for i in range(A.nrows()):
        s += "["
        for j in range(A.ncols()):
            s += str(A[i,j])
        s += "]\n"
    print s
AY0 = Matrix(GF(2), [[1,1,1,0,1],
                     [0,1,0,0,1],
                     [0,0,0,0,1],
                     [1,0,0,1,0],
                     [1,1,1,1,1],
                     [1,0,0,0,1]])
AY1 = Matrix(GF(2), [[1,0,1,1,0,1,1,1],
                     [0,1,0,0,1,0,0,1],
```

```
[0000001111000011111111100001111111111111101]
[0000000000111100000000111100000001111101001]
[0000001111000011110000111100001111000000001]
[0000000000111100011110000000000001111110010]
[0000011111111000000001111000011111111111111]
[0000000000111100001110000111111110000010001]
[1001100100010001000100010001000100010000000]
[0101010010001000100010001000100010001000000]
[0010110001000100010001000100010001000100000]
```

Computation 8

In [8]:
```
AY0 = Matrix(GF(2), 4, 0)
# Matrix(GF(2), 4, 0) is the empty binary matrix with 4 rows and
# 0 columns.
AY1 = Matrix(GF(2), [[1,0,0,0,1],
                     [0,1,0,0,1],
                     [0,0,1,0,1],
                     [0,0,0,1,1]])
A=template_test(4, AY1, AY0)
M=Matroid(field=GF(2), matrix=A)
#   This tests for a (PG(3,2)\e)-minor.
for S in Subsets(M.groundset(), M.rank() - 4):
    if len((M / S).simplify()) >= 14 and (M / S).rank() == 4:
        print S
        break
```
{17, 12, 7}

Computation 9

In [9]:
```
AY0 = Matrix(GF(2), 4, 0)
AY1 = Matrix(GF(2), [[1,0,0,0,1,0],
                     [0,1,0,0,1,0],
                     [0,0,1,0,0,1],
                     [0,0,0,1,0,1]])
A=template_test(4, AY1, AY0)
M=Matroid(field=GF(2), matrix=A)
#   This tests for a (PG(3,2)\e)-minor.
for S in Subsets(M.groundset(), M.rank() - 4):
    if len((M / S).simplify()) >= 14 and (M / S).rank() == 4:
        print S
        break
```
{25, 12, 7}

Computation 10

```
In [10]: AY0 = Matrix(GF(2), 3, 0)
         AY1 = Matrix(GF(2), [[1,0,0,1,0,1],
                              [0,1,0,1,1,0],
                              [0,0,1,0,1,1]])
         A=template_test(3, AY1, AY0)
         M=Matroid(field=GF(2), matrix=A)
         #   This tests for a (PG(3,2)\e)-minor.
         for S in Subsets(M.groundset(), M.rank() - 4):
             if len((M / S).simplify()) >= 14 and (M / S).rank() == 4:
                 print S
                 break
```

{4}

Computation 11

```
In [11]: AY0 = Matrix(GF(2), 3, 0)
         AY1 = Matrix(GF(2), [[1,0,0,1,1,1],
                              [0,1,0,1,0,1],
                              [0,0,1,0,1,1]])
         A=template_test(3, AY1, AY0)
         M=Matroid(field=GF(2), matrix=A)
         #   This tests for a (PG(3,2)\e)-minor.
         for S in Subsets(M.groundset(), M.rank() - 4):
             if len((M / S).simplify()) >= 14 and (M / S).rank() == 4:
                 print S
                 break
```

{4}

Computation 12

```
In [12]: AY0 = Matrix(GF(2), [[1,1,1],
                              [1,0,1],
                              [1,0,0],
                              [0,1,1],
                              [0,1,0],
                              [0,0,1]])
         AY1 = matrix.identity(6)
         A=template_test(3, AY1, AY0)
         M=Matroid(field=GF(2), matrix=A)
         #   This tests for a (PG(3,2)\e)-minor.
         for S in Subsets(M.groundset(), M.rank() - 4):
             if len((M / S).simplify()) >= 14 and (M / S).rank() == 4:
                 print S
                 break
```

{17, 11, 4, 23}

Computation 13

```
In [13]: AY0 = Matrix(GF(2), [[1,1,0],
                               [1,0,1],
                               [1,0,0],
                               [0,1,1],
                               [0,1,0],
                               [0,0,1]])
         AY1 = matrix.identity(6)
         A=template_test(3, AY1, AY0)
         M=Matroid(field=GF(2), matrix=A)
         #    This tests for a (PG(3,2)\e)-minor.
         for S in Subsets(M.groundset(), M.rank() - 4):
             if len((M / S).simplify()) >= 14 and (M / S).rank() == 4:
                 print S
                 break
```

{17, 11, 4, 23}

Computation 14

```
In [14]: AY0 = Matrix(GF(2), [[1,1,0],
                               [1,0,1],
                               [1,0,1],
                               [0,1,1],
                               [0,1,0]])
         AY1 = matrix.identity(5)
         A=template_test(3, AY1, AY0)
         M=Matroid(field=GF(2), matrix=A)
         #    This tests for a (PG(3,2)\e)-minor.
         for S in Subsets(M.groundset(), M.rank() - 4):
             if len((M / S).simplify()) >= 14 and (M / S).rank() == 4:
                 print S
                 break
```

{19, 11, 7}

Computation 15

```
In [15]: AY0 = Matrix(GF(2), [[1,1,1],
                               [1,0,1],
                               [1,0,1],
                               [0,1,1],
                               [0,1,0]])
         AY1 = matrix.identity(5)
         A=template_test(3, AY1, AY0)
         M=Matroid(field=GF(2), matrix=A)
         #    This tests for a (PG(3,2)\e)-minor.
         for S in Subsets(M.groundset(), M.rank() - 4):
             if len((M / S).simplify()) >= 14 and (M / S).rank() == 4:
                 print S
                 break
```

{19, 11, 7}

Computation 16

```
In [16]: AY0 = Matrix(GF(2), [[1,1,0],
                               [1,0,1],
                               [1,0,1],
                               [0,1,1],
                               [0,1,1]])
         AY1 = matrix.identity(5)
         A=template_test(3, AY1, AY0)
         M=Matroid(field=GF(2), matrix=A)

         print M.has_minor(L11)

         #   This tests for a (PG(3,2)_{-2})-minor.
         for S in Subsets(M.groundset(), M.rank() - 4):
             if len((M / S).simplify()) >= 13 and (M / S).rank() == 4:
                 print S
                 break
```

```
True
{8, 20, 4}
```

Computation 17

```
In [17]: AY0 = Matrix(GF(2), [[1,1,0],
                               [1,1,0],
                               [1,0,1],
                               [1,0,1],
                               [0,1,1]])
         AY1 = matrix.identity(5)
         A=template_test(3, AY1, AY0)
         M=Matroid(field=GF(2), matrix=A)

         print M.has_minor(L11)

         #   This tests for a (PG(3,2)_{-2})-minor.
         for S in Subsets(M.groundset(), M.rank() - 4):
             if len((M / S).simplify()) >= 13 and (M / S).rank() == 4:
                 print S
                 break
```

```
True
{8, 20, 4}
```

Computation 18

```
In [18]:  AY0 = Matrix(GF(2), [[1,1,1],
                               [1,1,0],
                               [1,0,1],
                               [1,0,1],
                               [0,1,1]])
          AY1 = matrix.identity(5)
          A=template_test(3, AY1, AY0)
          M=Matroid(field=GF(2), matrix=A)
          #   This tests for a (PG(3,2)\e)-minor.
          for S in Subsets(M.groundset(), M.rank() - 4):
              if len((M / S).simplify()) >= 14 and (M / S).rank() == 4:
                  print S
                  break
```

{10, 19, 14}


Computation 19

```
In [19]:  AY0 = Matrix(GF(2), [[1,1,1],
                               [1,1,0],
                               [1,0,1],
                               [1,0,1],
                               [0,1,0]])
          AY1 = matrix.identity(5)
          A=template_test(3, AY1, AY0)
          M=Matroid(field=GF(2), matrix=A)
          #   This tests for a (PG(3,2)\e)-minor.
          for S in Subsets(M.groundset(), M.rank() - 4):
              if len((M / S).simplify()) >= 14 and (M / S).rank() == 4:
                  print S
                  break
```

{10, 19, 14}


Computation 20

```
In [20]:  AY0 = Matrix(GF(2), [[1,1,1],
                               [1,1,0],
                               [1,0,1],
                               [1,0,0],
                               [0,1,1]])
          AY1 = matrix.identity(5)
          A=template_test(3, AY1, AY0)
          M=Matroid(field=GF(2), matrix=A)
          #   This tests for a (PG(3,2)\e)-minor.
          for S in Subsets(M.groundset(), M.rank() - 4):
              if len((M / S).simplify()) >= 14 and (M / S).rank() == 4:
                  print S
                  break
```

{17, 18, 4}

Computaion 21

```
In [21]: AY0 = Matrix(GF(2), [[1,1,1],
                              [1,1,0],
                              [1,0,1],
                              [1,0,0],
                              [0,1,0],
                              [0,0,1]])
         AY1 = matrix.identity(6)
         A=template_test(3, AY1, AY0)
         M=Matroid(field=GF(2), matrix=A)
         #  This tests for a (PG(3,2)\e)-minor.
         for S in Subsets(M.groundset(), M.rank() - 4):
             if len((M / S).simplify()) >= 14 and (M / S).rank() == 4:
                 print S
                 break
```

{17, 4, 21, 20}

Computation 22

```
In [22]: AY0 = Matrix(GF(2), [[1,1,1,0],
                              [1,1,0,1],
                              [1,0,1,1],
                              [0,1,1,1]])
         AY1 = matrix.identity(4)
         A=template_test(4, AY1, AY0)
         M=Matroid(field=GF(2), matrix=A)
         #  This tests for a (PG(3,2)\e)-minor.
         for S in Subsets(M.groundset(), M.rank() - 4):
             if len((M / S).simplify()) >= 14 and (M / S).rank() == 4:
                 print S
                 break
```

{17, 12, 7}

Computation 23

```
In [23]: AY0 = Matrix(GF(2), [[1,1,1,1],
                               [1,1,0,1],
                               [1,0,1,1],
                               [0,1,1,1]])
         AY1 = matrix.identity(4)
         A=template_test(4, AY1, AY0)
         M=Matroid(field=GF(2), matrix=A)
         #   This tests for a (PG(3,2)\e)-minor.
         for S in Subsets(M.groundset(), M.rank() - 4):
             if len((M / S).simplify()) >= 14 and (M / S).rank() == 4:
                 print S
                 break
```

{17, 12, 7}


Computation 24

```
In [24]: AY0 = Matrix(GF(2), [[1,1,0],
                               [1,1,0],
                               [1,0,1],
                               [1,0,1],
                               [0,1,1],
                               [0,1,1]])
         AY1 = matrix.identity(6)
         A=template_test(3, AY1, AY0)
         M=Matroid(field=GF(2), matrix=A)

         print M.has_minor(L11)

         #   This tests for a (PG(3,2)_{-2})-minor.
         for S in Subsets(M.groundset(), M.rank() - 4):
             if len((M / S).simplify()) >= 13 and (M / S).rank() == 4:
                 print S
                 break
```

True
{11, 23, 3, 7}


Computation 25

```
In [25]: AY0 = Matrix(GF(2), [[1,1,1],
                               [1,1,0],
                               [1,0,1],
                               [1,0,1],
                               [0,1,1],
                               [0,1,0]])
         AY1 = matrix.identity(6)
         A=template_test(3, AY1, AY0)
         M=Matroid(field=GF(2), matrix=A)
         #  This tests for a (PG(3,2)\e)-minor.
         for S in Subsets(M.groundset(), M.rank() - 4):
             if len((M / S).simplify()) >= 14 and (M / S).rank() == 4:
                 print S
                 break
```

{10, 3, 22, 14}

Computation 26

```
In [26]: AY0 = Matrix(GF(2), [[1,1,1],
                               [1,1,0],
                               [1,0,1],
                               [1,0,0],
                               [0,1,1],
                               [0,1,0]])
         AY1 = matrix.identity(6)
         A=template_test(3, AY1, AY0)
         M=Matroid(field=GF(2), matrix=A)
         #  This tests for a (PG(3,2)\e)-minor.
         for S in Subsets(M.groundset(), M.rank() - 4):
             if len((M / S).simplify()) >= 14 and (M / S).rank() == 4:
                 print S
                 break
```

{12, 11, 4, 22}

Computation 27

```
In [27]: AY0 = Matrix(GF(2), [[0],
                               [1],
                               [1],
                               [1]])
         AY1 = Matrix(GF(2), [[1,0,0,0,1,1],
                              [0,1,0,0,1,0],
                              [0,0,1,0,0,1],
                              [0,0,0,1,0,0]])
         A=template_test(3, AY1, AY0)
         M=Matroid(field=GF(2), matrix=A)
         #   This tests for a (PG(3,2)\e)-minor.
         for S in Subsets(M.groundset(), M.rank() - 4):
             if len((M / S).simplify()) >= 14 and (M / S).rank() == 4:
                 print S
                 break
```

{4, 21}

Computation 28

```
In [28]: AY0 = Matrix(GF(2), [[1],
                               [1],
                               [1],
                               [1]])
         AY1 = Matrix(GF(2), [[1,0,0,0,1,1],
                              [0,1,0,0,1,0],
                              [0,0,1,0,0,1],
                              [0,0,0,1,0,0]])
         A=template_test(3, AY1, AY0)
         M=Matroid(field=GF(2), matrix=A)
         #   This tests for a (PG(3,2)\e)-minor.
         for S in Subsets(M.groundset(), M.rank() - 4):
             if len((M / S).simplify()) >= 14 and (M / S).rank() == 4:
                 print S
                 break
```

{4, 21}

Computation 29

```
AY0 = Matrix(GF(2), [[0],
                     [0],
                     [1],
                     [1],
                     [1]])
AY1 = Matrix(GF(2), [[1,0,0,0,0,1,1],
                     [0,1,0,0,0,1,0],
                     [0,0,1,0,0,0,1],
                     [0,0,0,1,0,0,0],
                     [0,0,0,0,1,0,0]])
A=template_test(3, AY1, AY0)
M=Matroid(field=GF(2), matrix=A)
#   This tests for a (PG(3,2)\e)-minor.
for S in Subsets(M.groundset(), M.rank() - 4):
    if len((M / S).simplify()) >= 14 and (M / S).rank() == 4:
        print S
        break
```

{17, 13, 7}

Computation 30

```
AY0 = Matrix(GF(2), [[1],
                     [0],
                     [1],
                     [1],
                     [1]])
AY1 = Matrix(GF(2), [[1,0,0,0,0,1,1],
                     [0,1,0,0,0,1,0],
                     [0,0,1,0,0,0,1],
                     [0,0,0,1,0,0,0],
                     [0,0,0,0,1,0,0]])
A=template_test(3, AY1, AY0)
M=Matroid(field=GF(2), matrix=A)
#   This tests for a (PG(3,2)\e)-minor.
for S in Subsets(M.groundset(), M.rank() - 4):
    if len((M / S).simplify()) >= 14 and (M / S).rank() == 4:
        print S
        break
```

{17, 13, 7}

Computation 31

```
AY0 = Matrix(GF(2), [[0],
                     [0],
                     [0],
                     [1],
                     [1],
                     [1]])
AY1 = Matrix(GF(2), [[1,0,0,0,0,0,1,1],
                     [0,1,0,0,0,0,1,0],
                     [0,0,1,0,0,0,0,1],
                     [0,0,0,1,0,0,0,0],
                     [0,0,0,0,1,0,0,0],
                     [0,0,0,0,0,1,0,0]])
A=template_test(4, AY1, AY0)
M=Matroid(field=GF(2), matrix=A)
#  This tests for a (PG(3,2)\e)-minor.
for S in Subsets(M.groundset(), M.rank() - 4):
    if len((M / S).simplify()) >= 14 and (M / S).rank() == 4:
        print S
        break
```

{16, 10, 37, 38, 7}

Computation 32

```
AY0 = Matrix(GF(2), [[1],
                     [0],
                     [0],
                     [1],
                     [1],
                     [1]])
AY1 = Matrix(GF(2), [[1,0,0,0,0,0,1,1],
                     [0,1,0,0,0,0,1,0],
                     [0,0,1,0,0,0,0,1],
                     [0,0,0,1,0,0,0,0],
                     [0,0,0,0,1,0,0,0],
                     [0,0,0,0,0,1,0,0]])
A=template_test(4, AY1, AY0)
M=Matroid(field=GF(2), matrix=A)
#  This tests for a (PG(3,2)\e)-minor.
for S in Subsets(M.groundset(), M.rank() - 4):
    if len((M / S).simplify()) >= 14 and (M / S).rank() == 4:
        print S
        break
```

{16, 10, 37, 38, 7}

Computation 33

```
In [33]: AY0 = Matrix(GF(2), [[1],
                               [0],
                               [1],
                               [1]])
         AY1 = Matrix(GF(2), [[1,0,0,0,1,1],
                              [0,1,0,0,1,1],
                              [0,0,1,0,1,0],
                              [0,0,0,1,0,1]])
         A=template_test(3, AY1, AY0)
         M=Matroid(field=GF(2), matrix=A)
         #   This tests for a (PG(3,2)\e)-minor.
         for S in Subsets(M.groundset(), M.rank() - 4):
             if len((M / S).simplify()) >= 14 and (M / S).rank() == 4:
                 print S
                 break
```

{11, 4}

Computation 34

```
In [34]: AY0 = Matrix(GF(2), [[1],
                               [0],
                               [0],
                               [1],
                               [1]])
         AY1 = Matrix(GF(2), [[1,0,0,0,0,1,1],
                              [0,1,0,0,0,1,1],
                              [0,0,1,0,0,1,0],
                              [0,0,0,1,0,0,1],
                              [0,0,0,0,1,0,0]])
         A=template_test(3, AY1, AY0)
         M=Matroid(field=GF(2), matrix=A)
         #   This tests for a (PG(3,2)\e)-minor.
         for S in Subsets(M.groundset(), M.rank() - 4):
             if len((M / S).simplify()) >= 14 and (M / S).rank() == 4:
                 print S
                 break
```

{9, 4, 17}

Computation 35

In [35]: 
```
AY0 = Matrix(GF(2), [[1],
                     [0],
                     [0],
                     [0],
                     [1],
                     [1]])
AY1 = Matrix(GF(2), [[1,0,0,0,0,0,1,1],
                     [0,1,0,0,0,0,1,1],
                     [0,0,1,0,0,0,1,0],
                     [0,0,0,1,0,0,0,1],
                     [0,0,0,0,1,0,0,0],
                     [0,0,0,0,0,1,0,0]])
A=template_test(3, AY1, AY0)
M=Matroid(field=GF(2), matrix=A)
#   This tests for a (PG(3,2)\e)-minor.
for S in Subsets(M.groundset(), M.rank() - 4):
    if len((M / S).simplify()) >= 14 and (M / S).rank() == 4:
        print S
        break
```

{16, 9, 20, 13}

Computation 36

In [36]: 
```
AY0 = Matrix(GF(2), [[1,0],
                     [0,1],
                     [1,1],
                     [1,1]])
AY1 = Matrix(GF(2), [[1,0,0,0,1],
                     [0,1,0,0,1],
                     [0,0,1,0,0],
                     [0,0,0,1,0]])
A=template_test(3, AY1, AY0)
M=Matroid(field=GF(2), matrix=A)
#   This tests for a (PG(3,2)\e)-minor.
for S in Subsets(M.groundset(), M.rank() - 4):
    if len((M / S).simplify()) >= 14 and (M / S).rank() == 4:
        print S
        break
```

{10, 14}

Computation 37

```
AY0 = Matrix(GF(2), [[1,0],
                     [0,1],
                     [1,1],
                     [1,0],
                     [0,1]])
AY1 = Matrix(GF(2), [[1,0,0,0,0,1],
                     [0,1,0,0,0,1],
                     [0,0,1,0,0,0],
                     [0,0,0,1,0,0],
                     [0,0,0,0,1,0]])
A=template_test(3, AY1, AY0)
M=Matroid(field=GF(2), matrix=A)
#   This tests for a (PG(3,2)\e)-minor.
for S in Subsets(M.groundset(), M.rank() - 4):
    if len((M / S).simplify()) >= 14 and (M / S).rank() == 4:
        print S
        break
```

{4, 22, 14}

Computation 38

```
AY0 = Matrix(GF(2), [[1,1],
                     [1,0],
                     [1,1],
                     [1,1]])
AY1 = Matrix(GF(2), [[1,0,0,0,1],
                     [0,1,0,0,1],
                     [0,0,1,0,0],
                     [0,0,0,1,0]])
A=template_test(3, AY1, AY0)
M=Matroid(field=GF(2), matrix=A)
#   This tests for a (PG(3,2)\e)-minor.
for S in Subsets(M.groundset(), M.rank() - 4):
    if len((M / S).simplify()) >= 14 and (M / S).rank() == 4:
        print S
        break
```

{10, 14}

Computation 39

```
AY0 = Matrix(GF(2), [[1,1],
                     [1,0],
                     [1,0],
                     [1,1],
                     [0,1]])
AY1 = Matrix(GF(2), [[1,0,0,0,0,1],
                     [0,1,0,0,0,1],
                     [0,0,1,0,0,0],
                     [0,0,0,1,0,0],
                     [0,0,0,0,1,0]])
A=template_test(3, AY1, AY0)
M=Matroid(field=GF(2), matrix=A)
#  This tests for a (PG(3,2)\e)-minor.
for S in Subsets(M.groundset(), M.rank() - 4):
    if len((M / S).simplify()) >= 14 and (M / S).rank() == 4:
        print S
        break
```

{17, 4, 21}

Computation 40

```
AY0 = Matrix(GF(2), [[1,1],
                     [1,0],
                     [0,1],
                     [1,1]])
AY1 = Matrix(GF(2), [[1,0,0,0,1],
                     [0,1,0,0,1],
                     [0,0,1,0,1],
                     [0,0,0,1,0]])
A=template_test(3, AY1, AY0)
M=Matroid(field=GF(2), matrix=A)
#  This tests for a (PG(3,2)\e)-minor.
for S in Subsets(M.groundset(), M.rank() - 4):
    if len((M / S).simplify()) >= 14 and (M / S).rank() == 4:
        print S
        break
```

{4, 14}

Computation 41

```
In [41]:  AY0 = Matrix(GF(2), [[1,1],
                               [1,0],
                               [0,1],
                               [1,0],
                               [0,1]])
          AY1 = Matrix(GF(2), [[1,0,0,0,0,1],
                               [0,1,0,0,0,1],
                               [0,0,1,0,0,1],
                               [0,0,0,1,0,0],
                               [0,0,0,0,1,0]])
          A=template_test(3, AY1, AY0)
          M=Matroid(field=GF(2), matrix=A)
          #  This tests for a (PG(3,2)\e)-minor.
          for S in Subsets(M.groundset(), M.rank() - 4):
              if len((M / S).simplify()) >= 14 and (M / S).rank() == 4:
                  print S
                  break
```

{17, 4, 14}

Computation 42

```
In [42]:  AY0 = Matrix(GF(2),3,0)
          AY1 = Matrix(GF(2), [[1,0,0,1,1],
                               [0,1,0,1,0],
                               [0,0,1,0,1]])
          A=template_test(4, AY1, AY0)
          M=Matroid(field=GF(2), matrix=A)
          #  This tests for a (PG(3,2)_{-2})-minor.
          for S in Subsets(M.groundset(), M.rank() - 4):
              if len((M / S).simplify()) >= 13 and (M / S).rank() == 4:
                  print S
                  break
```

{0, 12}

Computation 43

```
In [43]:  AY0 = Matrix(GF(2),4,0)
          AY1 = Matrix(GF(2), [[1,0,0,0,1,1],
                               [0,1,0,0,1,1],
                               [0,0,1,0,1,0],
                               [0,0,0,1,0,1]])
          A=template_test(4, AY1, AY0)
          M=Matroid(field=GF(2), matrix=A)
          #  This tests for a (PG(3,2)_{-2})-minor.
          for S in Subsets(M.groundset(), M.rank() - 4):
              if len((M / S).simplify()) >= 13 and (M / S).rank() == 4:
                  print S
                  break
```

{0, 16, 6}

Computation 44

In [44]:
```
AY0 = Matrix(GF(2), [[1],
                     [1],
                     [1],
                     [1]])
AY1 = Matrix(GF(2), [[1,0,0,0,1],
                     [0,1,0,0,1],
                     [0,0,1,0,0],
                     [0,0,0,1,0]])
A=template_test(3, AY1, AY0)
M=Matroid(field=GF(2), matrix=A)
#   This tests for a (PG(3,2)_{-2})-minor.
for S in Subsets(M.groundset(), M.rank() - 4):
    if len((M / S).simplify()) >= 13 and (M / S).rank() == 4:
        print S
        break
```
{18, 4}

Computation 45

In [45]:
```
AY0 = Matrix(GF(2), [[1],
                     [0],
                     [1],
                     [1]])
AY1 = Matrix(GF(2), [[1,0,0,0,1],
                     [0,1,0,0,1],
                     [0,0,1,0,0],
                     [0,0,0,1,0]])
A=template_test(3, AY1, AY0)
M=Matroid(field=GF(2), matrix=A)
#   This tests for a (PG(3,2)_{-2})-minor.
for S in Subsets(M.groundset(), M.rank() - 4):
    if len((M / S).simplify()) >= 13 and (M / S).rank() == 4:
        print S
        break
```
{18, 7}

Computation 46

```
In [46]:  AY0 = Matrix(GF(2), [[1],
                               [1],
                               [0],
                               [1]])
          AY1 = Matrix(GF(2), [[1,0,0,0,1],
                               [0,1,0,0,1],
                               [0,0,1,0,1],
                               [0,0,0,1,0]])
          A=template_test(3, AY1, AY0)
          M=Matroid(field=GF(2), matrix=A)
          #   This tests for a (PG(3,2)_{-2})-minor.
          for S in Subsets(M.groundset(), M.rank() - 4):
              if len((M / S).simplify()) >= 13 and (M / S).rank() == 4:
                  print S
                  break
```

{4, 14}

Computation 47

```
In [47]:  AY0 = Matrix(GF(2), [[1,1,1],
                               [1,0,1],
                               [1,1,0],
                               [0,1,1]])
          AY1 = matrix.identity(4)
          A=template_test(4, AY1, AY0)
          M=Matroid(field=GF(2), matrix=A)
          #   This tests for a (PG(3,2)_{-2})-minor.
          for S in Subsets(M.groundset(), M.rank() - 4):
              if len((M / S).simplify()) >= 13 and (M / S).rank() == 4:
                  print S
                  break
```

{17, 12, 7}

Computation 48

```python
def matrix_from_template(r, Y1I, Y1A, Y1B, Y0A, Y0B, B_Jrows):
    """
    Each of Y1A, Y1B, Y0A, and Y0B either must have the same
    number of rows as Y1I or must be an empty matrix.
    """
    print "Template with C empty, Lambda trivial,Y1 consisting of m
atrices I, Y1A, Y1B, with " , B_Jrows, " all-one rows below B and t
hen"
    print "all-zero rows, Y0 consisting of matrices A, B with ", B_
Jrows, " all-one rows below B and then all-zero rows."
    c = Y1I.nrows()
    A = Matrix(GF(2), r + c, binomial(r+1,2) + (r+1) *Y1I.ncols() +
(r+1) * Y1A.ncols() + (r+1) * Y1B.ncols() + Y0A.ncols() + Y0B.ncols
())
    i = 0
    for j in range(c,r+c):
        A[j,i] = 1
        i = i + 1
    for j in range(c, r-1 + c):
        for k in range(j+1,r+c):
            A[j,i] = 1
            A[k,i] = 1
            i = i + 1
    for j in range(Y1I.ncols()):
        for k in range(c):
            A[k,i] = Y1I[k,j]
        i = i + 1
        for l in range(r):
            for k in range(c):
                A[k,i] = Y1I[k,j]
            A[c+l, i] = 1
            i = i + 1
    for j in range(Y1A.ncols()):
        for k in range(c):
            A[k,i] = Y1A[k,j]
        i = i + 1
        for l in range(r):
            for k in range(c):
                A[k,i] = Y1A[k,j]
            A[c+l, i] = 1
            i = i + 1
    for j in range(Y1B.ncols()):
        for k in range(c):
            A[k,i] = Y1B[k,j]
        for k in range(B_Jrows):
            A[c+k,i] = A[c+k,i] + 1
        i = i + 1
        for l in range(r):
            for k in range(c):
                A[k,i] = Y1B[k,j]
            A[c+l, i] = 1
            for k in range(B_Jrows):
                A[c+k,i] = A[c+k,i] + 1
            i = i + 1
    for j in range(Y0A.ncols()):
        for k in range(c):
            A[k,i] = Y0A[k,j]
```

Template with C empty, Lambda trivial,Y1 consisting of matrices I,
Y1A, Y1B, with  2  all-one rows below B and then
all-zero rows, Y0 consisting of matrices A, B with  2  all-one rows
below B and then all-zero rows.

```
[00000000001111000000000001111111111111111111111111110]
[00000000000000011110000011111111110000000000001010101]
[00000000000000000011111000001111100000111111110011001]
[10001110000100001000010000100001000101111011100001111]
[01001001100010000100001000010000100110111101100001111]
[00100101010001000010000100001000010000100001000000000]
[00010010110001000010000100001000010000100001000000000]
```

In [49]:
```
Y1I = identity_matrix(GF(2),1)
Y1A = Matrix(GF(2), [])
Y1B = Matrix(GF(2), [[1,0]])
Y0A = Matrix(GF(2), [])
Y0B = Matrix(GF(2), [])
A=matrix_from_template(4, Y1I, Y1A, Y1B, Y0A, Y0B, 3)
M = Matroid(A)
M.has_minor(H12)
```

Template with C empty, Lambda trivial,Y1 consisting of matrices I,
Y1A, Y1B, with  3  all-one rows below B and then
all-zero rows, Y0 consisting of matrices A, B with  3  all-one rows
below B and then all-zero rows.

Out[49]:  True

In [50]:
```
Y1I = identity_matrix(GF(2),2)
Y1A = Matrix(GF(2), [[1],
                     [1]])
Y1B = Matrix(GF(2), [[1],
                     [1]])
Y0A = Matrix(GF(2), [])
Y0B = Matrix(GF(2),[])
A=matrix_from_template(3, Y1I, Y1A, Y1B, Y0A, Y0B, 3)
M = Matroid(A)
M.has_minor(H12)
```

Template with C empty, Lambda trivial,Y1 consisting of matrices I,
Y1A, Y1B, with  3  all-one rows below B and then
all-zero rows, Y0 consisting of matrices A, B with  3  all-one rows
below B and then all-zero rows.

Out[50]:  True

```
Y1I = identity_matrix(GF(2),2)
Y1A = Matrix(GF(2), [[1],
                      [1]])
Y1B = Matrix(GF(2), [[0],
                      [0]])
Y0A = Matrix(GF(2), [])
Y0B = Matrix(GF(2), [])
A=matrix_from_template(3, Y1I, Y1A, Y1B, Y0A, Y0B, 3)
M = Matroid(A)
M.has_minor(H12)
```

Template with C empty, Lambda trivial,Y1 consisting of matrices I,
Y1A, Y1B, with  3  all-one rows below B and then
all-zero rows, Y0 consisting of matrices A, B with  3  all-one rows
below B and then all-zero rows.

Out[51]: True

Computation 52

```
Y1I = identity_matrix(GF(2),2)
Y1A = Matrix(GF(2), [[1],
                      [1]])
Y1B = Matrix(GF(2), [[1],
                      [0]])
Y0A = Matrix(GF(2), [])
Y0B = Matrix(GF(2), [])
A=matrix_from_template(3, Y1I, Y1A, Y1B, Y0A, Y0B, 3)
M = Matroid(A)
M.has_minor(H12)
```

Template with C empty, Lambda trivial,Y1 consisting of matrices I,
Y1A, Y1B, with  3  all-one rows below B and then
all-zero rows, Y0 consisting of matrices A, B with  3  all-one rows
below B and then all-zero rows.

Out[52]: True

Computation 53

In [53]:
```
Y1I = identity_matrix(GF(2),3)
Y1A = Matrix(GF(2), [[1,1],
                     [1,1],
                     [0,1]])
Y1B = Matrix(GF(2), [])
Y0A = Matrix(GF(2), [])
Y0B = Matrix(GF(2), [])
A=matrix_from_template(2, Y1I, Y1A, Y1B, Y0A, Y0B, 0)
M = Matroid(A)
M.has_minor(H12)
```

Template with C empty, Lambda trivial,Y1 consisting of matrices I,
Y1A, Y1B, with  0  all-one rows below B and then
all-zero rows, Y0 consisting of matrices A, B with  0  all-one rows
below B and then all-zero rows.

Out[53]: True


Computation 54

In [54]:
```
Y1I = identity_matrix(GF(2),3)
Y1A = Matrix(GF(2), [[1,0],
                     [1,1],
                     [0,1]])
Y1B = Matrix(GF(2), [])
Y0A = Matrix(GF(2), [])
Y0B = Matrix(GF(2), [])
A=matrix_from_template(2, Y1I, Y1A, Y1B, Y0A, Y0B, 0)
M = Matroid(A)
M.has_minor(H12)
```

Template with C empty, Lambda trivial,Y1 consisting of matrices I,
Y1A, Y1B, with  0  all-one rows below B and then
all-zero rows, Y0 consisting of matrices A, B with  0  all-one rows
below B and then all-zero rows.

Out[54]: True


Computation 55

```
Y1I = identity_matrix(GF(2),4)
Y1A = Matrix(GF(2), [[1,0],
                     [1,0],
                     [0,1],
                     [0,1]])
Y1B = Matrix(GF(2), [])
Y0A = Matrix(GF(2), [])
Y0B = Matrix(GF(2), [])
A=matrix_from_template(2, Y1I, Y1A, Y1B, Y0A, Y0B, 0)
M = Matroid(A)
M.has_minor(H12)
```

Template with C empty, Lambda trivial,Y1 consisting of matrices I,
Y1A, Y1B, with  0  all-one rows below B and then
all-zero rows, Y0 consisting of matrices A, B with  0  all-one rows
below B and then all-zero rows.

Out[55]: True

Computation 56

In [56]:
```
Y1I = identity_matrix(GF(2),2)
Y1A = Matrix(GF(2), [[1],
                     [1]])
Y1B = Matrix(GF(2), [])
Y0A = Matrix(GF(2), [])
Y0B = Matrix(GF(2), [[0],
                     [0]])
A=matrix_from_template(3, Y1I, Y1A, Y1B, Y0A, Y0B, 3)
M = Matroid(A)
M.has_minor(H12)
```

Template with C empty, Lambda trivial,Y1 consisting of matrices I,
Y1A, Y1B, with  3  all-one rows below B and then
all-zero rows, Y0 consisting of matrices A, B with  3  all-one rows
below B and then all-zero rows.

Out[56]: True

Computation 57

```
Y1I = identity_matrix(GF(2),2)
Y1A = Matrix(GF(2), [[1],
                     [1]])
Y1B = Matrix(GF(2), [])
Y0A = Matrix(GF(2), [])
Y0B = Matrix(GF(2), [[1],
                     [0]])
A=matrix_from_template(3, Y1I, Y1A, Y1B, Y0A, Y0B, 3)
M = Matroid(A)
M.has_minor(H12)
```

Template with C empty, Lambda trivial,Y1 consisting of matrices I,
Y1A, Y1B, with  3  all-one rows below B and then
all-zero rows, Y0 consisting of matrices A, B with  3  all-one rows
below B and then all-zero rows.

Out[57]: True

Computation 58

In [58]:

```
Y1I = identity_matrix(GF(2),2)
Y1A = Matrix(GF(2), [[1],
                     [1]])
Y1B = Matrix(GF(2), [])
Y0A = Matrix(GF(2), [])
Y0B = Matrix(GF(2), [[1],
                     [1]])
A=matrix_from_template(3, Y1I, Y1A, Y1B, Y0A, Y0B, 3)
M = Matroid(A)
M.has_minor(H12)
```

Template with C empty, Lambda trivial,Y1 consisting of matrices I,
Y1A, Y1B, with  3  all-one rows below B and then
all-zero rows, Y0 consisting of matrices A, B with  3  all-one rows
below B and then all-zero rows.

Out[58]: True

Computation 59

In [59]:
```
Y1I = identity_matrix(GF(2),3)
Y1A = Matrix(GF(2), [])
Y1B = Matrix(GF(2), [[0],
                     [0],
                     [0]])
Y0A = Matrix(GF(2), [[1],
                     [1],
                     [1]])
Y0B = Matrix(GF(2), [])
A=matrix_from_template(3, Y1I, Y1A, Y1B, Y0A, Y0B, 3)
M = Matroid(A)
M.has_minor(H12)
```

Template with C empty, Lambda trivial,Y1 consisting of matrices I,
Y1A, Y1B, with  3  all-one rows below B and then
all-zero rows, Y0 consisting of matrices A, B with  3  all-one rows
below B and then all-zero rows.

Out[59]: True

Computation 60

In [60]:
```
Y1I = identity_matrix(GF(2),3)
Y1A = Matrix(GF(2), [])
Y1B = Matrix(GF(2), [[1],
                     [0],
                     [0]])
Y0A = Matrix(GF(2), [[1],
                     [1],
                     [1]])
Y0B = Matrix(GF(2), [])
A=matrix_from_template(3, Y1I, Y1A, Y1B, Y0A, Y0B, 3)
M = Matroid(A)
M.has_minor(H12)
```

Template with C empty, Lambda trivial,Y1 consisting of matrices I,
Y1A, Y1B, with  3  all-one rows below B and then
all-zero rows, Y0 consisting of matrices A, B with  3  all-one rows
below B and then all-zero rows.

Out[60]: True

Computation 61

```
Y1I = identity_matrix(GF(2),3)
Y1A = Matrix(GF(2), [])
Y1B = Matrix(GF(2), [[1],
                     [1],
                     [0]])
Y0A = Matrix(GF(2), [[1],
                     [1],
                     [1]])
Y0B = Matrix(GF(2), [])
A=matrix_from_template(3, Y1I, Y1A, Y1B, Y0A, Y0B, 3)
M = Matroid(A)
M.has_minor(H12)
```

Template with C empty, Lambda trivial,Y1 consisting of matrices I,
Y1A, Y1B, with  3  all-one rows below B and then
all-zero rows, Y0 consisting of matrices A, B with  3  all-one rows
below B and then all-zero rows.

Out[61]: True

Computation 62

In [62]:

```
Y1I = identity_matrix(GF(2),3)
Y1A = Matrix(GF(2), [])
Y1B = Matrix(GF(2), [[1],
                     [1],
                     [1]])
Y0A = Matrix(GF(2), [[1],
                     [1],
                     [1]])
Y0B = Matrix(GF(2), [])
A=matrix_from_template(2, Y1I, Y1A, Y1B, Y0A, Y0B, 2)
M = Matroid(A)
M.has_minor(H12)
```

Template with C empty, Lambda trivial,Y1 consisting of matrices I,
Y1A, Y1B, with  2  all-one rows below B and then
all-zero rows, Y0 consisting of matrices A, B with  2  all-one rows
below B and then all-zero rows.

Out[62]: True

Computation 63

In [63]:
```
Y1I = identity_matrix(GF(2),2)
Y1A = Matrix(GF(2), [])
Y1B = Matrix(GF(2), [[0],
                     [0]])
Y0A = Matrix(GF(2), [])
Y0B = Matrix(GF(2), [[1],
                     [1]])
A=matrix_from_template(4, Y1I, Y1A, Y1B, Y0A, Y0B, 4)
M = Matroid(A)
M.has_minor(H12)
```

Template with C empty, Lambda trivial,Y1 consisting of matrices I,
Y1A, Y1B, with  4  all-one rows below B and then
all-zero rows, Y0 consisting of matrices A, B with  4  all-one rows
below B and then all-zero rows.

Out[63]: True

Computation 64

In [64]:
```
Y1I = identity_matrix(GF(2),2)
Y1A = Matrix(GF(2), [])
Y1B = Matrix(GF(2), [[1],
                     [0]])
Y0A = Matrix(GF(2), [])
Y0B = Matrix(GF(2), [[0],
                     [1]])
A=matrix_from_template(3, Y1I, Y1A, Y1B, Y0A, Y0B, 2)
M = Matroid(A)
M.has_minor(H12)
```

Template with C empty, Lambda trivial,Y1 consisting of matrices I,
Y1A, Y1B, with  2  all-one rows below B and then
all-zero rows, Y0 consisting of matrices A, B with  2  all-one rows
below B and then all-zero rows.

Out[64]: True

Computation 65

```
Y1I = identity_matrix(GF(2),2)
Y1A = Matrix(GF(2), [])
Y1B = Matrix(GF(2), [[1],
                     [1]])
Y0A = Matrix(GF(2), [])
Y0B = Matrix(GF(2), [[0],
                     [0]])
A=matrix_from_template(3, Y1I, Y1A, Y1B, Y0A, Y0B, 3)
M = Matroid(A)
M.has_minor(H12)
```

Template with C empty, Lambda trivial,Y1 consisting of matrices I,
Y1A, Y1B, with  3  all-one rows below B and then
all-zero rows, Y0 consisting of matrices A, B with  3  all-one rows
below B and then all-zero rows.

Out[65]: True


Computation 66

In [66]:

```
Y1I = identity_matrix(GF(2),6)
Y1A = Matrix(GF(2), [])
Y1B = Matrix(GF(2), [])
Y0A = Matrix(GF(2), [[1,0,1,0],
                     [1,1,0,0],
                     [1,1,1,1],
                     [1,1,1,1],
                     [1,1,1,1],
                     [0,1,1,0]])
Y0B = Matrix(GF(2), [])
A=matrix_from_template(1, Y1I, Y1A, Y1B, Y0A, Y0B, 0)
M = Matroid(A)
M.has_minor(H12)
```

Template with C empty, Lambda trivial,Y1 consisting of matrices I,
Y1A, Y1B, with  0  all-one rows below B and then
all-zero rows, Y0 consisting of matrices A, B with  0  all-one rows
below B and then all-zero rows.

Out[66]: True


Computation 67

In [67]:
```
Y1I = identity_matrix(GF(2),2)
Y1A = Matrix(GF(2), [])
Y1B = Matrix(GF(2), [])
Y0A = Matrix(GF(2), [])
Y0B = Matrix(GF(2),[[1,1,0,0],
                    [1,0,1,0]])
A=matrix_from_template(3, Y1I, Y1A, Y1B, Y0A, Y0B, 3)
M = Matroid(A)
M.has_minor(H12)
```

Template with C empty, Lambda trivial,Y1 consisting of matrices I,
Y1A, Y1B, with  3  all-one rows below B and then
all-zero rows, Y0 consisting of matrices A, B with  3  all-one rows
below B and then all-zero rows.

Out[67]: True


Computation 68

In [68]:
```
Y1I = identity_matrix(GF(2),5)
Y1A = Matrix(GF(2), [])
Y1B = Matrix(GF(2), [])
Y0A = Matrix(GF(2), [[1,1,1,1],
                     [1,1,1,1],
                     [1,1,0,0],
                     [1,0,1,0],
                     [1,0,0,1]])
Y0B = Matrix(GF(2), [])
A=matrix_from_template(2, Y1I, Y1A, Y1B, Y0A, Y0B, 0)
M = Matroid(A)
M.has_minor(H12)
```

Template with C empty, Lambda trivial,Y1 consisting of matrices I,
Y1A, Y1B, with  0  all-one rows below B and then
all-zero rows, Y0 consisting of matrices A, B with  0  all-one rows
below B and then all-zero rows.

Out[68]: True


In [69]:
```
version()
```

Out[69]: 'SageMath version 8.1, Release Date: 2017-12-07'

Mathematics Department, Louisiana State University, Baton Rouge, Louisiana
  *E-mail address*: `kgrace3@lsu.edu`

Mathematics Department, Louisiana State University, Baton Rouge, Louisiana
  *E-mail address*: `svanzwam@math.lsu.edu`


\end{document}